\pdfoutput=1
\documentclass[11pt]{article}
\usepackage[left=1in,right=1in,top=1in,bottom=1in]{geometry}
\usepackage{times}
\usepackage{expl3}
\usepackage{cite}
\usepackage[table]{xcolor}
\usepackage{multirow}
\usepackage{stackengine} 
\usepackage{hhline}
\usepackage{lipsum}
\usepackage{titlesec}
\usepackage{wrapfig}
\usepackage{enumerate}
\usepackage{epsfig}
\usepackage{amsmath}
\usepackage{tabularx}
\usepackage{array}
\usepackage{booktabs}
\usepackage{enumitem}
\usepackage{bbm}
\usepackage{calc}
\usepackage{graphicx}
\usepackage{amsmath}
\usepackage[title]{appendix}
\usepackage{amssymb}
\usepackage{epstopdf}
\usepackage{boldline}
\usepackage{arydshln}
\usepackage{calligra}
\usepackage{bm}
\usepackage{url}
\usepackage{blindtext}
\usepackage{accents}

\newcommand{\define}{\stackrel{\mbox{\tiny def}}{=}}

\newtheorem{theorem}{Theorem}
\newtheorem{proposition}{Proposition}
\newtheorem{corollary}{Corollary}
\newtheorem{lemma}{Lemma}

\usepackage{mathtools}
\usepackage{epstopdf}
\usepackage{balance}
\usepackage{thmtools}
\usepackage{thm-restate}
\usepackage{hyperref}
\usepackage{cleveref}
\usepackage[mathscr]{euscript}

\usepackage[ruled,vlined]{algorithm2e}
\include{pythonlisting}

\newcommand{\ostar}{\mathbin{\mathpalette\make@circled\star}}

\makeatletter
\newcommand{\removelatexerror}{\let\@latex@error\@gobble}
\makeatother
\makeatletter
\newcommand*{\rom}[1]{\expandafter\@slowromancap\romannumeral #1@}
\makeatother

\ExplSyntaxOn
\newcommand\latinabbrev[1]{
  \peek_meaning:NTF. {
    #1\@}%
  { \peek_catcode:NTF a {
      #1.\@ }%
    {#1.\@}}}
\ExplSyntaxOff

\titleclass{\subsubsubsection}{straight}[\subsubsection]

\begin{document}
\vspace{1cm}
\title{Generalized Double Operator Integrals: Finite Dimensions}
\vspace{1.8cm}
\author{Shih-Yu~Chang
\thanks{Shih-Yu Chang is with the Department of Applied Data Science,
San Jose State University, San Jose, CA, U. S. A. (e-mail: {\tt
shihyu.chang@sjsu.edu}). 
           }}

\maketitle

\begin{abstract}
The Double Operator Integral (DOI) framework provides a powerful tool for analyzing perturbations and interactions between self-adjoint operators in functional analysis and spectral theory. However, most existing DOI formulations rely on self-adjointness (Hermitian) or unitary assumptions, limiting their applicability to non-Hermitian settings. Motivated by advancements in non-Hermitian physics and operator theory, this paper introduces Generalized Double Operator Integrals (GDOIs), extending DOI theory to arbitrary non-Hermitian and non-normal matrices. We establish key algebraic properties of GDOIs, derive norm estimations, and develop a perturbation formula that leads to Lipschitz continuity estimates for operator functions. Additionally, we prove the continuity of GDOIs and explore applications in random matrix theory and functional analysis, including tail bounds and H\"older-type estimations. These results provide a unified and flexible integral framework for non-Hermitian spectral analysis, broadening the impact of DOI techniques in non-commutative analysis and mathematical physics. 
\end{abstract}

\begin{keywords}
Double Operator Integral (DOI), spectral mapping theorem, spectral analysis. Lipschitz continuity, non-Hermitian physics, norm estimations. 
\end{keywords}

\section{Introduction}\label{sec: Introduction}

﻿The  Double Operator Integral (DOI)  is a mathematical framework used to examine the interplay between  operators via an integration format. Introduced in functional analysis~\cite{birman2003double,de2002double,peller2016multiple}, DOI theory extends the idea of operator integrals, which generalize classical integral operators by way of incorporating two separate self-adjoint operators into the integration system~\cite{skripka2019multilinear}. Specifically, given self-adjoint operators $\bm{A}$ and $\bm{B}$ on a Hilbert space and a bivariable function $\beta(x, y)$, the DOI takes the form related to spectral measures associated with $\bm{A}$ and $\bm{B}$. This formulation is mainly beneficial when analyzing the perturbation of functions of operators, as it provides a structured approach to characterizing higher-order interactions between non-commutative operators. DOI concept performs a critical function in knowledge how operator capabilities exchange under perturbations, making it a key tool in modern spectral concept and non-commutative geometry and analysis~\cite{connes1994noncommutative}.

﻿The DOI framework has been broadly studied and carried out in numerous areas of mathematical analysis, especially in spectral concept, quantum mechanics, and perturbation theory~\cite{lord2017quantum,coine2022perturbation}. One of its fundamental applications is in  Lifshitz–Krein trace formulation, which describe how spectral shifts occur underneath perturbations of self-adjoint operators~\cite{chattopadhyay2018trace}. DOI techniques are also used in  operator Lipschitz estimates, supporting quantify how operator features reply to perturbations in non-commutative settings. In quantum physics, DOI strategies had been hired to investigate the stability of quantum structures beneath perturbations, in particular in models concerning Schrödinger operators. Moreover, DOI has applications in matrix analysis, random matrix/tensor theory, and sign processing, in which  integral representations of operator functions are crucial for studying transformations and stability properties~\cite{skripka2019stability,chang2022randomMOI,chang2022randomDTI,chang2022randomPDT}. Its capability to address non-commutative structures makes it a precious device in advanced functional analysis, and mathematical physics~\cite{hiai2010matrix}.

Most existing research on  Double Operator Integral (DOI) theory  assumes that the parameter operators are either  self-adjoint (Hermitian matrices) or unitary matrices. This assumption is largely due to the well-defined spectral properties of such operators, which enable the use of orthogonal spectral decompositions. The spectral theorem provides a natural framework for defining operator functions through integration against spectral measures, simplifying the study of perturbations, functional calculus, and operator Lipschitz estimates. In particular, self-adjoint and unitary operators exhibit eigenvector orthogonality, which facilitates explicit computations and stability analysis within DOI formulations. However, many real-world systems, particularly in physics and engineering, do not adhere to these constraints. Inspired by recent advances in  non-Hermitian physics~\cite{ashida2020non}, there is a growing need to extend DOI theory to  non-self-adjoint (non-Hermitian) operators, allowing for a broader class of transformations beyond traditional Hermitian or unitary assumptions.  

Generalizing DOI to arbitrary non-Hermitian parameter operators introduces significant theoretical and computational challenges. Unlike self-adjoint operators, non-Hermitian matrices can have  complex eigenvalues, non-orthogonal eigenvectors, and spectral instabilities, making spectral measure-based approaches less straightforward. These challenges necessitate alternative techniques such as  pseudospectral analysis, non-orthogonal functional calculus, and contour integral representations  to properly define and analyze DOI in the non-Hermitian setting. By removing the restriction of self-adjointness or unitarity, DOI can be applied to emerging fields such as  non-Hermitian quantum mechanics, open quantum systems, dissipative dynamics, and control theory. Expanding DOI to accommodate non-Hermitian operators would not only deepen its mathematical foundations but also bridge the gap between operator theory and modern applications in physics and engineering, where non-normal operators naturally arise. This extension represents a significant step toward a more general and flexible integral framework for operator functions in non-commutative geometry and analysis~\cite{kosaki1984applications}.

This paper develops a generalized framework for  Double Operator Integrals (DOIs)  by extending their applicability beyond the conventional self-adjoint and unitary settings. We first revisit the traditional definition of DOIs and demonstrate that they are special cases of the \emph{Spectral Mapping Theorem}, as established in \cite{chang2024operatorChar}. This insight unifies DOI theory with fundamental results in spectral analysis and provides a more comprehensive theoretical foundation. The discussion in Section~\ref{sec: Conventional DOIs are Special Cases of Spectral Mapping Theorem} highlights this connection and sets the stage for our generalization.  

To extend the DOI framework, we introduce  Generalized Double Operator Integrals (GDOIs), which remove the assumption that parameter operators must be self-adjoint or unitary. This formulation allows the integral framework to be applied to  non-Hermitian and non-normal operators, broadening its relevance to modern mathematical physics and functional analysis. We further explore the algebraic properties of GDOIs, establishing key functional and structural characteristics that distinguish them from conventional DOIs. The theoretical formulation and algebraic properties of GDOIs are presented in Section~\ref{sec: Generalized Double Operator Integrals and Their Algebraic Properties}.  

A crucial aspect of this work is the development of  norm estimations  for GDOIs. In Section~\ref{sec: Norm Estimations}, we derive both upper and lower bound estimates for the norm of GDOIs, providing rigorous analytical tools to quantify their behavior. These estimations are fundamental for assessing the stability and boundedness of operator functions in non-Hermitian settings. In addition to norm analysis, we establish a  perturbation formula  for operator functions of the form \( f(\bm{X}_1) - f(\bm{X}_2) \), using GDOIs as the underlying transformation mechanism. This leads to the derivation of  Lipschitz continuity estimates, which play a key role in understanding the sensitivity of operator functions under perturbations. These results are presented in Section~\ref{sec:Perturbation Formula and Lipschitz Estimation}.  

To ensure the well-posedness of GDOIs, we rigorously prove their  continuity properties  in Section~\ref{sec: Continuity}. This result guarantees that the GDOI operator \( T_{\beta}^{\bm{X}_1,\bm{X}_2}(\bm{Y}) \) behaves smoothly under small variations in the input matrices, reinforcing the robustness of the proposed framework. Finally, we illustrate two important applications of GDOIs in Section~\ref{sec: Applications}. The first application establishes  tail bounds for Lipschitz estimations in random matrix theory, providing probabilistic control over deviations in operator functions. The second extends  Lipschitz estimations to Hölder-type bounds, offering a more general approach to norm control beyond linear constraints.  

By formulating a rigorous GDOI framework, developing norm and perturbation analyses, and demonstrating its applicability in both deterministic and probabilistic settings, this paper significantly extends the scope of DOI theory. These advancements contribute to a deeper understanding of operator functions in  non-Hermitian spectral analysis  and open new directions for research in mathematical physics, functional analysis, and random matrix theory.

\section{Conventional DOIs are Special Cases of Spectral Mapping Theorem}\label{sec: Conventional DOIs are Special Cases of Spectral Mapping Theorem}

Let us review conventional DOI definitions. Given a function $\beta: \mathbb{R}^2 \rightarrow \mathbb{C}$, two Hermitian matrices $\bm{X}_1, \bm{X}_2 \in\mathbb{C}^{n \times n}$, and any matrix $\bm{Y}\in\mathbb{C}^{n \times n}$. From spectral mapping theorem, we have
\begin{eqnarray}\label{eq0:  conv DOI def}
\bm{X}_1&=&\sum\limits_{k_1=1}^{K_1}\sum\limits_{i_1=1}^{\alpha^{\mathrm{G}}_{k_1}}\lambda_{k_1}\bm{P}_{k_1,i_1};\nonumber \\
\bm{X}_2&=&\sum\limits_{k_2=1}^{K_2}\sum\limits_{i_2=1}^{\alpha^{\mathrm{G}}_{k_2}}\lambda_{k_2}\bm{P}_{k_2,i_2};\nonumber \\
\bm{Y}&=&\sum\limits_{k_3=1}^{K_3}\sum\limits_{i_3=1}^{\alpha^{\mathrm{G}}_{k_3}}\lambda_{k_3}\bm{P}_{k_3,i_3}+\sum\limits_{k_3=1}^{K_3}\sum\limits_{i_3=1}^{\alpha^{\mathrm{G}}_{k_3}}\bm{N}_{k_3,i_3},
\end{eqnarray}
where $K_1, K_2, K_3$ are the numbers of distinct eigenvalues of the matrices $\bm{X}_1, \bm{X}_2, \bm{Y}$,  $\alpha^{\mathrm{G}}_{k_1}, \alpha^{\mathrm{G}}_{k_2}, \alpha^{\mathrm{G}}_{k_3}$ are the geometry multiplicities of distinct eigenvalues $\lambda_{k_1}, \lambda_{k_2}, \lambda_{k_3}$ of the matrices $\bm{X}_1, \bm{X}_2, \bm{Y}$, and $\bm{P}_{k_1,i_1}, \bm{P}_{k_2,i_2}, \bm{P}_{k_3,i_3}$ are the projector matrices corresponding to the $i_j$-th geometric component of the $k_j$-th eigenvalue of the matrices $\bm{X}_1 (j=1)$ and $\bm{X}_2 (j=2)$. $\bm{N}_{k_3,i_3}$ is the nilpotent matrix corresponding to the $i_3$-th geometric component of the $k_3$-th eigenvalue of the matrix $\bm{Y}$. Let $\alpha^{\mathrm{A}}_{k_1}, \alpha^{\mathrm{A}}_{k_2}, \alpha^{\mathrm{A}}_{k_3}$ be the algebraic multiplicities of distinct eigenvalues $\lambda_{k_1}, \lambda_{k_2}, \lambda_{k_3}$ of the matrices $\bm{X}_1, \bm{X}_2, \bm{Y}$, we have $\sum\limits_{k_1=1}^{K_1}\alpha^{\mathrm{A}}_{k_1}=n$, $\sum\limits_{k_2=1}^{K_2}\alpha^{\mathrm{A}}_{k_2}=n$, and $\sum\limits_{k_3=1}^{K_3}\alpha^{\mathrm{A}}_{k_3}=n$.

The DOI is a matrix, denoted by $T_{\beta}^{\bm{X}_1,\bm{X}_2}(\bm{Y})$, which can be expressed as~\cite{skripka2019multilinear}:
\begin{eqnarray}\label{eq1:  conv DOI def}
T_{\beta}^{\bm{X}_1,\bm{X}_2}(\bm{Y})\define\sum\limits_{k_1=1}^{K_1}\sum\limits_{k_2=1}^{K_2}\sum\limits_{i_1=1}^{\alpha_{k_1}^{(\mathrm{G})}}\sum\limits_{i_2=1}^{\alpha_{k_2}^{(\mathrm{G})}}\beta(\lambda_{k_1}, \lambda_{k_2})\bm{P}_{k_1,i_1}\bm{Y}\bm{P}_{k_2,i_2}.
\end{eqnarray}
From the decomposition of the matrix $\bm{Y}$ given by Eq.~\eqref{eq0:  conv DOI def}, Eq.~\eqref{eq1:  conv DOI def} can further be expressed as
\begin{eqnarray}\label{eq2:  conv DOI def}
T_{\beta}^{\bm{X}_1,\bm{X}_2}(\bm{Y})&=&\sum\limits_{k_1=1}^{K_1}\sum\limits_{k_2=1}^{K_2}\sum\limits_{i_1=1}^{\alpha_{k_1}^{(\mathrm{G})}}\sum\limits_{i_2=1}^{\alpha_{k_2}^{(\mathrm{G})}}\beta(\lambda_{k_1}, \lambda_{k_2})\bm{P}_{k_1,i_1}\left(\sum\limits_{k_3=1}^{K_3}\sum\limits_{i_3=1}^{\alpha^{\mathrm{G}}_{k_3}}\lambda_{k_3}\bm{P}_{k_3,i_3}+\sum\limits_{k_3=1}^{K_3}\sum\limits_{i_3=1}^{\alpha^{\mathrm{G}}_{k_3}}\bm{N}_{k_3,i_3}\right)\bm{P}_{k_2,i_2}\nonumber \\
&=&\sum\limits_{k_1=1}^{K_1}\sum\limits_{k_2=1}^{K_2}\sum\limits_{k_3=1}^{K_3}\sum\limits_{i_1=1}^{\alpha_{k_1}^{(\mathrm{G})}}\sum\limits_{i_2=1}^{\alpha_{k_2}^{(\mathrm{G})}}\sum\limits_{i_3=1}^{\alpha_{k_3}^{(\mathrm{G})}}\beta(\lambda_{k_1}, \lambda_{k_2})\lambda_{k_3}\bm{P}_{k_1,i_1}\bm{P}_{k_3,i_3}\bm{P}_{k_2,i_2}\nonumber \\
&&+\sum\limits_{k_1=1}^{K_1}\sum\limits_{k_2=1}^{K_2}\sum\limits_{k_3=1}^{K_3}\sum\limits_{i_1=1}^{\alpha_{k_1}^{(\mathrm{G})}}\sum\limits_{i_2=1}^{\alpha_{k_2}^{(\mathrm{G})}}\sum\limits_{i_3=1}^{\alpha_{k_3}^{(\mathrm{G})}}\beta(\lambda_{k_1}, \lambda_{k_2})\bm{P}_{k_1,i_1}\bm{N}_{k_3,i_3}\bm{P}_{k_2,i_2}
\end{eqnarray}

Let us recall Theorem 3 in~\cite{chang2024operatorChar}. Before presenting this theorem, we review several special ntations related to this Theorem 3 in~\cite{chang2024operatorChar}. Given $r$ positive integers $q_1,q_2,\ldots,q_r$, we define $\alpha_{\kappa}(q_1,\ldots,q_r)$ to be the selection of these $r$ arguments $q_1,\ldots,q_r$ to $\kappa$ arguments, i.e., we have
\begin{eqnarray}
\alpha_{\kappa}(q_1,\ldots,q_r)&=&\{q_{\iota_1},q_{\iota_2},\ldots,q_{\iota_\kappa}\}.
\end{eqnarray}
We use $\mbox{Ind}(\alpha_{\kappa}(q_1,\ldots,q_r))$ to obtain indices of those $\kappa$ positive integers $\{q_{\iota_1},q_{\iota_2},\ldots,q_{\iota_\kappa}\}$, i.e., we have
\begin{eqnarray}
\mbox{Ind}(\alpha_{\kappa}(q_1,\ldots,q_r))&=&\{\iota_1,\iota_2,\ldots,\iota_\kappa\}.
\end{eqnarray}
We use $\alpha_{\kappa}(q_1,\ldots,q_r)=1$ to represent $q_{\iota_1}=1,q_{\iota_2}=1,\ldots,q_{\iota_\kappa}=1$. We also use \\
$m_{k_{\mbox{Ind}(\alpha_{\kappa}(q_1,\ldots,q_r))},i_{\mbox{Ind}(\alpha_{\kappa}(q_1,\ldots,q_r))}}-1$ to represent $m_{k_{\iota_1},i_{\iota_1}}-1,m_{k_{\iota_2},i_{\iota_2}}-1,\ldots,m_{k_{\iota_\kappa},i_{\iota_\kappa}}-1$, where $m_{k_{\iota_j},i_{\iota_j}}$ is the order for the nilpotent matrix $\bm{N}_{k_{\iota_j},i_{\iota_j}}$, i.e., $\bm{N}^\ell_{k_{\iota_j},i_{\iota_j}} = \bm{0}$, for $\ell \geq m_{k_{\iota_j},i_{\iota_j}}$ and $j=1,2,\ldots,\kappa$.

Then, Theorem 3 in~\cite{chang2024operatorChar} is given below.
\begin{theorem}\label{thm: Spectral Mapping Theorem for r Variables}
Given an analytic function $f(z_1,z_2,\ldots,z_r)$ within the domain for $|z_l| < R_l$, and the matrix $\bm{X}_l$ with the dimension $m$ and $K_l$ distinct eigenvalues $\lambda_{k_l}$ for $k_l=1,2,\ldots,K_l$ such that
\begin{eqnarray}\label{eq1-1: thm: Spectral Mapping Theorem for r Variables}
\bm{X}_l&=&\sum\limits_{k_l=1}^{K_l}\sum\limits_{i_l=1}^{\alpha_{k_l}^{\mathrm{G}}} \lambda_{k_l} \bm{P}_{k_l,i_l}+
\sum\limits_{k_l=1}^{K_l}\sum\limits_{i_l=1}^{\alpha_{k_l}^{\mathrm{G}}} \bm{N}_{k_l,i_l},
\end{eqnarray}
where $\left\vert\lambda_{k_l}\right\vert<R_l$ for $l=1,2,\ldots,r$.

Then, we have
\begin{eqnarray}\label{eq2: thm: Spectral Mapping Theorem for kappa Variables}
\lefteqn{f(\bm{X}_1,\ldots,\bm{X}_r)=}\nonumber \\
&& \sum\limits_{k_1=\ldots=k_r=1}^{K_1,\ldots,K_r} \sum\limits_{i_1=\ldots=i_r=1}^{\alpha_{k_1}^{(\mathrm{G})},\ldots,\alpha_{k_r}^{(\mathrm{G})}}
f(\lambda_{k_1},\ldots,\lambda_{k_r})\bm{P}_{k_1,i_1}\ldots\bm{P}_{k_r,i_r}\nonumber \\
&&+\sum\limits_{k_1=\ldots=k_r=1}^{K_1,\ldots,K_r} \sum\limits_{i_1=\ldots=i_r=1}^{\alpha_{k_1}^{(\mathrm{G})},\ldots,\alpha_{k_r}^{(\mathrm{G})}}\sum\limits_{\kappa=1}^{r-1}\sum\limits_{\alpha_\kappa(q_1,\ldots,q_r)}\Bigg(\sum\limits_{\alpha_{\kappa}(q_1,\ldots,q_r)=1}^{m_{k_{\mbox{Ind}(\alpha_{\kappa}(q_1,\ldots,q_r))},i_{\mbox{Ind}(\alpha_{\kappa}(q_1,\ldots,q_r))}}-1}\nonumber \\
&&~~~~~ \frac{f^{\alpha_{\kappa}(q_1,\ldots,q_r)}(\lambda_{k_1},\ldots,\lambda_{k_r})}{q_{\iota_1}!q_{\iota_2}!\ldots q_{\iota_\kappa}!}\times \prod\limits_{\substack{\beta =\mbox{Ind}(\alpha_{\kappa}(q_1,\ldots,q_r)), \bm{Y}=\bm{N}^{q_\beta}_{k_\beta,i_\beta} \\ \beta \neq \mbox{Ind}(\alpha_{\kappa}(q_1,\ldots,q_r)), \bm{Y}=\bm{P}_{k_\beta,i_\beta} }
}^{r} \bm{Y}\Bigg) 
\nonumber \\
&&+\sum\limits_{k_1=\ldots=k_r=1}^{K_1,\ldots,K_r} \sum\limits_{i_1=\ldots=i_r=1}^{\alpha_{k_1}^{(\mathrm{G})},\ldots,\alpha_{k_r}^{(\mathrm{G})}} \sum\limits_{q_1=\ldots=q_r=1}^{m_{k_1,i_1}-1,\ldots,m_{k_r,i_r}-1}
\frac{f^{(q_1,\ldots,q_r)}(\lambda_{k_1},\ldots,\lambda_{k_r})}{q_1!\cdots q_r!}\bm{N}^{q_1}_{k_1,i_1}\ldots\bm{N}^{q_r}_{k_r,i_r}
\end{eqnarray}
where we have
\begin{itemize}
\item $\sum\limits_{\alpha_{\kappa}(q_1,\ldots,q_r)}$ is the summation running over all selection of $\alpha_{\kappa}(q_1,\ldots,q_r)$ given $\kappa$;
\item $f^{\alpha_{\kappa}(q_1,\ldots,q_r)}(\lambda_1,\ldots,\lambda_r)$ represents the partial derivatives with respect to variables with indices \\
$\iota_1,\iota_2,\ldots,\iota_\kappa$ and the orders of derivatives given by $q_{\iota_1},q_{\iota_2},\ldots,q_{\iota_\kappa}$. 
\end{itemize}
\end{theorem}

By setting $r=3$ and $f(z_1,z_2,z_3) = \beta(z_1,z_3)z_2$, where $z_1=\lambda_{k_1,i_1},z_2=\lambda_{k_3,i_3}$ and $z_3=\lambda_{k_2,i_2}$, in Theorem~\ref{thm: Spectral Mapping Theorem for r Variables}, Eq.~\eqref{eq2: thm: Spectral Mapping Theorem for kappa Variables} can be reduced to Eq.~\eqref{eq2:  conv DOI def} because 
\begin{eqnarray}
\frac{\partial f(z_1,z_2,z_3)}{\partial z_2}&=&\beta(z_1,z_3);\nonumber \\
\frac{\partial^\ell f(z_1,z_2,z_3)}{\partial^\ell z_2}&=&0,
\end{eqnarray}
where $\ell > 1$. Therefore, we have 
\begin{eqnarray}\label{eq2-1:  conv DOI def}
T_{\beta}^{\bm{X}_1,\bm{X}_2}(\bm{Y})&=&f(\bm{X}_1,\bm{Y},\bm{X}_2).
\end{eqnarray}

If we set $r=3$ and $f(z_1,z_2,z_3) = \beta(z_2,z_3)z_1$, where $z_1=\lambda_{k_3,i_3},z_2=\lambda_{k_1,i_1}$ and $z_3=\lambda_{k_2,i_2}$, in Theorem~\ref{thm: Spectral Mapping Theorem for r Variables}, we can have the first variety of the convention DOI by changing the position for the variable matrix $\bm{Y}$ in DOI as  
\begin{eqnarray}\label{eq3:  conv DOI def}
T_{\beta}^{',\bm{X}_1,\bm{X}_2}(\bm{Y})\define\sum\limits_{k_1=1}^{K_1}\sum\limits_{k_2=1}^{K_2}\sum\limits_{i_1=1}^{\alpha_{k_1}^{(\mathrm{G})}}\sum\limits_{i_2=1}^{\alpha_{k_2}^{(\mathrm{G})}}\beta(\lambda_{k_1}, \lambda_{k_2})\bm{Y}\bm{P}_{k_1,i_1}\bm{P}_{k_2,i_2},
\end{eqnarray}
because 
\begin{eqnarray}
\frac{\partial f(z_1,z_2,z_3)}{\partial z_1}&=&\beta(\lambda_{k_1,i_1},\lambda_{k_2,i_2});\nonumber \\
\frac{\partial^\ell f(z_1,z_2,z_2)}{\partial^\ell z_1}&=&0,
\end{eqnarray}
where $\ell > 1$. Therefore, we have 
\begin{eqnarray}\label{eq3-1:  conv DOI def}
T_{\beta}^{',\bm{X}_1,\bm{X}_2}(\bm{Y})&=&f(\bm{Y},\bm{X}_1,\bm{X}_2).
\end{eqnarray}

Similarly, if we set $r=3$ and $f(z_1,z_2,z_3) = \beta(z_1,z_2)z_3$, where $z_1=\lambda_{k_1,i_1},z_2=\lambda_{k_2,i_2}$ and $z_3=\lambda_{k_3,i_3}$, in Theorem~\ref{thm: Spectral Mapping Theorem for r Variables}, we can have another variety of the conventional DOI by changing the position for the variable matrix $\bm{Y}$ in DOI as   
\begin{eqnarray}\label{eq4:  conv DOI def}
T_{\beta}^{'',\bm{X}_1,\bm{X}_2}(\bm{Y})\define\sum\limits_{k_1=1}^{K_1}\sum\limits_{k_2=1}^{K_2}\sum\limits_{i_1=1}^{\alpha_{k_1}^{(\mathrm{G})}}\sum\limits_{i_2=1}^{\alpha_{k_2}^{(\mathrm{G})}}\beta(\lambda_{k_1}, \lambda_{k_2})\bm{P}_{k_1,i_1}\bm{P}_{k_2,i_2}\bm{Y},
\end{eqnarray}
because 
\begin{eqnarray}
\frac{\partial f(z_1,z_2,z_3)}{\partial z_3}&=&\beta(\lambda_{k_1,i_1},\lambda_{k_2,i_2});\nonumber \\
\frac{\partial^\ell f(z_1,z_2,z_2)}{\partial^\ell z_3}&=&0,
\end{eqnarray}
where $\ell > 1$. Therefore, we have 
\begin{eqnarray}\label{eq4-1:  conv DOI def}
T_{\beta}^{'',\bm{X}_1,\bm{X}_2}(\bm{Y})&=&f(\bm{X}_1,\bm{X}_2,\bm{Y}).
\end{eqnarray}

From Eq.~\eqref{eq2-1: conv DOI def}, Eq.~\eqref{eq3-1: conv DOI def}, and Eq.~\eqref{eq4-1: conv DOI def}, it follows that the application of the multivariable operator spectral mapping theorem, as discussed in~\cite{chang2024operatorChar}, allows for an extension of the conventional DOI definition by considering different arrangement relationships of the input matrix $\bm{Y}$ with \emph{parameter matrices} $\bm{X}_1$ and $\bm{X}_2$.

\section{Generalized Double Operator Integrals and Their Algebraic Properties}\label{sec: Generalized Double Operator Integrals and Their Algebraic Properties}

The generalized double operator integrals (GDOI) will be formulated in Section~\ref{sec: Generalized Double Operator Integrals}. Then, their algebraic properties will be explored in Section~\ref{sec: Algebraic Properties}.

\subsection{Generalized Double Operator Integrals}\label{sec: Generalized Double Operator Integrals}

Given a function $\beta: \mathbb{C}^2 \rightarrow \mathbb{C}$, two any  matrices $\bm{X}_1, \bm{X}_2 \in\mathbb{C}^{n \times n}$, and any matrix $\bm{Y}\in\mathbb{C}^{n \times n}$. From spectral mapping theorem, we have
\begin{eqnarray}\label{eq0:  GDOI def}
\bm{X}_1&=&\sum\limits_{k_1=1}^{K_1}\sum\limits_{i_1=1}^{\alpha^{\mathrm{G}}_{k_1}}\lambda_{k_1}\bm{P}_{k_1,i_1}+\sum\limits_{k_1=1}^{K_1}\sum\limits_{i_1=1}^{\alpha^{\mathrm{G}}_{k_1}}\bm{N}_{k_1,i_1};\nonumber \\
\bm{X}_2&=&\sum\limits_{k_2=1}^{K_2}\sum\limits_{i_2=1}^{\alpha^{\mathrm{G}}_{k_2}}\lambda_{k_2}\bm{P}_{k_2,i_2}+\sum\limits_{k_2=1}^{K_2}\sum\limits_{i_2=1}^{\alpha^{\mathrm{G}}_{k_2}}\bm{N}_{k_2,i_2};\nonumber \\
\bm{Y}&=&\sum\limits_{k_3=1}^{K_3}\sum\limits_{i_3=1}^{\alpha^{\mathrm{G}}_{k_3}}\lambda_{k_3}\bm{P}_{k_3,i_3}+\sum\limits_{k_3=1}^{K_3}\sum\limits_{i_3=1}^{\alpha^{\mathrm{G}}_{k_3}}\bm{N}_{k_3,i_3},
\end{eqnarray}
where $K_1, K_2, K_3$ are the numbers of distinct eigenvalues of the matrices $\bm{X}_1, \bm{X}_2, \bm{Y}$,  $\alpha^{\mathrm{G}}_{k_1}, \alpha^{\mathrm{G}}_{k_2}, \alpha^{\mathrm{G}}_{k_3}$ are the geometry multiplicities of distinct eigenvalues $\lambda_{k_1}, \lambda_{k_2}, \lambda_{k_3}$ of the matrices $\bm{X}_1, \bm{X}_2, \bm{Y}$, and $\bm{P}_{k_1,i_1}, \bm{P}_{k_2,i_2}, \bm{P}_{k_3,i_3}$ are the projector matrices corresponding to the $i_j$-th geometric component of the $k_j$-th eigenvalue of the matrices  $\bm{X}_1 (j=1)$, $\bm{X}_2 (j=2)$ and $\bm{Y} (j=3)$, respectively. $\bm{N}_{k_1,i_1}$, $\bm{N}_{k_2,i_2}$ and $\bm{N}_{k_3,i_3}$ are the nilpotent matrices corresponding to the $i_j$-th geometric component of the $k_j$-th eigenvalue of the matrices  $\bm{X}_1 (j=1)$, $\bm{X}_2 (j=2)$ and $\bm{Y} (j=3)$, respectively. Let $\alpha^{\mathrm{A}}_{k_1}, \alpha^{\mathrm{A}}_{k_2}, \alpha^{\mathrm{A}}_{k_3}$ are the algebraic multiplicities of distinct eigenvalues $\lambda_{k_1}, \lambda_{k_2}, \lambda_{k_3}$ of the matrices $\bm{X}_1, \bm{X}_2, \bm{Y}$, we have $\sum\limits_{k_1=1}^{K_1}\alpha^{\mathrm{A}}_{k_1}=n$, $\sum\limits_{k_2=1}^{K_2}\alpha^{\mathrm{A}}_{k_2}=n$, and $\sum\limits_{k_3=1}^{K_3}\alpha^{\mathrm{A}}_{k_3}=n$.

We adopt Theorem 2 in~\cite{chang2024operatorChar} as below.
\begin{theorem}\label{thm: Spectral Mapping Theorem for Two Variables}
Given an analytic function $f(z_1,z_2)$ within the domain for $|z_1| < R_1$ and $|z_2| < R_2$, the first matrix $\bm{X}_1$ with the dimension $m$ and $K_1$ distinct eigenvalues $\lambda_{k_1}$ for $k_1=1,2,\ldots,K_1$ such that
\begin{eqnarray}\label{eq1-1: thm: Spectral Mapping Theorem for Two Variables}
\bm{X}_1&=&\sum\limits_{k_1=1}^{K_1}\sum\limits_{i_1=1}^{\alpha_{k_1}^{\mathrm{G}}} \lambda_{k_1} \bm{P}_{k_1,i_1}+
\sum\limits_{k_1=1}^{K_1}\sum\limits_{i_1=1}^{\alpha_{k_1}^{\mathrm{G}}} \bm{N}_{k_1,i_1},
\end{eqnarray}
where $\left\vert\lambda_{k_1}\right\vert<R_1$, and second matrix $\bm{X}_2$ with the dimension $m$ and $K_2$ distinct eigenvalues $\lambda_{k_2}$ for $k_2=1,2,\ldots,K_2$ such that
\begin{eqnarray}\label{eq1-2: thm: Spectral Mapping Theorem for Two Variables}
\bm{X}_2&=&\sum\limits_{k_2=1}^{K_2}\sum\limits_{i_2=1}^{\alpha_{k_2}^{\mathrm{G}}} \lambda_{k_2} \bm{P}_{k_2,i_2}+
\sum\limits_{k_2=1}^{K_2}\sum\limits_{i_2=1}^{\alpha_{k_2}^{\mathrm{G}}} \bm{N}_{k_2,i_2},
\end{eqnarray}
where $\left\vert\lambda_{k_2}\right\vert<R_2$.

Then, we have
\begin{eqnarray}\label{eq2: thm: Spectral Mapping Theorem for Two Variables}
f(\bm{X}_1, \bm{X}_2)&=&\sum\limits_{k_1=1}^{K_1}\sum\limits_{k_2=1}^{K_2}\sum\limits_{i_1=1}^{\alpha_{k_1}^{(\mathrm{G})}}\sum\limits_{i_2=1}^{\alpha_{k_2}^{(\mathrm{G})}}f(\lambda_{k_1}, \lambda_{k_2})\bm{P}_{k_1,i_1}\bm{P}_{k_2,i_2} \nonumber \\
&&+\sum\limits_{k_1=1}^{K_1}\sum\limits_{k_2=1}^{K_2}\sum\limits_{i_1=1}^{\alpha_{k_1}^{(\mathrm{G})}}\sum\limits_{i_2=1}^{\alpha_{k_2}^{(\mathrm{G})}}\sum_{q_2=1}^{m_{k_2,i_2}-1}\frac{f^{(-,q_2)}(\lambda_{k_1},\lambda_{k_2})}{q_2!}\bm{P}_{k_1,i_1}\bm{N}_{k_2,i_2}^{q_2} \nonumber \\
&&+\sum\limits_{k_1=1}^{K_1}\sum\limits_{k_2=1}^{K_2}\sum\limits_{i_1=1}^{\alpha_{k_1}^{(\mathrm{G})}}\sum\limits_{i_2=1}^{\alpha_{k_2}^{(\mathrm{G})}}\sum_{q_1=1}^{m_{k_1,i_1}-1}\frac{f^{(q_1,-)}(\lambda_{k_1},\lambda_{k_2})}{q_1!}\bm{N}_{k_1,i_1}^{q_1}\bm{P}_{k_2,i_2}  \nonumber \\
&&+\sum\limits_{k_1=1}^{K_1}\sum\limits_{k_2=1}^{K_2}\sum\limits_{i_1=1}^{\alpha_{k_1}^{(\mathrm{G})}}\sum\limits_{i_2=1}^{\alpha_{k_2}^{(\mathrm{G})}}\sum_{q_1=1}^{m_{k_1,i_1}-1}\sum_{q_2=1}^{m_{k_2,i_2}-1}\frac{f^{(q_1,q_2)}(\lambda_{k_1},\lambda_{k_2})}{q_1!q_2!}\bm{N}_{k_1,i_1}^{q_1}\bm{N}_{k_2,i_2}^{q_2},
\end{eqnarray}
where $m_{k_1,i_1}$ and $m_{k_2,i_2}$ are  orders for the nilpotent matrices $\bm{N}_{k_1,i_1}$ and $\bm{N}_{k_2,i_2}$, respectively, i.e., $\bm{N}^\ell_{k_j,i_j} = \bm{0}$, for $\ell \geq m_{k_j,i_j}$ and $j=1,2$. 
\end{theorem}

From Theorem~\ref{thm: Spectral Mapping Theorem for Two Variables}, the GDOI is a matrix, denoted by $T_{\beta}^{\bm{X}_1,\bm{X}_2}(\bm{Y})$ with any \emph{parameter matrices} $\bm{X}_1$ and $\bm{X}_2$, which is defined by
\begin{eqnarray}\label{eq1:  GDOI def}
T_{\beta}^{\bm{X}_1,\bm{X}_2}(\bm{Y})&\define&\sum\limits_{k_1=1}^{K_1}\sum\limits_{k_2=1}^{K_2}\sum\limits_{i_1=1}^{\alpha_{k_1}^{(\mathrm{G})}}\sum\limits_{i_2=1}^{\alpha_{k_2}^{(\mathrm{G})}}\beta(\lambda_{k_1}, \lambda_{k_2})\bm{P}_{k_1,i_1}\bm{Y}\bm{P}_{k_2,i_2} \nonumber \\
&&+\sum\limits_{k_1=1}^{K_1}\sum\limits_{k_2=1}^{K_2}\sum\limits_{i_1=1}^{\alpha_{k_1}^{(\mathrm{G})}}\sum\limits_{i_2=1}^{\alpha_{k_2}^{(\mathrm{G})}}\sum_{q_2=1}^{m_{k_2,i_2}-1}\frac{\beta^{(-,q_2)}(\lambda_{k_1},\lambda_{k_2})}{q_2!}\bm{P}_{k_1,i_1}\bm{Y}\bm{N}_{k_2,i_2}^{q_2} \nonumber \\
&&+\sum\limits_{k_1=1}^{K_1}\sum\limits_{k_2=1}^{K_2}\sum\limits_{i_1=1}^{\alpha_{k_1}^{(\mathrm{G})}}\sum\limits_{i_2=1}^{\alpha_{k_2}^{(\mathrm{G})}}\sum_{q_1=1}^{m_{k_1,i_1}-1}\frac{\beta^{(q_1,-)}(\lambda_{k_1},\lambda_{k_2})}{q_1!}\bm{N}_{k_1,i_1}^{q_1}\bm{Y}\bm{P}_{k_2,i_2}  \nonumber \\
&&+\sum\limits_{k_1=1}^{K_1}\sum\limits_{k_2=1}^{K_2}\sum\limits_{i_1=1}^{\alpha_{k_1}^{(\mathrm{G})}}\sum\limits_{i_2=1}^{\alpha_{k_2}^{(\mathrm{G})}}\sum_{q_1=1}^{m_{k_1,i_1}-1}\sum_{q_2=1}^{m_{k_2,i_2}-1}\frac{\beta^{(q_1,q_2)}(\lambda_{k_1},\lambda_{k_2})}{q_1!q_2!}\bm{N}_{k_1,i_1}^{q_1}\bm{Y}\bm{N}_{k_2,i_2}^{q_2}.
\end{eqnarray}
From Eq.~\eqref{eq1:  GDOI def}, if matrices $\bm{X}_1$ and $\bm{X}_2$ are Hermitian, the definition of $T_{\beta}^{\bm{X}_1,\bm{X}_2}(\bm{Y})$ given by Eq.~\eqref{eq1:  GDOI def} is reduced to the conventional DOI definition given by Eq.~\eqref{eq1:  conv DOI def}.

\subsection{Algebraic Properties}\label{sec: Algebraic Properties}

In this section, we will establish the algebraic properties of the operator $T_{\beta}^{\bm{X}_1,\bm{X}_2}(\bm{Y})$ defined by Eq.~\eqref{eq1:  GDOI def}. 

If the matrices $\bm{X}_1$ and $\bm{X}_2$ are decomposed as:
\begin{eqnarray}\label{eq:X1 decomp p and n parts}
\bm{X}_1&=&\sum\limits_{k_1=1}^{K_1}\sum\limits_{i_1=1}^{\alpha_{k_1}^{\mathrm{G}}} \lambda_{k_1} \bm{P}_{k_1,i_1}+
\sum\limits_{k_1=1}^{K_1}\sum\limits_{i_1=1}^{\alpha_{k_1}^{\mathrm{G}}} \bm{N}_{k_1,i_1}\nonumber \\
&\define&\bm{X}_{1,\bm{P}}+\bm{X}_{1,\bm{N}},
\end{eqnarray}
and 
\begin{eqnarray}\label{eq:X2 decomp p and n parts}
\bm{X}_2&=&\sum\limits_{k_2=1}^{K_2}\sum\limits_{i_2=1}^{\alpha_{k_2}^{\mathrm{G}}} \lambda_{k_2} \bm{P}_{k_2,i_2}+
\sum\limits_{k_2=1}^{K_2}\sum\limits_{i_2=1}^{\alpha_{k_2}^{\mathrm{G}}} \bm{N}_{k_2,i_2}\nonumber \\
&\define&\bm{X}_{2,\bm{P}}+\bm{X}_{2,\bm{N}},
\end{eqnarray}
then, from the definition of $T_{\beta}^{\bm{X}_1,\bm{X}_2}(\bm{Y})$ given by Eq.~\eqref{eq1: GDOI def}, we have the following decomposition proposition with respect to parameters matrices $\bm{X}_1$ and $\bm{X}_2$ immediately.
\begin{proposition}
Given matrices $\bm{X}_1$ and $\bm{X}_2$, which are decomposed as Eq.~\eqref{eq:X1 decomp p and n parts} and Eq.~\eqref{eq:X2 decomp p and n parts}, respectively, then, we have
\begin{eqnarray}
T_{\beta}^{\bm{X}_1,\bm{X}_2}(\bm{Y})&=&T_{\beta}^{\bm{X}_{1,\bm{P}},\bm{X}_{2,\bm{P}}}(\bm{Y})+T_{\beta}^{\bm{X}_{1,\bm{P}},\bm{X}_{2,\bm{N}}}(\bm{Y})+T_{\beta}^{\bm{X}_{1,\bm{N}},\bm{X}_{2,\bm{P}}}(\bm{Y})+T_{\beta}^{\bm{X}_{1,\bm{N}},\bm{X}_{2,\bm{N}}}(\bm{Y}).
\end{eqnarray}
\end{proposition}

Below, we will focus on algebraic properties with respect to the function $\beta$. We will begin with the following Lemma~\ref{lma: ind of Proj and Nilp} about the independence of product among projection matrices and nilpotent matrices.

\begin{lemma}\label{lma: ind of Proj and Nilp}
Let $\bm{X}_1$ and $\bm{X}_2$ are two matrices with spectral decomposition given by Eq.~\eqref{eq1-1: thm: Spectral Mapping Theorem for Two Variables} and Eq.~\eqref{eq1-2: thm: Spectral Mapping Theorem for Two Variables}, respectively. We have the linear independence of the following four categories of matrices:
\begin{eqnarray}\label{eq1: lma: ind of Proj and Nilp}
1. ~~\bm{P}_{k_1,i_1} \bm{Y} \bm{P}_{k_2,i_2};\nonumber\\ 
2. ~\bm{P}_{k_1,i_1} \bm{Y} \bm{N}_{k_2,i_2}^{q_2};\nonumber\\  
3. ~\bm{N}_{k_1,i_1}^{q_1} \bm{Y} \bm{P}_{k_2,i_2};\nonumber\\ 
4. ~\bm{N}_{k_1,i_1}^{q_1} \bm{Y} \bm{N}_{k_2,i_2}^{q_2},   
\end{eqnarray}
where $\bm{Y} \neq \bm{0}$, $1 \leq q_1 < m_{k_1,i_1}$ and $1 \leq q_2 < m_{k_2,i_2}$.
\end{lemma} 
\textbf{Proof:}
We define the set $\mathcal{S}_{\bm{X}_1, \bm{X}_2, \bm{Y}}$ of matrices as
\begin{eqnarray}\label{eq2: lma: ind of Proj and Nilp}
\mathcal{S}_{\bm{X}_1, \bm{X}_2, \bm{Y}}&\define&\Bigg\{
\sum\limits_{k_1=1}^{K_1}\sum\limits_{k_2=1}^{K_2}\sum\limits_{i_1=1}^{\alpha_{k_1}^{(\mathrm{G})}}\sum\limits_{i_2=1}^{\alpha_{k_2}^{(\mathrm{G})}}c^{PP}_{k_1,k_2,i_1,i_2}\bm{P}_{k_1,i_1}\bm{Y}\bm{P}_{k_2,i_2} 
\nonumber \\
&&+\sum\limits_{k_1=1}^{K_1}\sum\limits_{k_2=1}^{K_2}\sum\limits_{i_1=1}^{\alpha_{k_1}^{(\mathrm{G})}}\sum\limits_{i_2=1}^{\alpha_{k_2}^{(\mathrm{G})}}\sum_{q_2=1}^{m_{k_2,i_2}-1}c^{PN}_{k_1,k_2,i_1,i_2,q_2}\bm{P}_{k_1,i_1}\bm{Y}\bm{N}_{k_2,i_2}^{q_2}\nonumber \\
&&+ \sum\limits_{k_1=1}^{K_1}\sum\limits_{k_2=1}^{K_2}\sum\limits_{i_1=1}^{\alpha_{k_1}^{(\mathrm{G})}}\sum\limits_{i_2=1}^{\alpha_{k_2}^{(\mathrm{G})}}\sum_{q_1=1}^{m_{k_1,i_1}-1}c^{NP}_{k_1,k_2,i_1,i_2,q_1}\bm{N}_{k_1,i_1}^{q_1}\bm{Y}\bm{P}_{k_2,i_2}\nonumber \\
&&+\sum\limits_{k_1=1}^{K_1}\sum\limits_{k_2=1}^{K_2}\sum\limits_{i_1=1}^{\alpha_{k_1}^{(\mathrm{G})}}\sum\limits_{i_2=1}^{\alpha_{k_2}^{(\mathrm{G})}}\sum_{q_1=1}^{m_{k_1,i_1}-1}\sum_{q_2=1}^{m_{k_2,i_2}-1}c^{NN}_{k_1,k_2,i_1,i_2,q_1,q_2}\bm{N}_{k_1,i_1}^{q_1}\bm{Y}\bm{N}_{k_2,i_2}^{q_2}\Bigg\},
\end{eqnarray}
where $c^{PP}_{k_1,k_2,i_1,i_2}$, $c^{PN}_{k_1,k_2,i_1,i_2,q_2}$, $c^{NP}_{k_1,k_2,i_1,i_2,q_1}$ and $c^{NN}_{k_1,k_2,i_1,i_2,q_1,q_2}$ are complex scalers.  

For $j=1,2$, we have
\begin{eqnarray}\label{eq3: lma: ind of Proj and Nilp}
\bm{P}_{k_j,i_j}\bm{P}_{k'_j,i'_j}&=&\bm{P}_{k_j,i_j}\delta(k_j,k'_j)\delta(i_j,i'_j), \nonumber \\
\bm{P}_{k'_j,i'_j}\bm{N}_{k_j,i_j}&=&\bm{N}_{k_j,i_j}\bm{P}_{k'_j,i'_j}=\bm{N}_{k_j,i_j}\delta(k_j,k'_j)\delta(i_j,i'_j), \nonumber \\
\bm{N}_{k_j,i_j}\bm{N}_{k'_j,i'_j}&=&\bm{N}^2_{k_j,i_j}\delta(k_j,k'_j)\delta(i_j,i'_j);
\end{eqnarray}   
therefore, \(\bm{P}_{k_j,i_j}\) and \(\bm{N}_{k_j,i_j}\) act on different generalized eigenspaces, their behaviors are distinct. More specifically:
\begin{itemize}
\item The term \(\bm{P}_{k_1,i_1} \bm{Y} \bm{P}_{k_2,i_2}\) isolates components where both left and right transformations remain in the eigenspaces.
\item The term \(\bm{P}_{k_1,i_1} \bm{Y} \bm{N}_{k_2,i_2}^{q_2}\) involves a right multiplication by a nilpotent matrix, affecting only part of the generalized eigenspace.
\item Similarly, \(\bm{N}_{k_1,i_1}^{q_1} \bm{Y} \bm{P}_{k_2,i_2}\) applies nilpotent transformations on the left.
\item The term \(\bm{N}_{k_1,i_1}^{q_1} \bm{Y} \bm{N}_{k_2,i_2}^{q_2}\) applies nilpotent transformations on both sides.
\end{itemize}
Since nilpotent matrices act non-trivially in their respective Jordan blocks and projectors restrict transformations to specific eigenspaces, these terms span different transformation spaces if $\bm{Y}\neq\bm{0}$. Thus, they cannot be expressed as linear combinations of each other, proving that four categories of matrices given by Eq.~\eqref{eq1: lma: ind of Proj and Nilp} are linearly independent.
$\hfill\Box$

Given a domain $\mathcal{D} \in \mathbb{C}^2$ with $(\lambda_{k_1}, \lambda_{k_2})\in \mathcal{D}$ with a bi-variable analytic function defined over $\mathcal{D}$, namely $\beta(z_1,z_2)$,  we use $\mathcal{S}_{\bm{X}_1,\bm{X}_2,\bm{Y}}(\mathcal{D})$ to represent the following set:
\begin{eqnarray}\label{eq0: lma: GDOI linear homomorphism}
\mathcal{S}_{\bm{X}_1, \bm{X}_2, \bm{Y}}(\mathcal{D})&\define&\Bigg\{
\sum\limits_{k_1=1}^{K_1}\sum\limits_{k_2=1}^{K_2}\sum\limits_{i_1=1}^{\alpha_{k_1}^{(\mathrm{G})}}\sum\limits_{i_2=1}^{\alpha_{k_2}^{(\mathrm{G})}}c^{PP}_{k_1,k_2,i_1,i_2}\bm{P}_{k_1,i_1}\bm{Y}\bm{P}_{k_2,i_2} 
\nonumber \\
&&+\sum\limits_{k_1=1}^{K_1}\sum\limits_{k_2=1}^{K_2}\sum\limits_{i_1=1}^{\alpha_{k_1}^{(\mathrm{G})}}\sum\limits_{i_2=1}^{\alpha_{k_2}^{(\mathrm{G})}}\sum_{q_2=1}^{m_{k_2,i_2}-1}c^{PN}_{k_1,k_2,i_1,i_2,q_2}\bm{P}_{k_1,i_1}\bm{Y}\bm{N}_{k_2,i_2}^{q_2}\nonumber \\
&&+ \sum\limits_{k_1=1}^{K_1}\sum\limits_{k_2=1}^{K_2}\sum\limits_{i_1=1}^{\alpha_{k_1}^{(\mathrm{G})}}\sum\limits_{i_2=1}^{\alpha_{k_2}^{(\mathrm{G})}}\sum_{q_1=1}^{m_{k_1,i_1}-1}c^{NP}_{k_1,k_2,i_1,i_2,q_1}\bm{N}_{k_1,i_1}^{q_1}\bm{Y}\bm{P}_{k_2,i_2}\nonumber \\
&&+\sum\limits_{k_1=1}^{K_1}\sum\limits_{k_2=1}^{K_2}\sum\limits_{i_1=1}^{\alpha_{k_1}^{(\mathrm{G})}}\sum\limits_{i_2=1}^{\alpha_{k_2}^{(\mathrm{G})}}\sum_{q_1=1}^{m_{k_1,i_1}-1}\sum_{q_2=1}^{m_{k_2,i_2}-1}c^{NN}_{k_1,k_2,i_1,i_2,q_1,q_2}\bm{N}_{k_1,i_1}^{q_1}\bm{Y}\bm{N}_{k_2,i_2}^{q_2}\Bigg\},
\end{eqnarray}
such that $c^{PP}_{k_1,k_2,i_1,i_2}$, $c^{PN}_{k_1,k_2,i_1,i_2,q_2}$, $c^{NP}_{k_1,k_2,i_1,i_2,q_1}$ and $c^{NN}_{k_1,k_2,i_1,i_2,q_1,q_2}$ satisfy the following:
 \begin{eqnarray}\label{eq0-1: lma: GDOI linear homomorphism}
c^{PP}_{k_1,k_2,i_1,i_2}&=&\beta(\lambda_{k_1},\lambda_{k_2});\nonumber \\
c^{PN}_{k_1,k_2,i_1,i_2,q_2}&=&\frac{\beta^{(-,q_2)}(\lambda_{k_1},\lambda_{k_2})}{q_2!};\nonumber \\
c^{NP}_{k_1,k_2,i_1,i_2,q_1}&=&\frac{\beta^{(q_1,-)}(\lambda_{k_1},\lambda_{k_2})}{q_1!};\nonumber \\
c^{NN}_{k_1,k_2,i_1,i_2,q_1,q_2}&=&\frac{\beta^{(q_1,q_2)}(\lambda_{k_1},\lambda_{k_2})}{q_1!q_2!}.
\end{eqnarray}

We also define the opeartion $\circ$ between $T^{\bm{X}_1,\bm{X}_2}_\beta$ and $T^{\bm{X}_1,\bm{X}_2}_\gamma$ as 
\begin{eqnarray}
T^{\bm{X}_1,\bm{X}_2}_\beta \circ T^{\bm{X}_1,\bm{X}_2}_\gamma \define T^{\bm{X}_1,\bm{X}_2}_\beta(T^{\bm{X}_1,\bm{X}_2}_\gamma).
\end{eqnarray}

Then, we have the following Lemma~\ref{lma: GDOI linear homomorphism} about linear homomorphism property of the operator $T^{\bm{X}_1,\bm{X}_2}_\beta$. 
\begin{lemma}\label{lma: GDOI linear homomorphism}
Let $\bm{X}_1$ and $\bm{X}_2$ be arbitrary matrices and consider the operator $T^{\bm{X}_1,\bm{X}_2}_\beta$, which is defined by Eq.~\eqref{eq1:  GDOI def}. Given $\beta(z_1,z_2)$ as a bi-variable analytic function over the domain $\mathcal{D}$, the mapping $f: \beta \rightarrow T^{\bm{X}_1,\bm{X}_2}_{\beta}$ is a linear homomorphism. 
\end{lemma}
\textbf{Proof:}
To show that the mapping $f$ is a linear homomorphism, the following properties have to be estblished:
\begin{eqnarray}\label{eq1-1: lma: GDOI linear homomorphism}
T^{\bm{X}_1,\bm{X}_2}_{c_1\beta_1+c_2\beta_2}&=&c_1 T^{\bm{X}_1,\bm{X}_2}_{\beta_1}+ c_2 T^{\bm{X}_1,\bm{X}_2}_{\beta_2},
\end{eqnarray}
where $c_1$ and $c_2$ are two scalers; and 
\begin{eqnarray}\label{eq1-2: lma: GDOI linear homomorphism}
T^{\bm{X}_1,\bm{X}_2}_{\beta\gamma}&=&T^{\bm{X}_1,\bm{X}_2}_{\beta}\circ T^{\bm{X}_1,\bm{X}_2}_{\gamma}.
\end{eqnarray}

\textbf{Proof of Eq.~\eqref{eq1-1: lma: GDOI linear homomorphism}}

Since we have
\begin{eqnarray}\label{eq2: lma: GDOI linear homomorphism}
\lefteqn{T_{c_1\beta_1+c_2\beta_2}^{\bm{X}_1,\bm{X}_2}(\bm{Y})}\nonumber \\
&=_1&\sum\limits_{k_1=1}^{K_1}\sum\limits_{k_2=1}^{K_2}\sum\limits_{i_1=1}^{\alpha_{k_1}^{(\mathrm{G})}}\sum\limits_{i_2=1}^{\alpha_{k_2}^{(\mathrm{G})}}[c_1\beta_1(\lambda_{k_1}, \lambda_{k_2})+c_2\beta_2(\lambda_{k_1}, \lambda_{k_2})]\bm{P}_{k_1,i_1}\bm{Y}\bm{P}_{k_2,i_2} \nonumber \\
&&+\sum\limits_{k_1=1}^{K_1}\sum\limits_{k_2=1}^{K_2}\sum\limits_{i_1=1}^{\alpha_{k_1}^{(\mathrm{G})}}\sum\limits_{i_2=1}^{\alpha_{k_2}^{(\mathrm{G})}}\sum_{q_2=1}^{m_{k_2,i_2}-1}\frac{[c_1\beta_1^{(-,q_2)}(\lambda_{k_1},\lambda_{k_2})+c_2\beta_2^{(-,q_2)}(\lambda_{k_1},\lambda_{k_2})]}{q_2!}\bm{P}_{k_1,i_1}\bm{Y}\bm{N}_{k_2,i_2}^{q_2} \nonumber \\
&&+\sum\limits_{k_1=1}^{K_1}\sum\limits_{k_2=1}^{K_2}\sum\limits_{i_1=1}^{\alpha_{k_1}^{(\mathrm{G})}}\sum\limits_{i_2=1}^{\alpha_{k_2}^{(\mathrm{G})}}\sum_{q_1=1}^{m_{k_1,i_1}-1}\frac{[c_1\beta_1^{(q_1,-)}(\lambda_{k_1},\lambda_{k_2}) + c_2 \beta_2^{(q_1,-)}(\lambda_{k_1},\lambda_{k_2})]}{q_1!}\bm{N}_{k_1,i_1}^{q_1}\bm{Y}\bm{P}_{k_2,i_2}  \nonumber \\
&&+\sum\limits_{k_1=1}^{K_1}\sum\limits_{k_2=1}^{K_2}\sum\limits_{i_1=1}^{\alpha_{k_1}^{(\mathrm{G})}}\sum\limits_{i_2=1}^{\alpha_{k_2}^{(\mathrm{G})}}\sum_{q_1=1}^{m_{k_1,i_1}-1}\sum_{q_2=1}^{m_{k_2,i_2}-1}\frac{[c_1\beta_1^{(q_1,q_2)}(\lambda_{k_1},\lambda_{k_2}) + c_2\beta_2^{(q_1,q_2)}(\lambda_{k_1},\lambda_{k_2})]}{q_1!q_2!}\bm{N}_{k_1,i_1}^{q_1}\bm{Y}\bm{N}_{k_2,i_2}^{q_2}\nonumber \\
&=&c_1 T^{\bm{X}_1,\bm{X}_2}_{\beta_1}(\bm{Y})+ c_2 T^{\bm{X}_1,\bm{X}_2}_{\beta_2}(\bm{Y}),
\end{eqnarray}
where we apply the linearity of the partial derivative in $=_1$. Then, Eq.~\eqref{eq1-1: lma: GDOI linear homomorphism} is established.

\textbf{Proof of Eq.~\eqref{eq1-2: lma: GDOI linear homomorphism}}

The term $T^{\bm{X}_1,\bm{X}_2}_{\beta\gamma}$ can be expressed by
\begin{eqnarray}\label{eq3: lma: GDOI linear homomorphism}
\lefteqn{T_{\beta\gamma}^{\bm{X}_1,\bm{X}_2}(\bm{Y})}\nonumber \\
&=_1&\underbrace{\sum\limits_{k_1=1}^{K_1}\sum\limits_{k_2=1}^{K_2}\sum\limits_{i_1=1}^{\alpha_{k_1}^{(\mathrm{G})}}\sum\limits_{i_2=1}^{\alpha_{k_2}^{(\mathrm{G})}}[\beta(\lambda_{k_1}, \lambda_{k_2})\gamma(\lambda_{k_1}, \lambda_{k_2})]\bm{P}_{k_1,i_1}\bm{Y}\bm{P}_{k_2,i_2}}_{\mbox{Part 1}} \nonumber \\
&&+\underbrace{\sum\limits_{k_1=1}^{K_1}\sum\limits_{k_2=1}^{K_2}\sum\limits_{i_1=1}^{\alpha_{k_1}^{(\mathrm{G})}}\sum\limits_{i_2=1}^{\alpha_{k_2}^{(\mathrm{G})}}\sum_{q_2=1}^{m_{k_2,i_2}-1}\frac{[\beta(\lambda_{k_1}, \lambda_{k_2})\gamma(\lambda_{k_1}, \lambda_{k_2})]^{(-,q_2)}}{q_2!}\bm{P}_{k_1,i_1}\bm{Y}\bm{N}_{k_2,i_2}^{q_2}}_{\mbox{Part 2}} \nonumber \\
&&+\underbrace{\sum\limits_{k_1=1}^{K_1}\sum\limits_{k_2=1}^{K_2}\sum\limits_{i_1=1}^{\alpha_{k_1}^{(\mathrm{G})}}\sum\limits_{i_2=1}^{\alpha_{k_2}^{(\mathrm{G})}}\sum_{q_1=1}^{m_{k_1,i_1}-1}\frac{[\beta(\lambda_{k_1}, \lambda_{k_2})\gamma(\lambda_{k_1}, \lambda_{k_2})]^{(q_1,-)}}{q_1!}\bm{N}_{k_1,i_1}^{q_1}\bm{Y}\bm{P}_{k_2,i_2}}_{\mbox{Part 3}} \nonumber \\
&&+\underbrace{\sum\limits_{k_1=1}^{K_1}\sum\limits_{k_2=1}^{K_2}\sum\limits_{i_1=1}^{\alpha_{k_1}^{(\mathrm{G})}}\sum\limits_{i_2=1}^{\alpha_{k_2}^{(\mathrm{G})}}\sum_{q_1=1}^{m_{k_1,i_1}-1}\sum_{q_2=1}^{m_{k_2,i_2}-1}\frac{[\beta(\lambda_{k_1}, \lambda_{k_2})\gamma(\lambda_{k_1}, \lambda_{k_2})]^{(q_1,q_2)}}{q_1!q_2!}\bm{N}_{k_1,i_1}^{q_1}\bm{Y}\bm{N}_{k_2,i_2}^{q_2}}_{\mbox{Part 4}}.
\end{eqnarray}

On the other hand, we can express $T^{\bm{X}_1,\bm{X}_2}_{\beta}\circ T^{\bm{X}_1,\bm{X}_2}_{\gamma}$ as
\begin{eqnarray}\label{eq3: lma: GDOI linear homomorphism}
 \lefteqn{T^{\bm{X}_1,\bm{X}_2}_{\beta}\circ T^{\bm{X}_1,\bm{X}_2}_{\gamma}=T^{\bm{X}_1,\bm{X}_2}_{\beta}( T^{\bm{X}_1,\bm{X}_2}_{\gamma})}\nonumber \\
&=&\sum\limits_{k_1=1}^{K_1}\sum\limits_{k_2=1}^{K_2}\sum\limits_{i_1=1}^{\alpha_{k_1}^{(\mathrm{G})}}\sum\limits_{i_2=1}^{\alpha_{k_2}^{(\mathrm{G})}}\beta(\lambda_{k_1}, \lambda_{k_2})\bm{P}_{k_1,i_1}[T^{\bm{X}_1,\bm{X}_2}_{\gamma}(\bm{Y})]\bm{P}_{k_2,i_2} \nonumber \\
&&+\sum\limits_{k_1=1}^{K_1}\sum\limits_{k_2=1}^{K_2}\sum\limits_{i_1=1}^{\alpha_{k_1}^{(\mathrm{G})}}\sum\limits_{i_2=1}^{\alpha_{k_2}^{(\mathrm{G})}}\sum_{q_2=1}^{m_{k_2,i_2}-1}\frac{\beta^{(-,q_2)}(\lambda_{k_1},\lambda_{k_2})}{q_2!}\bm{P}_{k_1,i_1}[T^{\bm{X}_1,\bm{X}_2}_{\gamma}(\bm{Y})]\bm{N}_{k_2,i_2}^{q_2} \nonumber \\
&&+\sum\limits_{k_1=1}^{K_1}\sum\limits_{k_2=1}^{K_2}\sum\limits_{i_1=1}^{\alpha_{k_1}^{(\mathrm{G})}}\sum\limits_{i_2=1}^{\alpha_{k_2}^{(\mathrm{G})}}\sum_{q_1=1}^{m_{k_1,i_1}-1}\frac{\beta^{(q_1,-)}(\lambda_{k_1},\lambda_{k_2})}{q_1!}\bm{N}_{k_1,i_1}^{q_1}[T^{\bm{X}_1,\bm{X}_2}_{\gamma}(\bm{Y})]\bm{P}_{k_2,i_2}  \nonumber \\
&&+\sum\limits_{k_1=1}^{K_1}\sum\limits_{k_2=1}^{K_2}\sum\limits_{i_1=1}^{\alpha_{k_1}^{(\mathrm{G})}}\sum\limits_{i_2=1}^{\alpha_{k_2}^{(\mathrm{G})}}\sum_{q_1=1}^{m_{k_1,i_1}-1}\sum_{q_2=1}^{m_{k_2,i_2}-1}\frac{\beta^{(q_1,q_2)}(\lambda_{k_1},\lambda_{k_2})}{q_1!q_2!}\bm{N}_{k_1,i_1}^{q_1}[T^{\bm{X}_1,\bm{X}_2}_{\gamma}(\bm{Y})]\bm{N}_{k_2,i_2}^{q_2}.
\end{eqnarray}
 
We arrange the expansion of R.H.S. of Eq.~\eqref{eq3: lma: GDOI linear homomorphism} according to $\bm{P}_{k_1,i_1}\bm{Y}\bm{P}_{k_2,i_2}$, $\bm{P}_{k_1,i_1}\bm{Y}\bm{N}_{k_2,i_2}^{q_2}$, $\bm{N}_{k_1,i_1}^{q_1}\bm{Y}\bm{P}_{k_2,i_2}$ and $\bm{N}_{k_1,i_1}^{q_1}\bm{Y}\bm{N}_{k_2,i_2}^{q_2}$ to obtain the following four expressions :
\begin{eqnarray}\label{eq4-1: lma: GDOI linear homomorphism}
\sum\limits_{k_1=1}^{K_1}\sum\limits_{k_2=1}^{K_2}\sum\limits_{i_1=1}^{\alpha_{k_1}^{(\mathrm{G})}}\sum\limits_{i_2=1}^{\alpha_{k_2}^{(\mathrm{G})}}[\beta(\lambda_{k_1}, \lambda_{k_2})\gamma(\lambda_{k_1}, \lambda_{k_2})]\bm{P}_{k_1,i_1}\bm{Y}\bm{P}_{k_2,i_2}, 
\end{eqnarray}
\begin{eqnarray}\label{eq4-2: lma: GDOI linear homomorphism}
\sum\limits_{k_1=1}^{K_1}\sum\limits_{k_2=1}^{K_2}\sum\limits_{i_1=1}^{\alpha_{k_1}^{(\mathrm{G})}}\sum\limits_{i_2=1}^{\alpha_{k_2}^{(\mathrm{G})}}\sum_{q_2=1}^{m_{k_2,i_2}-1}\frac{\beta(\lambda_{k_1},\lambda_{k_2})\gamma^{(-,q_2)}(\lambda_{k_1},\lambda_{k_2})}{q_2!}\bm{P}_{k_1,i_1}\bm{Y}\bm{N}_{k_2,i_2}^{q_2}\nonumber \\
+\sum\limits_{k_1=1}^{K_1}\sum\limits_{k_2=1}^{K_2}\sum\limits_{i_1=1}^{\alpha_{k_1}^{(\mathrm{G})}}\sum\limits_{i_2=1}^{\alpha_{k_2}^{(\mathrm{G})}}\sum_{q_2=1}^{m_{k_2,i_2}-1}\frac{\beta^{(-,q_2)}(\lambda_{k_1},\lambda_{k_2})\gamma(\lambda_{k_1},\lambda_{k_2})}{q_2!}\bm{P}_{k_1,i_1}\bm{Y}\bm{N}_{k_2,i_2}^{q_2}\nonumber \\
+\sum\limits_{k_1=1}^{K_1}\sum\limits_{k_2=1}^{K_2}\sum\limits_{i_1=1}^{\alpha_{k_1}^{(\mathrm{G})}}\sum\limits_{i_2=1}^{\alpha_{k_2}^{(\mathrm{G})}}\sum_{q'_2+ q''_2=1}^{q'_2+q''_2 = m_{k_2,i_2}-1}\frac{\beta^{(-,q'_2)}(\lambda_{k_1},\lambda_{k_2})}{q'_2!}\frac{\gamma^{(-,q''_2)}(\lambda_{k_1},\lambda_{k_2})}{q''_2!}\bm{P}_{k_1,i_1}\bm{Y}\bm{N}_{k_2,i_2}^{q'_2+q''_2},
\end{eqnarray}
\begin{eqnarray}\label{eq4-3: lma: GDOI linear homomorphism}
\sum\limits_{k_1=1}^{K_1}\sum\limits_{k_2=1}^{K_2}\sum\limits_{i_1=1}^{\alpha_{k_1}^{(\mathrm{G})}}\sum\limits_{i_2=1}^{\alpha_{k_2}^{(\mathrm{G})}}\sum_{q_1=1}^{m_{k_1,i_1}-1}\frac{\beta(\lambda_{k_1},\lambda_{k_2})\gamma^{(q_1,-)}(\lambda_{k_1},\lambda_{k_2})}{q_1!}\bm{N}^{q_1}_{k_1,i_1}\bm{Y}\bm{P}_{k_2,i_2}\nonumber \\
+\sum\limits_{k_1=1}^{K_1}\sum\limits_{k_2=1}^{K_2}\sum\limits_{i_1=1}^{\alpha_{k_1}^{(\mathrm{G})}}\sum\limits_{i_2=1}^{\alpha_{k_2}^{(\mathrm{G})}}\sum_{q_1=1}^{m_{k_1,i_1}-1}\frac{\beta^{(q_1,-)}(\lambda_{k_1},\lambda_{k_2})\gamma(\lambda_{k_1},\lambda_{k_2})}{q_1!}\bm{N}^{q_1}_{k_1,i_1}\bm{Y}\bm{P}_{k_2,i_2}\nonumber \\
+\sum\limits_{k_1=1}^{K_1}\sum\limits_{k_2=1}^{K_2}\sum\limits_{i_1=1}^{\alpha_{k_1}^{(\mathrm{G})}}\sum\limits_{i_2=1}^{\alpha_{k_2}^{(\mathrm{G})}}\sum_{q'_1+ q''_1=1}^{q'_1+q''_1 = m_{k_1,i_1}-1}\frac{\beta^{(q'_1,-)}(\lambda_{k_1},\lambda_{k_2})}{q'_1!}\frac{\gamma^{(q''_1,-)}(\lambda_{k_1},\lambda_{k_2})}{q''_1!}\bm{N}^{q'_1+q''_1}_{k_1,i_1}\bm{Y}\bm{P}_{k_2,i_2},    
\end{eqnarray}
\begin{eqnarray}\label{eq4-4: lma: GDOI linear homomorphism}
\sum\limits_{q_1=1}^{m_{k_1,i_1}-1}\sum\limits_{q_2=1}^{m_{k_2,i_2}-1}\frac{\beta^{(q_1,q_2)}(\lambda_{k_1},\lambda_{k_1})\gamma(\lambda_{k_1},\lambda_{k_1})}{q_1! q_2!}\bm{N}^{q_1}_{k_1,i_1}\bm{Y}\bm{N}^{q_2}_{k_2,i_2}\nonumber \\
+\sum\limits_{q_1=1}^{m_{k_1,i_1}-1}\sum\limits_{q_2=1}^{m_{k_2,i_2}-1}\frac{\beta^{(q_1,-)}(\lambda_{k_1},\lambda_{k_1})\gamma^{(-, q_2)}(\lambda_{k_2},\lambda_{k_2})}{q_1! q_2!}\bm{N}_{k_1,i_1}^{q_1}\bm{Y}\bm{N}^{q_2}_{k_2,i_2}\nonumber \\
+\sum\limits_{q_1=1}^{m_{k_1,i_1}-1}\sum\limits_{q_2=1}^{m_{k_2,i_2}-1}\frac{\beta^{(-,q_2)}(\lambda_{k_1},\lambda_{k_1})\gamma^{(q_1,-)}(\lambda_{k_2},\lambda_{k_2})}{q_1! q_2!}\bm{N}_{k_1,i_1}^{q_1}\bm{Y}\bm{N}^{q_2}_{k_2,i_2}\nonumber \\
+\sum\limits_{q_1=1}^{m_{k_1,i_1}-1}\sum\limits_{q_2=1}^{m_{k_2,i_2}-1}\frac{\beta(\lambda_{k_1},\lambda_{k_1})\gamma^{(q_1,q_2)}(\lambda_{k_1},\lambda_{k_1})}{q_1! q_2!}\bm{N}_{k_1,i_1}^{q_1}\bm{Y}\bm{N}^{q_2}_{k_2,i_2}\nonumber \\
+\sum\limits_{q_1=1}^{m_{k_1,i_1}-1}\sum\limits_{q'_2+ q''_2=1}^{q'_2 + q''_2=m_{k_2,i_2}-1}\frac{\beta^{(q_1, q'_2)}(\lambda_{k_1},\lambda_{k_1})\gamma^{(-, q''_2)}(\lambda_{k_2},\lambda_{k_2})}{q_1! q'_2! q''_2!}\bm{N}_{k_1,i_1}^{q_1}\bm{Y}\bm{N}^{q'_2+q''_2}_{k_2,i_2}\nonumber \\
+\sum\limits_{q_1=1}^{m_{k_1,i_1}-1}\sum\limits_{q'_2+q''_2=1}^{q'_2 + q''_2=m_{k_2,i_2}-1}\frac{\beta^{(-, q'_2)}(\lambda_{k_1},\lambda_{k_1})\gamma^{(q_1, q''_2)}(\lambda_{k_2},\lambda_{k_2}) }{q_1! q'_2! q''_2!}\bm{N}_{k_1,i_1}^{q_1}\bm{Y}\bm{N}^{q'_2+q''_2}_{k_2,i_2}\nonumber \\
+\sum\limits_{q'_1+ q''_1=1}^{q'_1 + q''_1=m_{k_1,i_1}-1}\sum\limits_{q_2=1}^{m_{k_2,i_2}-1}\frac{\beta^{(q'_1, q_2)}(\lambda_{k_1},\lambda_{k_1})\gamma^{(q''_1, -)}(\lambda_{k_2},\lambda_{k_2})}{q'_1! q_2! q''_1!}\bm{N}^{q'_1+q''_1}_{k_1,i_1}\bm{Y}\bm{N}_{k_2,i_2}^{q_2}\nonumber \\
+\sum\limits_{q'_1+q''_1=1}^{q'_1 + q''_1=m_{k_1,i_1}-1}\sum\limits_{q_2=1}^{m_{k_2,i_2}-1}\frac{\beta^{(q'_1, -)}(\lambda_{k_1},\lambda_{k_1})\gamma^{(q''_1,q_2)}(\lambda_{k_2},\lambda_{k_2})}{q'_1! q_2! q''_1!}\bm{N}^{q'_1+q''_1}_{k_1,i_1}\bm{Y}\bm{N}_{k_2,i_2}^{q_2}\nonumber \\
+\sum\limits_{q'_1+ q''_1=1}^{q'_1 + q''_1=m_{k_1,i_1}-1}\sum\limits_{q'_2+q''_2=1}^{q'_2 + q''_2=m_{k_2,i_2}-1}\frac{\beta^{(q'_1, q'_2)}(\lambda_{k_1},\lambda_{k_1})\gamma^{(q''_1,q''_2)}(\lambda_{k_2},\lambda_{k_2})}{q'_1! q'_2! q''_1!q''_2!}\bm{N}^{q'_1+q''_1}_{k_1,i_1}\bm{Y}\bm{N}_{k_2,i_2}^{q'_2+q''_2}.
\end{eqnarray}

Recall the Leibniz rule for partial derivatives of a product of two functions \( f(x, y) \) and \( g(x, y) \) is:
\begin{eqnarray}\label{eq5: lma: GDOI linear homomorphism}
(f(x, y)g(x, y))^{(q_1, q_2)} = \sum_{k_1=0}^{q_1} \sum_{k_2=0}^{q_2} \binom{q_1}{k_1} \binom{q_2}{k_2} f^{(k_1, k_2)}(x, y) g^{(q_1 - k_1, q_2 - k_2)}(x, y),
\end{eqnarray}
where \( f^{(k_1, k_2)} \) denotes the partial derivative of \( f(x, y) \) with respect to \( x \) for \( k_1 \) times and with respect to \( y \) for \( k_2 \) times. 

It is clear that the Part 1 in Eq.~\eqref{eq3: lma: GDOI linear homomorphism} is identical to Eq.~\eqref{eq4-1: lma: GDOI linear homomorphism}. For the reamining Part 2, Part 3, and Part 4 in Eq.~\eqref{eq3: lma: GDOI linear homomorphism}, they are identical to Eq.~\eqref{eq4-2: lma: GDOI linear homomorphism}, Eq.~\eqref{eq4-3: lma: GDOI linear homomorphism}, and Eq.~\eqref{eq4-4: lma: GDOI linear homomorphism}, respectively by applying the Leibniz rule for partial derivatives given by Eq.~\eqref{eq5: lma: GDOI linear homomorphism}. Therefore, Eq.~\eqref{eq1-2: lma: GDOI linear homomorphism} is also valid.
$\hfill\Box$

The linear homomorphism property given by Lemma~\ref{lma: GDOI linear homomorphism} can be enhanced to linear isomorphism if more stronger conditions are provided to the function $\beta(z_1,z_2)$ and the space of $T^{\bm{X}_1,\bm{X}_2}_{\beta}$.

Let us consider a fintie set of pair of complex numbers, denoted by $\mathcal{P}$, which is given by
\begin{eqnarray}\label{eq1: telescope func}
\mathcal{P} = \{ ( z_1^{(i)},  z_2^{(i)}) \}_{i=1}^{N} \subset \mathcal{D} \subset \mathbb{C}^2.
\end{eqnarray}
Each point \( ( z_1^{(i)},  z_2^{(i)}) \) in the domain \( \mathcal{D} \) has an associated highest derivative order \( (p_i, q_i) \), meaning we impose the conditions:
\begin{eqnarray}\label{eq2: telescope func}
\frac{\partial^p}{\partial z_1^p} \frac{\partial^q}{\partial z_2^q} f( z_1^{(i)},  z_2^{(i)}) = \frac{\partial^p}{\partial z_1^p} \frac{\partial^q}{\partial z_2^q} g( z_1^{(i)},  z_2^{(i)}), \quad \forall \ 0 \leq p \leq p_i, \ 0 \leq q \leq q_i, \ i=1, \dots, N,
\end{eqnarray}
where $f,g$ are two analytic functions defined over $\mathcal{D}$. This defines a new hierarchy of function spaces indexed by multi-indices:
\begin{eqnarray}\label{eq2: telescope func}
\mathcal{F}_{(m_1,n_1), (m_2,n_2), \dots, (m_N, n_N)},
\end{eqnarray}
which consists of all analytic functions in \( \mathcal{D} \) that are uniquely determined by their values and derivatives at all points in \( \mathcal{P} \).

We construct a \emph{nested hierarchy} of function spaces:
\begin{eqnarray}\label{eq2: telescope func}
\mathcal{F}_{(0,0),\dots, (0,0)} \subseteq \mathcal{F}_{(1,0), (0,1), \dots} \subseteq \dots \subseteq \mathcal{F}_{(m_1,n_1), (m_2,n_2),\dots} \subseteq \mathcal{F}_{\infty, \infty}
\end{eqnarray}
where we have 
\begin{enumerate}
\item \( \mathcal{F}_{(0,0), \dots, (0,0)} \): Functions uniquely determined only by function values at each \( ( z_1^{(i)},  z_2^{(i)}) \).
\item \( \mathcal{F}_{(m_1,n_1), (m_2,n_2), \dots} \): Functions uniquely determined by derivatives up to \( (m_i, n_i) \) at each \( ( z_1^{(i)},  z_2^{(i)}) \).
\item \( \mathcal{F}_{\infty, \infty} \): The full space of analytic functions uniquely determined by infinite-order Taylor series at all given points.
\end{enumerate}
Since each \( (m_i, n_i) \) determines how much local information is captured, these spaces naturally telescope, meaning higher-order spaces contain lower-order ones. 

For the uniqueness conditions for the family of functions given by Eq.~\eqref{eq2: telescope func}, we will consider two cases (A) and (B) discussed below. 

(A) Uniqueness in \( \mathcal{F}_{\infty,\infty} \): \\ 

If two functions \( f, g \) satisfy:
\begin{eqnarray}\label{eq3: telescope func}
\frac{\partial^m}{\partial z_1^m} \frac{\partial^n}{\partial z_2^n} f( z_1^{(i)},  z_2^{(i)}) = \frac{\partial^m}{\partial z_1^m} \frac{\partial^n}{\partial z_2^n} g( z_1^{(i)},  z_2^{(i)}), \quad \forall m, n, i,
\end{eqnarray}
then \( f \equiv g \) everywhere in \( \mathcal{D} \), because the Taylor series expansions at each \( ( z_1^{(i)},  z_2^{(i)}) \) uniquely determine the function globally.

(B)  If \( (m_i, n_i) \) are Finite: \\

If \( (m_i, n_i) \) are finite for all \( i \), then there may be multiple functions satisfying the constraints. However, if we restrict $f,g$ as finite degree bivariate polynomials, then \( f \equiv g \) everywhere in \( \mathcal{D} \). Hence, the space \( \mathcal{F}_{(m_1,n_1),\dots, (m_N, n_N)} \) will contain polynomials of finite degree in \( z_1, z_2 \) if restricted to polynomial spaces. The solutions can be written as:
\begin{eqnarray}\label{eq4: telescope func}
  f(z_1, z_2) = P(z_1, z_2) + g(z_1, z_2),
\end{eqnarray}
where \( P(z_1, z_2) \) is a bivariate polynomial satisfying the constraints from their values and derivatives at all points in \( \mathcal{P} \)., and \( g(z_1, z_2) \) is any analytic function vanishing at \( ( z_1^{(i)},  z_2^{(i)}) \) up to the prescribed derivative orders. Therefore, we will reformulation the function space given by Eq.~\eqref{eq2: telescope func} by using quotient spaces of banishing parts.  

Let \( \mathcal{O}(\mathcal{D}) \) be the space of all analytic functions in a domain \( \mathcal{D} \subset \mathbb{C}^2 \). First, we define:  
\begin{enumerate}
\item The subspace \( \mathcal{F}_{(m_1, n_1), \ldots, (m_N, n_N)} \) as the set of analytic functions that satisfy derivative constraints up to order \( (m_i, n_i) \) at points \( ( z_1^{(i)},  z_2^{(i)}) \).
\item The subspace \( \mathcal{N}_{(m_1, n_1), \ldots, (m_N, n_N)} \) as the set of all analytic functions that vanish (including derivatives up to \( (m_i, n_i) \)) at each \( ( z_1^{(i)},  z_2^{(i)}) \).
\end{enumerate}
Then, we can define the quotient space:
\begin{eqnarray}\label{eq5: telescope func}
\mathcal{Q}_{(m_1, n_1), \ldots, (m_N, n_N)}&=&\mathcal{O}(D)\bigg/\mathcal{N}_{(m_1, n_1), \ldots, (m_N, n_N)},
\end{eqnarray}
which classifies analytic functions modulo functions that vanish at the given constraints.

If we restrict to polynomial spaces, then the quotient space \( \mathcal{Q} \) is finite-dimensional. The space contains polynomials of degree at most \( (m, n) \), which are uniquely determined by their values and derivatives at the prescribed points. On the other hand, if we allow general analytic functions, the quotient space still classifies functions uniquely up to an equivalence class of vanishing functions. Then, any function in \( \mathcal{O}(D) \) can be written as:
\begin{eqnarray}\label{eq6: telescope func}
     f(z_1, z_2) = P(z_1, z_2) + g(z_1, z_2)
\end{eqnarray}
where \( P(z_1, z_2) \) is a canonical representative (typically a polynomial satisfying constraints), and \( g(z_1, z_2) \) belongs to \( \mathcal{N} \), meaning it vanishes up to prescribed derivatives at given points.

Since the function spaces telescope in terms of derivative order, the quotient spaces also inherit a hierarchical structure:
\begin{eqnarray}\label{eq7: telescope func}
\mathcal{Q}_{(0,0),\dots, (0,0)} \subseteq \mathcal{Q}_{(1,0), (0,1), \dots} \subseteq \dots \subseteq \mathcal{Q}_{(m_1,n_1), (m_2,n_2),\dots} \subseteq \mathcal{Q}_{\infty, \infty},
\end{eqnarray}
where we have
\begin{enumerate}
\item \( \mathcal{Q}_{(0,0),..., (0,0)} \) consists of functions determined only by function values.
\item \( \mathcal{Q}_{(m_1, n_1), \ldots, (m_N, n_N)} \) consists of functions determined up to higher derivative orders.
\item \( \mathcal{Q}_{\infty, \infty} \) is trivial (i.e., contains only one equivalence class) because all analytic functions are uniquely determined.
\end{enumerate}

Therefore, using quotient space notation, we can classify function spaces as:
\begin{eqnarray}\label{eq8: telescope func}
\mathcal{Q}_{(m_1, n_1), \ldots, (m_N, n_N)} = \mathcal{O}(D)\bigg/\mathcal{N}_{(m_1, n_1),\ldots, (m_N, n_N)},
\end{eqnarray}
where the numerator \( \mathcal{O}(D) \) is the space of analytic functions, and the denominator \( \mathcal{N} \) captures the constraints, i.e., the space of functions that vanish up to prescribed derivative orders.

Given two matrices $\bm{X}_1$ and $\bm{X}_2$ with their nilpotent orders $m_{k_1,i_1}$ and $m_{k_2,i_2}$, we define the following symbols, denoted by $\mathfrak{P}(\bm{X}_1,\bm{X}_2)$, to represent highest derivative order pairs with respect to each pair of eigenvalues of the matrix $\bm{X}_1$ and eigenvalues of the matrix $\bm{X}_2$ in $T_{\beta}^{\bm{X}_1,\bm{X}_2}$:
\begin{eqnarray}\label{eq8: telescope func}
\mathfrak{P}(\bm{X}_1,\bm{X}_2)&\define&\prod\limits_{\{k_1\}\times\{k_2\}}(\max\limits_{i_1}m_{k_1,i_1},\max\limits_{i_2}m_{k_2,i_2}),
\end{eqnarray}
where $\{k_1\}$ and $\{k_2\}$ are indices for eigenvalues of the matrix $\bm{X}_1$ and eigenvalues of the matrix $\bm{X}_2$, respectively. 

The main purpose of this section is to present the following Theorem~\ref{thm: GDOI algebraic properties} used to characterize the algebraic properties of  $T^{\bm{X}_1,\bm{X}_2}_\beta$.
\begin{theorem}\label{thm: GDOI algebraic properties}
Let $\bm{X}_1$ and $\bm{X}_2$ be arbitrary matrices with their eigenvalues in the domain $\mathcal{D}$. Consider the operator $T^{\bm{X}_1,\bm{X}_2}_\beta$, which is defined by Eq.~\eqref{eq1:  GDOI def}. Given $\beta(z_1,z_2)$ as a bivariate analytic function over the domain $\mathcal{D}$, the mapping $f: \beta \rightarrow T^{\bm{X}_1,\bm{X}_2}_{\beta}$ is a linear isomorphism if $\beta \in \mathcal{Q}_{\mathfrak{P}(\bm{X}_1,\bm{X}_2)}$ and $T^{\bm{X}_1,\bm{X}_2}_{\beta}(\bm{Y})\in\mathcal{S}_{\bm{X}_1, \bm{X}_2, \bm{Y}}(\mathcal{D})$.
\end{theorem}
\textbf{Proof:}
From Lemma~\ref{lma: GDOI linear homomorphism}, we know that the mapping is a linear homomorphism. This theorem is proved if the injective and surjective are satisfied by the mapping $f$.  

To show that the mapping $f$ is injective, we first assume that $T^{\bm{X}_1,\bm{X}_2}_{\beta}=T^{\bm{X}_1,\bm{X}_2}_{\gamma}$, i.e., $\bm{0} = T^{\bm{X}_1,\bm{X}_2}_{\beta}-T^{\bm{X}_1,\bm{X}_2}_{\gamma}$. Then, we have
\begin{eqnarray}\label{eq6: thm: GDOI algebraic properties}
\bm{0}&=&\sum\limits_{k_1=1}^{K_1}\sum\limits_{k_2=1}^{K_2}\sum\limits_{i_1=1}^{\alpha_{k_1}^{(\mathrm{G})}}\sum\limits_{i_2=1}^{\alpha_{k_2}^{(\mathrm{G})}}(\beta(\lambda_{k_1}, \lambda_{k_2})-\gamma(\lambda_{k_1}, \lambda_{k_2}))\bm{P}_{k_1,i_1}\bm{Y}\bm{P}_{k_2,i_2} \nonumber \\
&&+\sum\limits_{k_1=1}^{K_1}\sum\limits_{k_2=1}^{K_2}\sum\limits_{i_1=1}^{\alpha_{k_1}^{(\mathrm{G})}}\sum\limits_{i_2=1}^{\alpha_{k_2}^{(\mathrm{G})}}\sum_{q_2=1}^{m_{k_2,i_2}-1}\frac{\beta^{(-,q_2)}(\lambda_{k_1},\lambda_{k_2}) - \gamma^{(-,q_2)}(\lambda_{k_1},\lambda_{k_2})}{q_2!}\bm{P}_{k_1,i_1}\bm{Y}\bm{N}_{k_2,i_2}^{q_2} \nonumber \\
&&+\sum\limits_{k_1=1}^{K_1}\sum\limits_{k_2=1}^{K_2}\sum\limits_{i_1=1}^{\alpha_{k_1}^{(\mathrm{G})}}\sum\limits_{i_2=1}^{\alpha_{k_2}^{(\mathrm{G})}}\sum_{q_1=1}^{m_{k_1,i_1}-1}\frac{\beta^{(q_1,-)}(\lambda_{k_1},\lambda_{k_2})- \gamma^{(q_1,-)}(\lambda_{k_1},\lambda_{k_2})}{q_1!}\bm{N}_{k_1,i_1}^{q_1}\bm{Y}\bm{P}_{k_2,i_2}  \nonumber \\
&&+\sum\limits_{k_1=1}^{K_1}\sum\limits_{k_2=1}^{K_2}\sum\limits_{i_1=1}^{\alpha_{k_1}^{(\mathrm{G})}}\sum\limits_{i_2=1}^{\alpha_{k_2}^{(\mathrm{G})}}\sum_{q_1=1}^{m_{k_1,i_1}-1}\sum_{q_2=1}^{m_{k_2,i_2}-1}\frac{\beta^{(q_1,q_2)}(\lambda_{k_1},\lambda_{k_2}) - \gamma^{(q_1,q_2)}(\lambda_{k_1},\lambda_{k_2})}{q_1!q_2!}\bm{N}_{k_1,i_1}^{q_1}\bm{Y}\bm{N}_{k_2,i_2}^{q_2}.
\end{eqnarray}
From Lemma~\ref{lma: ind of Proj and Nilp} and the definition of the function space $\mathcal{Q}_{\mathfrak{P}(\bm{X}_1,\bm{X}_2)}$, we have $\beta=\gamma$ in the domain $\mathcal{D}$. 

The surjective of the mapping $f$ is easy to check from the definition of $\mathcal{S}_{\bm{X}_1, \bm{X}_2, \bm{Y}}(\mathcal{D})$ given by Eq.~\eqref{eq0: lma: GDOI linear homomorphism} and Eq.~\eqref{eq0-1: lma: GDOI linear homomorphism}. 
$\hfill\Box$

\section{Norm Estimations}\label{sec: Norm Estimations}

In this section, we will provide norm estimation for the GDOI: $T_{\beta}^{\bm{X}_1,\bm{X}_2}(\bm{Y})$, defined by Eq.~\eqref{eq1:  GDOI def}. We have the following Theorem~\ref{thm: GDOI norm} about the upper and the lower estimations for the norm of $T_{\beta}^{\bm{X}_1,\bm{X}_2}(\bm{Y})$. We consider Frobenius norm here, however, the approach can be extended easily to other matrix norms.

We will begin by presenting Lemma~\ref{lma:conv triangle for Frob norm} for the converse triangle inequality for the Frobenius norm, denoted by $\left\Vert \cdot \right\Vert$.
\begin{lemma}\label{lma:conv triangle for Frob norm}
Given $n$ matrices $\bm{A}_1, \bm{A}_2, \ldots, \bm{A}_n$ such that $\left\Vert\bm{A}_1\right\Vert\geq\left\Vert\bm{A}_2\right\Vert\geq\ldots\geq\left\Vert\bm{A}_n\right\Vert$, then, we have
\begin{eqnarray}\label{eq1:lma:conv triangle for Frob norm}
\left\Vert\sum\limits_{i=1}^n \bm{A}_i \right\Vert \geq \max\left[0, \left\Vert\bm{A}_1\right\Vert - \sum\limits_{i=2}^n \left\Vert\bm{A}_i\right\Vert\right]
\end{eqnarray}

\end{lemma}
\textbf{Proof:}
We will prove this lemma by induction. We will consder $n=2$ case first. For any two matrices \( \bm{A}_1, \bm{A}_2 \) of the same size, the  Frobenius norm satisfies the triangle inequality :  

\[
\|\bm{A}_1 + \bm{A}_2\| \leq \|\bm{A}_1\| + \|\bm{A}_2\|.
\]

Since the Frobenius norm satisfies the  parallelogram law :  

\[
\|\bm{A}_1 + \bm{A}_2\|^2 + \|\bm{A}_1 - \bm{A}_2\|^2 = 2\|\bm{A}_1\|^2 + 2\|\bm{A}_2\|^2,
\]

we apply the  polarization identity :  

\[
\|\bm{A}_1 + \bm{A}_2\|^2 = \|\bm{A}_1\|^2 + \|\bm{A}_2\|^2 + 2 \langle \bm{A}_1, \bm{A}_2 \rangle_F,
\]

where the  Frobenius inner product  is:  

\[
\langle \bm{A}_1, \bm{A}_2 \rangle_F = \mathrm{Tr}(\bm{A}_1^* \bm{A}_2).
\]

Taking the square root, we get:  

\[
\|\bm{A}_1 + \bm{A}_2\| = \sqrt{\|\bm{A}_1\|^2 + \|\bm{A}_2\|^2 + 2 \text{Re} \langle \bm{A}_1, \bm{A}_2 \rangle_F}.
\]

Using the  angle \( \theta \) between \( \bm{A}_1 \) and \( \bm{A}_2 \) in Frobenius inner product space, where:  

\[
\cos\theta = \frac{\langle \bm{A}_1, \bm{A}_2 \rangle_F}{\|\bm{A}_1\| \|\bm{A}_2\|},
\]

we can rewrite:  

\[
\|\bm{A}_1 + \bm{A}_2\| \geq \sqrt{\|\bm{A}_1\|^2 + \|\bm{A}_2\|^2 - 2 \|\bm{A}_1\| \|\bm{A}_2\|}.
\]

Since the worst case occurs when \( \bm{A}_1 \) and \( \bm{A}_2 \) are  negatively aligned, we get:  

\[
\|\bm{A}_1 + \bm{A}_2\| \geq |\|\bm{A}_1\| - \|\bm{A}_2\|| = \|\bm{A}_1\| - \|\bm{A}_2\|,
\]
where the last equality comes from the assumption that $\|\bm{A}_1\| \geq \|\bm{A}_2\|$.

For \( n = 3 \) with three matrices \( \bm{A}_1, \bm{A}_2, \bm{A}_3 \), we apply the Frobenius norm properties:

\[
\|\bm{A}_1 + \bm{A}_2 + \bm{A}_3\|  \geq \big| \|\bm{A}_1\|  - \|\bm{A}_2 + \bm{A}_3\|  \big|.
\]

 Case 1: \( \|\bm{A}_1\|  \geq \|\bm{A}_2 + \bm{A}_3\|  \)   \\

If  

\[
\|\bm{A}_1\|  \geq \|\bm{A}_2 + \bm{A}_3\| \mbox{~but $ \|\bm{A}_2 + \bm{A}_3\| \leq \|\bm{A}_2\|  + \|\bm{A}_3\|$},
\]

then
\begin{eqnarray}\label{eq2:lma:conv triangle for Frob norm}
\|\bm{A}_1 + \bm{A}_2 + \bm{A}_3\|  \geq \|\bm{A}_1\|  - \|\bm{A}_2 + \bm{A}_3\|  \geq \|\bm{A}_1\|  - \|\bm{A}_2\|  - \|\bm{A}_3\|.
\end{eqnarray}

 Case 2: \( \|\bm{A}_1\|  \leq \|\bm{A}_2 + \bm{A}_3\|  \)   \\

Then, using the reverse triangle inequality:

\[
\|\bm{A}_1 + \bm{A}_2 + \bm{A}_3\|  \geq \|\bm{A}_2 + \bm{A}_3\|  - \|\bm{A}_1\|.
\]

Since  

\[
\|\bm{A}_2 + \bm{A}_3\|  \geq \big| \|\bm{A}_2\|  - \|\bm{A}_3\|  \big|,
\]

we get
\begin{eqnarray}\label{eq3:lma:conv triangle for Frob norm}
\|\bm{A}_1 + \bm{A}_2 + \bm{A}_3\|  \geq \|\bm{A}_2\|  - \|\bm{A}_3\|  - \|\bm{A}_1\|.
\end{eqnarray}
By combining Eq.~\eqref{eq2:lma:conv triangle for Frob norm} and Eq.~\eqref{eq3:lma:conv triangle for Frob norm}, we have
\begin{eqnarray}
\|\bm{A}_1 + \bm{A}_2 + \bm{A}_3\| &\geq&\max\big[\|\bm{A}_2\|  - \|\bm{A}_1\|, \|\bm{A}_1\|  - \|\bm{A}_2\|\big] - \|\bm{A}_3\| \nonumber \\
&\geq&\|\bm{A}_1\|  - \|\bm{A}_2\|  - \|\bm{A}_3\|.
\end{eqnarray}
Therefore, we still have Eq.~\eqref{eq1:lma:conv triangle for Frob norm} valid for $n=3$.

For \( n \geq 4 \), applying the same reasoning iteratively, we obtain:

\[
\|\bm{A}_1 + \bm{A}_2 + \dots + \bm{A}_n\|  \geq \max\left[0, \left\Vert\bm{A}_1\right\Vert - \sum\limits_{i=2}^n \left\Vert\bm{A}_i\right\Vert\right].
\]

This result shows how the  Frobenius norm follows a similar lower bound structure as the absolute value case, leveraging the  reverse triangle inequality.
$\hfill\Box$

\begin{theorem}\label{thm: GDOI norm}
We have the upper bound for the Frobenius norm of $T_{\beta}^{\bm{X}_1,\bm{X}_2}(\bm{Y})$, which is given by
\begin{eqnarray}\label{eq1: thm: GDOI norm}
\lefteqn{\left\Vert T_{\beta}^{\bm{X}_1,\bm{X}_2}(\bm{Y}) \right\Vert\leq}\nonumber \\
&&\left[\max\limits_{\lambda_1 \in \Lambda_{\bm{X}_1},\lambda_2 \in \Lambda_{\bm{X}_2}} \left\vert\beta(\lambda_1,\lambda_2)\right\vert\right]\left\Vert\bm{Y}\right\Vert \nonumber \\
&&+\sum\limits_{k_2=1}^{K_2}\sum\limits_{i_2=1}^{\alpha_{k_2}^{(\mathrm{G})}}\sum_{q_2=1}^{m_{k_2,i_2}-1}\left[\max\limits_{\lambda_1 \in \Lambda_{\bm{X}_1},\lambda_2 \in \Lambda_{\bm{X}_2}}\left\vert\frac{\beta^{(-,q_2)}(\lambda_{1},\lambda_{2})}{q_2!}\right\vert\right]\left\Vert\bm{N}_{k_2,i_2}^{q_2}\right\Vert\left\Vert\bm{Y}\right\Vert\nonumber \\
&&+\sum\limits_{k_1=1}^{K_1}\sum\limits_{i_1=1}^{\alpha_{k_1}^{(\mathrm{G})}}\sum_{q_1=1}^{m_{k_1,i_1}-1}\left[\max\limits_{\lambda_1 \in \Lambda_{\bm{X}_1},\lambda_2 \in \Lambda_{\bm{X}_2}}\left\vert\frac{\beta^{(q_1,-)}(\lambda_{1},\lambda_{2})}{q_1!}\right\vert\right]\left\Vert\bm{N}_{k_1,i_1}^{q_1}\right\Vert\left\Vert\bm{Y}\right\Vert  \nonumber \\
&&+ \sum\limits_{k_1=1}^{K_1}\sum\limits_{k_2=1}^{K_2}\sum\limits_{i_1=1}^{\alpha_{k_1}^{(\mathrm{G})}}\sum\limits_{i_2=1}^{\alpha_{k_2}^{(\mathrm{G})}}\sum_{q_1=1}^{m_{k_1,i_1}-1}\sum_{q_2=1}^{m_{k_2,i_2}-1}\left[\max\limits_{\lambda_1 \in \Lambda_{\bm{X}_1},\lambda_2 \in \Lambda_{\bm{X}_2}}\left\vert\frac{\beta^{(q_1,q_2)}(\lambda_{1},\lambda_{2})}{q_1!q_2!}\right\vert\right]\nonumber\\
&&~~~\times\left\Vert\bm{N}_{k_1,i_1}^{q_1}\right\Vert\left\Vert\bm{N}_{k_2,i_2}^{q_2}\right\Vert\left\Vert\bm{Y}\right\Vert.
\end{eqnarray}
where $\Lambda_{\bm{X}_1}$ and $\Lambda_{\bm{X}_2}$ are spectrums of the matrix $\bm{X}_1$ and the matrix $\bm{X}_2$, respectively.

On the other hand, let us define the following matrices
\begin{eqnarray}\label{eq1-1: thm: GDOI norm}
\bm{A}_1&\define& \sum\limits_{k_1=1}^{K_1}\sum\limits_{k_2=1}^{K_2}\sum\limits_{i_1=1}^{\alpha_{k_1}^{(\mathrm{G})}}\sum\limits_{i_2=1}^{\alpha_{k_2}^{(\mathrm{G})}}\beta(\lambda_{k_1}, \lambda_{k_2})\bm{P}_{k_1,i_1}\bm{Y}\bm{P}_{k_2,i_2}\nonumber \\
\bm{A}_2&\define& \sum\limits_{k_1=1}^{K_1}\sum\limits_{k_2=1}^{K_2}\sum\limits_{i_1=1}^{\alpha_{k_1}^{(\mathrm{G})}}\sum\limits_{i_2=1}^{\alpha_{k_2}^{(\mathrm{G})}}\sum_{q_2=1}^{m_{k_2,i_2}-1}\frac{\beta^{(-,q_2)}(\lambda_{k_1},\lambda_{k_2})}{q_2!}\bm{P}_{k_1,i_1}\bm{Y}\bm{N}_{k_2,i_2}^{q_2}\nonumber \\
\bm{A}_3&\define&\sum\limits_{k_1=1}^{K_1}\sum\limits_{k_2=1}^{K_2}\sum\limits_{i_1=1}^{\alpha_{k_1}^{(\mathrm{G})}}\sum\limits_{i_2=1}^{\alpha_{k_2}^{(\mathrm{G})}}\sum_{q_1=1}^{m_{k_1,i_1}-1}\frac{\beta^{(q_1,-)}(\lambda_{k_1},\lambda_{k_2})}{q_1!}\bm{N}_{k_1,i_1}^{q_1}\bm{Y}\bm{P}_{k_2,i_2}\nonumber \\
\bm{A}_4&\define&\sum\limits_{k_1=1}^{K_1}\sum\limits_{k_2=1}^{K_2}\sum\limits_{i_1=1}^{\alpha_{k_1}^{(\mathrm{G})}}\sum\limits_{i_2=1}^{\alpha_{k_2}^{(\mathrm{G})}}\sum_{q_1=1}^{m_{k_1,i_1}-1}\sum_{q_2=1}^{m_{k_2,i_2}-1}\frac{\beta^{(q_1,q_2)}(\lambda_{k_1},\lambda_{k_2})}{q_1!q_2!}\bm{N}_{k_1,i_1}^{q_1}\bm{Y}\bm{N}_{k_2,i_2}^{q_2}.
\end{eqnarray}
then, we have the lower bound for the Frobenius norm of $T_{\beta}^{\bm{X}_1,\bm{X}_2}(\bm{Y})$, which is given by
\begin{eqnarray}\label{eq2: thm: GDOI norm}
\left\Vert T_{\beta}^{\bm{X}_1,\bm{X}_2}(\bm{Y}) \right\Vert\geq \max\left[0, \left\Vert\bm{A}_{\sigma(1)}\right\Vert - \sum\limits_{i=2}^4 \left\Vert\bm{A}_{\sigma(i)}\right\Vert\right],
\end{eqnarray}
where $\sigma$ is the permutation of matrices $\bm{A}_i$ for $i=1,2,3,4$ such that $\left\Vert\bm{A}_{\sigma(1)}\right\Vert \geq \left\Vert\bm{A}_{\sigma(2)}\right\Vert \geq\left\Vert\bm{A}_{\sigma(3)}\right\Vert \geq\left\Vert\bm{A}_{\sigma(4)}\right\Vert \geq$. 

Further, if we have $\left[\min\limits_{\lambda_1 \in \Lambda_{\bm{X}_1},\lambda_2 \in \Lambda_{\bm{X}_2}} \beta(\lambda_1,\lambda_2)\right]\left\Vert\bm{Y}\right\Vert\geq \left\Vert\bm{A}_{2}\right\Vert+\left\Vert\bm{A}_{3}\right\Vert + \left\Vert\bm{A}_{4}\right\Vert$,  the lower bound for the Frobenius norm of $T_{\beta}^{\bm{X}_1,\bm{X}_2}(\bm{Y})$ can be expressed by
\begin{eqnarray}\label{eq3: thm: GDOI norm}
\left\Vert T_{\beta}^{\bm{X}_1,\bm{X}_2}(\bm{Y}) \right\Vert\geq \left[\min\limits_{\lambda_1 \in \Lambda_{\bm{X}_1},\lambda_2 \in \Lambda_{\bm{X}_2}} \beta(\lambda_1,\lambda_2)\right]\left\Vert\bm{Y}\right\Vert - \sum\limits_{i=2}^4 \left\Vert\bm{A}_{i}\right\Vert.
\end{eqnarray}
\end{theorem}
\textbf{Proof:}
From the definition of $T_{\beta}^{\bm{X}_1,\bm{X}_2}(\bm{Y})$ and the triangle inequality, we have 
\begin{eqnarray}\label{eq4: thm: GDOI norm}
\left\Vert T_{\beta}^{\bm{X}_1,\bm{X}_2}(\bm{Y}) \right\Vert&\leq&\underbrace{\left\Vert\sum\limits_{k_1=1}^{K_1}\sum\limits_{k_2=1}^{K_2}\sum\limits_{i_1=1}^{\alpha_{k_1}^{(\mathrm{G})}}\sum\limits_{i_2=1}^{\alpha_{k_2}^{(\mathrm{G})}}\beta(\lambda_{k_1}, \lambda_{k_2})\bm{P}_{k_1,i_1}\bm{Y}\bm{P}_{k_2,i_2}\right\Vert}_{\mbox{Part I}}\nonumber \\
&&+\underbrace{\left\Vert\sum\limits_{k_1=1}^{K_1}\sum\limits_{k_2=1}^{K_2}\sum\limits_{i_1=1}^{\alpha_{k_1}^{(\mathrm{G})}}\sum\limits_{i_2=1}^{\alpha_{k_2}^{(\mathrm{G})}}\sum_{q_2=1}^{m_{k_2,i_2}-1}\frac{\beta^{(-,q_2)}(\lambda_{k_1},\lambda_{k_2})}{q_2!}\bm{P}_{k_1,i_1}\bm{Y}\bm{N}_{k_2,i_2}^{q_2}\right\Vert}_{\mbox{Part II}}\nonumber \\
&&+\underbrace{\left\Vert\sum\limits_{k_1=1}^{K_1}\sum\limits_{k_2=1}^{K_2}\sum\limits_{i_1=1}^{\alpha_{k_1}^{(\mathrm{G})}}\sum\limits_{i_2=1}^{\alpha_{k_2}^{(\mathrm{G})}}\sum_{q_1=1}^{m_{k_1,i_1}-1}\frac{\beta^{(q_1,-)}(\lambda_{k_1},\lambda_{k_2})}{q_1!}\bm{N}_{k_1,i_1}^{q_1}\bm{Y}\bm{P}_{k_2,i_2}\right\Vert}_{\mbox{Part III}}\nonumber \\
&&+\underbrace{\left\Vert\sum\limits_{k_1=1}^{K_1}\sum\limits_{k_2=1}^{K_2}\sum\limits_{i_1=1}^{\alpha_{k_1}^{(\mathrm{G})}}\sum\limits_{i_2=1}^{\alpha_{k_2}^{(\mathrm{G})}}\sum_{q_1=1}^{m_{k_1,i_1}-1}\sum_{q_2=1}^{m_{k_2,i_2}-1}\frac{\beta^{(q_1,q_2)}(\lambda_{k_1},\lambda_{k_2})}{q_1!q_2!}\bm{N}_{k_1,i_1}^{q_1}\bm{Y}\bm{N}_{k_2,i_2}^{q_2}\right\Vert}_{\mbox{Part IV}}.
\end{eqnarray}

For Part I, we have 
\begin{eqnarray}\label{eq4-1: thm: GDOI norm}
\lefteqn{\underbrace{\left\Vert\sum\limits_{k_1=1}^{K_1}\sum\limits_{k_2=1}^{K_2}\sum\limits_{i_1=1}^{\alpha_{k_1}^{(\mathrm{G})}}\sum\limits_{i_2=1}^{\alpha_{k_2}^{(\mathrm{G})}}\beta(\lambda_{k_1}, \lambda_{k_2})\bm{P}_{k_1,i_1}\bm{Y}\bm{P}_{k_2,i_2}\right\Vert}_{\mbox{Part I}}}\nonumber \\
&\leq&\left[\max\limits_{\lambda_1 \in \Lambda_{\bm{X}_1},\lambda_2 \in \Lambda_{\bm{X}_2}} \left\vert\beta(\lambda_1,\lambda_2)\right\vert\right]\left\Vert\sum\limits_{k_1=1}^{K_1}\sum\limits_{k_2=1}^{K_2}\sum\limits_{i_1=1}^{\alpha_{k_1}^{(\mathrm{G})}}\sum\limits_{i_2=1}^{\alpha_{k_2}^{(\mathrm{G})}}\bm{P}_{k_1,i_1}\bm{Y}\bm{P}_{k_2,i_2}\right\Vert\nonumber \\
&=&\left[\max\limits_{\lambda_1 \in \Lambda_{\bm{X}_1},\lambda_2 \in \Lambda_{\bm{X}_2}} \left\vert\beta(\lambda_1,\lambda_2)\right\vert\right]\left\Vert\bm{Y}\right\Vert,
\end{eqnarray}
where we use the fact that $\sum\limits_{k_1=1}^{K_1}\sum\limits_{i_1=1}^{\alpha_{k_1}^{(\mathrm{G})}}\bm{P}_{k_1,i_1} = \bm{I}$ and $\sum\limits_{k_2=1}^{K_2}\sum\limits_{i_2=1}^{\alpha_{k_2}^{(\mathrm{G})}}\bm{P}_{k_2,i_2} = \bm{I}$ at the last equality. 

For Part II, we have 
\begin{eqnarray}\label{eq4-2: thm: GDOI norm}
\lefteqn{\underbrace{\left\Vert\sum\limits_{k_1=1}^{K_1}\sum\limits_{k_2=1}^{K_2}\sum\limits_{i_1=1}^{\alpha_{k_1}^{(\mathrm{G})}}\sum\limits_{i_2=1}^{\alpha_{k_2}^{(\mathrm{G})}}\sum_{q_2=1}^{m_{k_2,i_2}-1}\frac{\beta^{(-,q_2)}(\lambda_{k_1},\lambda_{k_2})}{q_2!}\bm{P}_{k_1,i_1}\bm{Y}\bm{N}_{k_2,i_2}^{q_2}\right\Vert}_{\mbox{Part II}}}\nonumber \\
&\leq&\sum\limits_{k_2=1}^{K_2}\sum\limits_{i_2=1}^{\alpha_{k_2}^{(\mathrm{G})}}\sum_{q_2=1}^{m_{k_2,i_2}-1}\left[\max\limits_{\lambda_1 \in \Lambda_{\bm{X}_1},\lambda_2 \in \Lambda_{\bm{X}_2}} \left\vert\frac{\beta^{(-,q_2)}(\lambda_{1},\lambda_{2})}{q_2!}\right\vert\right]\left\Vert\sum\limits_{k_1=1}^{K_1}\sum\limits_{i_1=1}^{\alpha_{k_1}^{(\mathrm{G})}}\bm{P}_{k_1,i_1}\bm{Y}\bm{N}_{k_2,i_2}^{q_2}\right\Vert\nonumber \\
&\leq_1&\sum\limits_{k_2=1}^{K_2}\sum\limits_{i_2=1}^{\alpha_{k_2}^{(\mathrm{G})}}\sum_{q_2=1}^{m_{k_2,i_2}-1}\left[\max\limits_{\lambda_1 \in \Lambda_{\bm{X}_1},\lambda_2 \in \Lambda_{\bm{X}_2}}\left\vert\frac{\beta^{(-,q_2)}(\lambda_{1},\lambda_{2})}{q_2!}\right\vert\right]\left\Vert\bm{N}_{k_2,i_2}^{q_2}\right\Vert\left\Vert\bm{Y}\right\Vert,
\end{eqnarray}
where we use the fact that $\sum\limits_{k_1=1}^{K_1}\sum\limits_{i_1=1}^{\alpha_{k_1}^{(\mathrm{G})}}\bm{P}_{k_1,i_1} = \bm{I}$ and the norm multiplicative inequality at the inequality $\leq_1$.

For Part III, we have 
\begin{eqnarray}\label{eq4-3: thm: GDOI norm}
\lefteqn{\underbrace{\left\Vert\sum\limits_{k_1=1}^{K_1}\sum\limits_{k_2=1}^{K_2}\sum\limits_{i_1=1}^{\alpha_{k_1}^{(\mathrm{G})}}\sum\limits_{i_2=1}^{\alpha_{k_2}^{(\mathrm{G})}}\sum_{q_1=1}^{m_{k_1,i_1}-1}\frac{\beta^{(q_1,-)}(\lambda_{k_1},\lambda_{k_2})}{q_1!}\bm{N}_{k_1,i_1}^{q_1}\bm{Y}\bm{P}_{k_2,i_2}\right\Vert}_{\mbox{Part III}}}\nonumber \\
&\leq&\sum\limits_{k_1=1}^{K_1}\sum\limits_{i_1=1}^{\alpha_{k_1}^{(\mathrm{G})}}\sum_{q_1=1}^{m_{k_1,i_1}-1}\left[\max\limits_{\lambda_1 \in \Lambda_{\bm{X}_1},\lambda_2 \in \Lambda_{\bm{X}_2}} \left\vert\frac{\beta^{(q_1,-)}(\lambda_{1},\lambda_{2})}{q_1!}\right\vert\right]\left\Vert\sum\limits_{k_2=1}^{K_2}\sum\limits_{i_2=1}^{\alpha_{k_2}^{(\mathrm{G})}}\bm{N}_{k_1,i_1}^{q_1}\bm{Y}\bm{P}_{k_2,i_2}\right\Vert\nonumber \\
&\leq_1&\sum\limits_{k_1=1}^{K_1}\sum\limits_{i_1=1}^{\alpha_{k_1}^{(\mathrm{G})}}\sum_{q_1=1}^{m_{k_1,i_1}-1}\left[\max\limits_{\lambda_1 \in \Lambda_{\bm{X}_1},\lambda_2 \in \Lambda_{\bm{X}_2}}\left\vert\frac{\beta^{(q_1,-)}(\lambda_{1},\lambda_{2})}{q_1!}\right\vert\right]\left\Vert\bm{N}_{k_1,i_1}^{q_1}\right\Vert\left\Vert\bm{Y}\right\Vert,
\end{eqnarray}
where we use the fact that $\sum\limits_{k_2=1}^{K_2}\sum\limits_{i_2=1}^{\alpha_{k_2}^{(\mathrm{G})}}\bm{P}_{k_2,i_2} = \bm{I}$ and the norm multiplicative inequality at the  inequality $\leq_1$.

For Part IV, we have 
\begin{eqnarray}\label{eq4-4: thm: GDOI norm}
\lefteqn{\underbrace{\left\Vert\sum\limits_{k_1=1}^{K_1}\sum\limits_{k_2=1}^{K_2}\sum\limits_{i_1=1}^{\alpha_{k_1}^{(\mathrm{G})}}\sum\limits_{i_2=1}^{\alpha_{k_2}^{(\mathrm{G})}}\sum_{q_1=1}^{m_{k_1,i_1}-1}\sum_{q_2=1}^{m_{k_2,i_2}-1}\frac{\beta^{(q_1,q_2)}(\lambda_{k_1},\lambda_{k_2})}{q_1!q_2!}\bm{N}_{k_1,i_1}^{q_1}\bm{Y}\bm{N}_{k_2,i_2}^{q_2}\right\Vert}_{\mbox{Part IV}}}\nonumber \\
&\leq& \sum\limits_{k_1=1}^{K_1}\sum\limits_{k_2=1}^{K_2}\sum\limits_{i_1=1}^{\alpha_{k_1}^{(\mathrm{G})}}\sum\limits_{i_2=1}^{\alpha_{k_2}^{(\mathrm{G})}}\sum_{q_1=1}^{m_{k_1,i_1}-1}\sum_{q_2=1}^{m_{k_2,i_2}-1}\left[\max\limits_{\lambda_1 \in \Lambda_{\bm{X}_1},\lambda_2 \in \Lambda_{\bm{X}_2}}\left\vert\frac{\beta^{(q_1,q_2)}(\lambda_{1},\lambda_{2})}{q_1!q_2!}\right\vert\right]\left\Vert\bm{N}_{k_1,i_1}^{q_1}\bm{Y}\bm{N}_{k_2,i_2}^{q_2}\right\Vert \nonumber \\
&\leq_1& \sum\limits_{k_1=1}^{K_1}\sum\limits_{k_2=1}^{K_2}\sum\limits_{i_1=1}^{\alpha_{k_1}^{(\mathrm{G})}}\sum\limits_{i_2=1}^{\alpha_{k_2}^{(\mathrm{G})}}\sum_{q_1=1}^{m_{k_1,i_1}-1}\sum_{q_2=1}^{m_{k_2,i_2}-1}\left[\max\limits_{\lambda_1 \in \Lambda_{\bm{X}_1},\lambda_2 \in \Lambda_{\bm{X}_2}}\left\vert\frac{\beta^{(q_1,q_2)}(\lambda_{1},\lambda_{2})}{q_1!q_2!}\right\vert\right]\nonumber\\
&&~~~\times\left\Vert\bm{N}_{k_1,i_1}^{q_1}\right\Vert\left\Vert\bm{N}_{k_2,i_2}^{q_2}\right\Vert\left\Vert\bm{Y}\right\Vert.
\end{eqnarray}
where we use the norm multiplicative inequality again at the inequality $\leq_1$. By combining previous four parts from Eq.~\eqref{eq4-1: thm: GDOI norm} to Eq.~\eqref{eq4-4: thm: GDOI norm}, we have the desired upper bound given by Eq.~\eqref{eq1: thm: GDOI norm}

For the lower bound of $\left\Vert T_{\beta}^{\bm{X}_1,\bm{X}_2}(\bm{Y})\right\Vert$, we have Eq.~\eqref{eq2: thm: GDOI norm} immediatedly from Lemma~\ref{lma:conv triangle for Frob norm}.

If we have $\left[\min\limits_{\lambda_1 \in \Lambda_{\bm{X}_1},\lambda_2 \in \Lambda_{\bm{X}_2}} \beta(\lambda_1,\lambda_2)\right]\left\Vert\bm{Y}\right\Vert\geq \left\Vert\bm{A}_{2}\right\Vert+\left\Vert\bm{A}_{3}\right\Vert + \left\Vert\bm{A}_{4}\right\Vert$ and Lemma~\ref{lma:conv triangle for Frob norm}, we have
\begin{eqnarray}\label{eq5: thm: GDOI norm}
\left\Vert T_{\beta}^{\bm{X}_1,\bm{X}_2}(\bm{Y})\right\Vert&\geq&\left\Vert\bm{A}_1\right\Vert-(\left\Vert\bm{A}_{2}\right\Vert+\left\Vert\bm{A}_{3}\right\Vert + \left\Vert\bm{A}_{4}\right\Vert) \nonumber \\
&\geq&\left[\min\limits_{\lambda_1 \in \Lambda_{\bm{X}_1},\lambda_2 \in \Lambda_{\bm{X}_2}} \beta(\lambda_1,\lambda_2)\right]\left\Vert\sum\limits_{k_1=1}^{K_1}\sum\limits_{k_2=1}^{K_2}\sum\limits_{i_1=1}^{\alpha_{k_1}^{(\mathrm{G})}}\sum\limits_{i_2=1}^{\alpha_{k_2}^{(\mathrm{G})}}\bm{P}_{k_1,i_1}\bm{Y}\bm{P}_{k_2,i_2}\right\Vert\nonumber \\
& &-  (\left\Vert\bm{A}_{2}\right\Vert+\left\Vert\bm{A}_{3}\right\Vert + \left\Vert\bm{A}_{4}\right\Vert) \nonumber \\
&=&\left[\min\limits_{\lambda_1 \in \Lambda_{\bm{X}_1},\lambda_2 \in \Lambda_{\bm{X}_2}} \beta(\lambda_1,\lambda_2)\right]\left\Vert\bm{Y}\right\Vert - (\left\Vert\bm{A}_{2}\right\Vert+\left\Vert\bm{A}_{3}\right\Vert + \left\Vert\bm{A}_{4}\right\Vert),
\end{eqnarray}
which is the lower bound of $\left\Vert T_{\beta}^{\bm{X}_1,\bm{X}_2}(\bm{Y})\right\Vert$ given by Eq.~\eqref{eq3: thm: GDOI norm}.  
$\hfill\Box$

\section{Perturbation Formula and Lipschitz Estimations}\label{sec:Perturbation Formula and Lipschitz Estimation}

The purpose of this section is to derive the perturbation formula for the difference between $f(\bm{X}_1)$ and $f(\bm{X}_2)$ via GDOI and establish Lipschitz estimations for the upper bound and the lower bound.

\subsection{Perturbation Formuila and Nilpotent Part Difference Characterization}\label{sec: Perturbation Formuila and Nilpotant Part Difference Characterization}

We have the following Theorem~\ref{thm:pert formula} to determine the difference between $f(\bm{X}_1)$ and $f(\bm{X}_2)$ via GDOI.
\begin{theorem}\label{thm:pert formula}
Given an analytic function $f(z)$ within the domain for $|z| < R$, the first matrix $\bm{X}_1$ with the dimension $m$ and $K_1$ distinct eigenvalues $\lambda_{k_1}$ for $k_1=1,2,\ldots,K_1$ such that
\begin{eqnarray}\label{eq1-1:thm:pert formula}
\bm{X}_1&=&\sum\limits_{k_1=1}^{K_1}\sum\limits_{i_1=1}^{\alpha_{k_1}^{\mathrm{G}}} \lambda_{k_1} \bm{P}_{k_1,i_1}+
\sum\limits_{k_1=1}^{K_1}\sum\limits_{i_1=1}^{\alpha_{k_1}^{\mathrm{G}}} \bm{N}_{k_1,i_1},
\end{eqnarray}
where $\left\vert\lambda_{k_1}\right\vert<R$, and second matrix $\bm{X}_2$ with the dimension $m$ and $K_2$ distinct eigenvalues $\lambda_{k_2}$ for $k_2=1,2,\ldots,K_2$ such that
\begin{eqnarray}\label{eq1-2:thm:pert formula}
\bm{X}_2&=&\sum\limits_{k_2=1}^{K_2}\sum\limits_{i_2=1}^{\alpha_{k_2}^{\mathrm{G}}} \lambda_{k_2} \bm{P}_{k_2,i_2}+
\sum\limits_{k_2=1}^{K_2}\sum\limits_{i_2=1}^{\alpha_{k_2}^{\mathrm{G}}} \bm{N}_{k_2,i_2},
\end{eqnarray}
where $\left\vert\lambda_{k_2}\right\vert<R$. We also assume that $\lambda_{k_1} \neq \lambda_{k_2}$ for any $k_1$ and $k_2$. Then, we have
\begin{eqnarray}\label{eq2:thm:pert formula}
f(\bm{X}_1)\bm{Y}-\bm{Y}f(\bm{X}_2)&=&T_{\frac{f(x_1) - f(x_2)}{x_1 -x_2}}^{\bm{X}_{1},\bm{X}_{2}}(\bm{X}_{1}\bm{Y}-\bm{Y}\bm{X}_{2}). 
\end{eqnarray}
\end{theorem}
\textbf{Proof:}
We set two projection functions $\pi_1(x_1,x_2)$ and $\pi_2(x_1,x_2)$ as
\begin{eqnarray}\label{eq3:thm:pert formula}
\pi_1(x_1,x_2)=x_1,~~\pi_2(x_1,x_2)&=&x_2.
\end{eqnarray}

Then, from GDOI definition given by Eq.~\eqref{eq1:  GDOI def}, we have
\begin{eqnarray}\label{eq4-1:thm:pert formula}
T_{\pi_1}^{\bm{X}_1,\bm{X}_2}(\bm{Y})&=&\sum\limits_{k_1=1}^{K_1}\sum\limits_{k_2=1}^{K_2}\sum\limits_{i_1=1}^{\alpha_{k_1}^{(\mathrm{G})}}\sum\limits_{i_2=1}^{\alpha_{k_2}^{(\mathrm{G})}}\lambda_{k_1}\bm{P}_{k_1,i_1}\bm{Y}\bm{P}_{k_2,i_2} \nonumber \\
&&+\sum\limits_{k_1=1}^{K_1}\sum\limits_{k_2=1}^{K_2}\sum\limits_{i_1=1}^{\alpha_{k_1}^{(\mathrm{G})}}\sum\limits_{i_2=1}^{\alpha_{k_2}^{(\mathrm{G})}}\bm{N}_{k_1,i_1}\bm{Y}\bm{P}_{k_2,i_2}  \nonumber \\
&=&\bm{X}_1 \bm{Y}.
\end{eqnarray}
Similarly, we have
\begin{eqnarray}\label{eq4-2:thm:pert formula}
T_{\pi_2}^{\bm{X}_1,\bm{X}_2}(\bm{Y})&=&\bm{Y}\bm{X}_2.
\end{eqnarray}
Besides, we also have 
\begin{eqnarray}\label{eq5-1:thm:pert formula}
T_{f \circ \pi_1}^{\bm{X}_1,\bm{X}_2}(\bm{Y})&=&\sum\limits_{k_1=1}^{K_1}\sum\limits_{k_2=1}^{K_2}\sum\limits_{i_1=1}^{\alpha_{k_1}^{(\mathrm{G})}}\sum\limits_{i_2=1}^{\alpha_{k_2}^{(\mathrm{G})}}f(\lambda_{k_1})\bm{P}_{k_1,i_1}\bm{Y}\bm{P}_{k_2,i_2} \nonumber \\
&&+\sum\limits_{k_1=1}^{K_1}\sum\limits_{k_2=1}^{K_2}\sum\limits_{i_1=1}^{\alpha_{k_1}^{(\mathrm{G})}}\sum\limits_{i_2=1}^{\alpha_{k_2}^{(\mathrm{G})}}\sum_{q_1=1}^{m_{k_1,i_1}-1}\frac{f^{(q_1)}(\lambda_{k_1})}{q_1!}\bm{N}_{k_1,i_1}^{q_1}\bm{Y}\bm{P}_{k_2,i_2}  \nonumber \\
&=&f(\bm{X}_1)\bm{Y}.
\end{eqnarray}
Similarly, we have
\begin{eqnarray}\label{eq5-2:thm:pert formula}
T_{f \circ \pi_2}^{\bm{X}_1,\bm{X}_2}(\bm{Y})&=&\bm{Y}f(\bm{X}_2).
\end{eqnarray}

Then, we have
\begin{eqnarray}\label{eq6:thm:pert formula}
f(\bm{X}_1)\bm{Y} - \bm{Y}f(\bm{X}_2)&=_1&T_{f \circ \pi_1}^{\bm{X}_1,\bm{X}_2}(\bm{Y}) - T_{f \circ \pi_2}^{\bm{X}_1,\bm{X}_2}(\bm{Y})\nonumber \\
&=&T_{f \circ \pi_1 - f \circ \pi_2}^{\bm{X}_1,\bm{X}_2}(\bm{Y})\nonumber \\
&=&T_{\frac{f(x_1) - f(x_2)}{x_1 - x_2} (\pi_1 - \pi_2)}^{\bm{X}_1,\bm{X}_2}(\bm{Y})\nonumber \\
&=&T_{\frac{f(x_1) - f(x_2)}{x_1 - x_2}\pi_1}^{\bm{X}_1,\bm{X}_2}(\bm{Y}) - T_{\frac{f(x_1) - f(x_2)}{x_1 - x_2}\pi_2}^{\bm{X}_1,\bm{X}_2}(\bm{Y})\nonumber \\
&=&T_{\frac{f(x_1) - f(x_2)}{x_1 - x_2}}^{\bm{X}_1,\bm{X}_2}(T_{\pi_1}^{\bm{X}_1,\bm{X}_2}(\bm{Y})) - T_{\frac{f(x_1) - f(x_2)}{x_1 - x_2}}^{\bm{X}_1,\bm{X}_2}(T_{\pi_2}^{\bm{X}_1,\bm{X}_2}(\bm{Y}))\nonumber \\
&=_2&T_{\frac{f(x_1) - f(x_2)}{x_1 - x_2}}^{\bm{X}_1,\bm{X}_2}(\bm{X}_1\bm{Y}-\bm{Y}\bm{X}_2),
\end{eqnarray}
where we apply Eq.~\eqref{eq5-1:thm:pert formula} and Eq.~\eqref{eq5-2:thm:pert formula} in $=_1$, we apply Eq.~\eqref{eq4-1:thm:pert formula} and Eq.~\eqref{eq4-2:thm:pert formula} in $=_2$ and other equalities come from Lemma~\ref{lma: GDOI linear homomorphism}.
$\hfill\Box$

\begin{corollary}\label{cor:pert formula I}
By the same setting as Theorem~\ref{thm:pert formula}, we have
\begin{eqnarray}\label{eq1:cor:pert formula I}
f(\bm{X}_1)-f(\bm{X}_2)&=&T_{\frac{f(x_1) - f(x_2)}{x_1 -x_2}}^{\bm{X}_{1},\bm{X}_{2}}(\bm{X}_{1}-\bm{X}_{2}). 
\end{eqnarray}
\end{corollary}
\textbf{Proof:}
From Thoerem~\ref{thm:pert formula}, this corollar is established by setting $\bm{Y}=\bm{I}$.
$\hfill\Box$

We have the following Theorem~\ref{thm:pert formula diff DOI} to characterize the difference between the GDOI of $T_{\frac{f(x_1) - f(x_2)}{x_1 -x_2}}^{\bm{X}_{1},\bm{X}_{2}}(\bm{X}_{1}-\bm{X}_{2})$ and the conventional DOI defined by projction parts of parameter matrices $\bm{X}_1$ and $\bm{X}_2$, i.e., \\
$T_{\frac{f(x_1) - f(x_2)}{x_1 -x_2}}^{\bm{X}_{1,\bm{P}},\bm{X}_{2,\bm{P}}}(\bm{X}_{1,\bm{P}}-\bm{X}_{2,\bm{P}})$ from Corollary~\ref{cor:pert formula I}
\begin{theorem}\label{thm:pert formula diff DOI}
Given an analytic function $f(z)$ within the domain for $|z| < R$, the first matrix $\bm{X}_1$ with the dimension $m$ and $K_1$ distinct eigenvalues $\lambda_{k_1}$ for $k_1=1,2,\ldots,K_1$ such that
\begin{eqnarray}\label{eq1-1:thm:pert formula diff DOI}
\bm{X}_1&=&\sum\limits_{k_1=1}^{K_1}\sum\limits_{i_1=1}^{\alpha_{k_1}^{\mathrm{G}}} \lambda_{k_1} \bm{P}_{k_1,i_1}+
\sum\limits_{k_1=1}^{K_1}\sum\limits_{i_1=1}^{\alpha_{k_1}^{\mathrm{G}}} \bm{N}_{k_1,i_1}\nonumber \\
&\define&\bm{X}_{1,\bm{P}}+\bm{X}_{1,\bm{N}},
\end{eqnarray}
where $\left\vert\lambda_{k_1}\right\vert<R$, and second matrix $\bm{X}_2$ with the dimension $m$ and $K_2$ distinct eigenvalues $\lambda_{k_2}$ for $k_2=1,2,\ldots,K_2$ such that
\begin{eqnarray}\label{eq1-2:thm:pert formula diff DOI}
\bm{X}_2&=&\sum\limits_{k_2=1}^{K_2}\sum\limits_{i_2=1}^{\alpha_{k_2}^{\mathrm{G}}} \lambda_{k_2} \bm{P}_{k_2,i_2}+
\sum\limits_{k_2=1}^{K_2}\sum\limits_{i_2=1}^{\alpha_{k_2}^{\mathrm{G}}} \bm{N}_{k_2,i_2}\nonumber \\
&\define&\bm{X}_{2,\bm{P}}+\bm{X}_{2,\bm{N}},
\end{eqnarray}
where $\left\vert\lambda_{k_2}\right\vert<R$. We also assume that $\lambda_{k_1} \neq \lambda_{k_2}$ for any $k_1$ and $k_2$. Then, we have
\begin{eqnarray}\label{eq2:thm:pert formula diff DOI}
T_{\frac{f(x_1) - f(x_2)}{x_1 -x_2}}^{\bm{X}_{1},\bm{X}_{2}}(\bm{X}_{1}-\bm{X}_{2})&=&T_{\frac{f(x_1) - f(x_2)}{x_1 -x_2}}^{\bm{X}_{1,\bm{P}},\bm{X}_{2,\bm{P}}}(\bm{X}_{1,\bm{P}}-\bm{X}_{2,\bm{P}})+\left[T_{f(x_1)}^{\bm{X}_{1,\bm{N}},\bm{X}_{2,\bm{P}}}(\bm{I})-T_{f(x_2)}^{\bm{X}_{1,\bm{P}},\bm{X}_{2,\bm{N}}}(\bm{I})\right].
\end{eqnarray}
\end{theorem}
\textbf{Proof:}
From Theorem 1 in~\cite{chang2024operatorChar}, we have 
\begin{eqnarray}\label{eq3-1:thm:pert formula diff DOI}
f(\bm{X}_1)&=&\sum\limits_{k_1=1}^{K_1} \left[\sum\limits_{i_1=1}^{\alpha_{k_1}^{(\mathrm{G})}}f(\lambda_{k_1})\bm{P}_{k_1,i_1}+\sum\limits_{i_1=1}^{\alpha_{k_1}^{(\mathrm{G})}}\sum\limits_{q_1=1}^{m_{k_1,i_1}-1}\frac{f^{(q_1)}(\lambda_{k_1})}{q_1!}\bm{N}_{k_1,i_1}^{q_1}\right],
\end{eqnarray}
and
\begin{eqnarray}\label{eq3-2:thm:pert formula diff DOI}
f(\bm{X}_2)&=&\sum\limits_{k_2=1}^{K_2} \left[\sum\limits_{i_2=1}^{\alpha_{k_2}^{(\mathrm{G})}}f(\lambda_{k_2})\bm{P}_{k_2,i_2}+\sum\limits_{i_2=1}^{\alpha_{k_2}^{(\mathrm{G})}}\sum\limits_{q_2=1}^{m_{k_2,i_2}-1}\frac{f^{(q_2)}(\lambda_{k_2})}{q_2!}\bm{N}_{k_2,i_2}^{q_2}\right].
\end{eqnarray}
Then, we can express $f(\bm{X}_1)-f(\bm{X}_2)$ as
\begin{eqnarray}\label{eq4:thm:pert formula diff DOI}
\lefteqn{f(\bm{X}_1)-f(\bm{X}_2)}\nonumber \\
&=&\underbrace{\left[\sum\limits_{k_1=1}^{K_1}\sum\limits_{i_1=1}^{\alpha_{k_1}^{(\mathrm{G})}}f(\lambda_{k_1})\bm{P}_{k_1,i_1}- \sum\limits_{k_2=1}^{K_2}\sum\limits_{i_2=1}^{\alpha_{k_2}^{(\mathrm{G})}}f(\lambda_{k_2})\bm{P}_{k_2,i_2}\right]}_{\mbox{Part I}}\nonumber \\
&&+\underbrace{\left[\sum\limits_{k_1=1}^{K_1}\sum\limits_{i_1=1}^{\alpha_{k_1}^{(\mathrm{G})}}\sum\limits_{q_1=1}^{m_{k_1,i_1}-1}\frac{f^{(q_1)}(\lambda_{k_1})}{q_1!}\bm{N}_{k_1,i_1}^{q_1}-\sum\limits_{k_2=1}^{K_2}\sum\limits_{i_2=1}^{\alpha_{k_2}^{(\mathrm{G})}}\sum\limits_{q_2=1}^{m_{k_2,i_2}-1}\frac{f^{(q_2)}(\lambda_{k_2})}{q_2!}\bm{N}_{k_2,i_2}^{q_2}\right]}_{\mbox{Part II}}.
\end{eqnarray}

From GDOI definition given by Eq.~\eqref{eq1:  GDOI def}, we have
\begin{eqnarray}\label{eq5:thm:pert formula diff DOI}
T_{\frac{f(x_1) - f(x_2)}{x_1 -x_2}}^{\bm{X}_{1,\bm{P}},\bm{X}_{2,\bm{P}}}(\bm{X}_{1,\bm{P}}-\bm{X}_{2,\bm{P}})&=&\sum\limits_{k_1=1}^{K_1}\sum\limits_{k_2=1}^{K_2}\sum\limits_{i_1=1}^{\alpha_{k_1}^{(\mathrm{G})}}\sum\limits_{i_2=1}^{\alpha_{k_2}^{(\mathrm{G})}}\frac{f(\lambda_{k_1}) - f(\lambda_{k_2})}{\lambda_{k_1} -\lambda_{k_2}}\bm{P}_{k_1,i_1}(\bm{X}_{1,\bm{P}}-\bm{X}_{2,\bm{P}})\bm{P}_{k_2,i_2} \nonumber \\
&=&\sum\limits_{k_1=1}^{K_1}\sum\limits_{k_2=1}^{K_2}\sum\limits_{i_1=1}^{\alpha_{k_1}^{(\mathrm{G})}}\sum\limits_{i_2=1}^{\alpha_{k_2}^{(\mathrm{G})}}\frac{f(\lambda_{k_1}) - f(\lambda_{k_2})}{\lambda_{k_1} -\lambda_{k_2}}\bm{P}_{k_1,i_1}\bm{X}_{1,\bm{P}}\bm{P}_{k_2,i_2} \nonumber \\
&&-\sum\limits_{k_1=1}^{K_1}\sum\limits_{k_2=1}^{K_2}\sum\limits_{i_1=1}^{\alpha_{k_1}^{(\mathrm{G})}}\sum\limits_{i_2=1}^{\alpha_{k_2}^{(\mathrm{G})}}\frac{f(\lambda_{k_1}) - f(\lambda_{k_2})}{\lambda_{k_1} -\lambda_{k_2}}\bm{P}_{k_1,i_1}\bm{X}_{2,\bm{P}}\bm{P}_{k_2,i_2} \nonumber \\
&=_1&\sum\limits_{k_1=1}^{K_1}\sum\limits_{k_2=1}^{K_2}\sum\limits_{i_1=1}^{\alpha_{k_1}^{(\mathrm{G})}}\sum\limits_{i_2=1}^{\alpha_{k_2}^{(\mathrm{G})}}\lambda_{k_1}\frac{f(\lambda_{k_1}) - f(\lambda_{k_2})}{\lambda_{k_1} -\lambda_{k_2}}\bm{P}_{k_1,i_1}\bm{P}_{k_2,i_2} \nonumber \\
&&-\sum\limits_{k_1=1}^{K_1}\sum\limits_{k_2=1}^{K_2}\sum\limits_{i_1=1}^{\alpha_{k_1}^{(\mathrm{G})}}\sum\limits_{i_2=1}^{\alpha_{k_2}^{(\mathrm{G})}}\lambda_{k_2}\frac{f(\lambda_{k_1}) - f(\lambda_{k_2})}{\lambda_{k_1} -\lambda_{k_2}}\bm{P}_{k_1,i_1}\bm{P}_{k_2,i_2} \nonumber \\
&=&\sum\limits_{k_1=1}^{K_1}\sum\limits_{k_2=1}^{K_2}\sum\limits_{i_1=1}^{\alpha_{k_1}^{(\mathrm{G})}}\sum\limits_{i_2=1}^{\alpha_{k_2}^{(\mathrm{G})}}f(\lambda_{k_1})\bm{P}_{k_1,i_1}\bm{P}_{k_2,i_2} \nonumber \\
&&-\sum\limits_{k_1=1}^{K_1}\sum\limits_{k_2=1}^{K_2}\sum\limits_{i_1=1}^{\alpha_{k_1}^{(\mathrm{G})}}\sum\limits_{i_2=1}^{\alpha_{k_2}^{(\mathrm{G})}}f(\lambda_{k_2})\bm{P}_{k_1,i_1}\bm{P}_{k_2,i_2} \nonumber \\
&=_2&\underbrace{\left[\sum\limits_{k_1=1}^{K_1}\sum\limits_{i_1=1}^{\alpha_{k_1}^{(\mathrm{G})}}f(\lambda_{k_1})\bm{P}_{k_1,i_1}- \sum\limits_{k_2=1}^{K_2}\sum\limits_{i_2=1}^{\alpha_{k_2}^{(\mathrm{G})}}f(\lambda_{k_2})\bm{P}_{k_2,i_2}\right]}_{\mbox{Part I}},
\end{eqnarray}
where we apply $\bm{X}_{1,\bm{P}}$ and $\bm{X}_{2,\bm{P}}$ definitions in $=_1$ with relations provided by Eq.~\eqref{eq3: lma: ind of Proj and Nilp}, and apply relations $\sum\limits_{k_1=1}^{K_1}\sum\limits_{i_1=1}^{\alpha_{k_1}^{(\mathrm{G})}}\bm{P}_{k_1,i_1} = \sum\limits_{k_2=1}^{K_2}\sum\limits_{i_2=1}^{\alpha_{k_2}^{(\mathrm{G})}}\bm{P}_{k_2,i_2} = \bm{I}$ in $=_2$. 

On the other hand, by applying the GDOI definition given by Eq.~\eqref{eq1:  GDOI def}, we also have
\begin{eqnarray}\label{eq6:thm:pert formula diff DOI}
\lefteqn{\left[T_{f(x_1)}^{\bm{X}_{1,\bm{N}},\bm{X}_{2,\bm{P}}}(\bm{I})-T_{f(x_2)}^{\bm{X}_{1,\bm{P}},\bm{X}_{2,\bm{N}}}(\bm{I})\right]}\nonumber \\
&=&\sum\limits_{k_1=1}^{K_1}\sum\limits_{i_1=1}^{\alpha_{k_1}^{(\mathrm{G})}}\sum\limits_{k_2=1}^{K_2}\sum\limits_{i_2=1}^{\alpha_{k_2}^{(\mathrm{G})}}\sum\limits_{q_1=1}^{m_{k_1,i_1}-1}\frac{f^{(q_1)}(\lambda_{k_1})}{q_1!}\bm{N}_{k_1,i_1}^{q_1}\bm{I}\bm{P}_{k_2,i_2}-\nonumber \\
&&\sum\limits_{k_1=1}^{K_1}\sum\limits_{i_1=1}^{\alpha_{k_1}^{(\mathrm{G})}}\sum\limits_{k_2=1}^{K_2}\sum\limits_{i_2=1}^{\alpha_{k_2}^{(\mathrm{G})}}\sum\limits_{q_2=1}^{m_{k_2,i_2}-1}\frac{f^{(q_2)}(\lambda_{k_2})}{q_2!}\bm{P}_{k_1,i_1}\bm{I}\bm{N}_{k_2,i_2}^{q_2}\nonumber \\
&=_1&\underbrace{\left[\sum\limits_{k_1=1}^{K_1}\sum\limits_{i_1=1}^{\alpha_{k_1}^{(\mathrm{G})}}\sum\limits_{q_1=1}^{m_{k_1,i_1}-1}\frac{f^{(q_1)}(\lambda_{k_1})}{q_1!}\bm{N}_{k_1,i_1}^{q_1}-\sum\limits_{k_2=1}^{K_2}\sum\limits_{i_2=1}^{\alpha_{k_2}^{(\mathrm{G})}}\sum\limits_{q_2=1}^{m_{k_2,i_2}-1}\frac{f^{(q_2)}(\lambda_{k_2})}{q_2!}\bm{N}_{k_2,i_2}^{q_2}\right]}_{\mbox{Part II}},
\end{eqnarray}
where we apply $\sum\limits_{k_1=1}^{K_1}\sum\limits_{i_1=1}^{\alpha_{k_1}^{(\mathrm{G})}}\bm{P}_{k_1,i_1} = \sum\limits_{k_2=1}^{K_2}\sum\limits_{i_2=1}^{\alpha_{k_2}^{(\mathrm{G})}}\bm{P}_{k_2,i_2} = \bm{I}$ in $=_1$.

This theorem is proved by combining Eq.~\eqref{eq4:thm:pert formula diff DOI}, Eq.~\eqref{eq5:thm:pert formula diff DOI} and Eq.~\eqref{eq6:thm:pert formula diff DOI}.
$\hfill\Box$

Comparing Theorem~\ref{thm:pert formula diff DOI} with the conventional perturbation formula under Hermitian matrices assumptions of $\bm{X}_1$ and $\bm{X}_2$, we have one extra term, which is $\left[T_{f(x_1)}^{\bm{X}_{1,\bm{N}},\bm{X}_{2,\bm{P}}}(\bm{I})-T_{f(x_2)}^{\bm{X}_{1,\bm{P}},\bm{X}_{2,\bm{N}}}(\bm{I})\right]$. In the remaining part of this subsection, we will explore properties of the following extra term
\begin{eqnarray}\label{eq:extra nilpotent term in pertf}
\mu(\bm{X}_1, \bm{X}_2, f)&\define&\left[T_{f(x_1)}^{\bm{X}_{1,\bm{N}},\bm{X}_{2,\bm{P}}}(\bm{I})-T_{f(x_2)}^{\bm{X}_{1,\bm{P}},\bm{X}_{2,\bm{N}}}(\bm{I})\right].
\end{eqnarray}

Below, we have to provide a Lemma~\ref{lma: tuple is total order} about the total-ordering relationship for a set of pairs $(\ell_1, \ell_2, r)$ under dictionary ordering, where $\ell_1,\ell_2$ are nonnegative integers and $r$ is nonnegative real numbers. The purpose of such total-ordering structure is to quantify the divergence between the matrix $\mu(\bm{X}_1, \bm{X}_2, f)$ and the zero matrix. 
\begin{lemma}\label{lma: tuple is total order}
The set 
\[
S = \{(\ell_1, \ell_2, r) \mid \ell_1, \ell_2 \in \mathbb{Z}_{\geq 0}, r \in \mathbb{R}_{\geq 0}\}
\]
can form a total-ordering set under lexicographical order.
\end{lemma}
\textbf{Proof:}
To prove thatwe need to define a relation \(\preceq\) on \(S\) that satisfies the total-orderin properties:
\begin{enumerate}
\item Reflexivity: \(a \preceq a\) for all \(a \in S\).
\item Antisymmetry: If \(a \preceq b\) and \(b \preceq a\), then \(a = b\).
\item Transitivity: If \(a \preceq b\) and \(b \preceq c\), then \(a \preceq c\).
\item Totality: For any \(a, b \in S\), either \(a \preceq b\) or \(b \preceq a\).
\end{enumerate}

For any \((\ell_1, \ell_2, r)\) and \((\ell_1', \ell_2', r')\) in \(S\), the lexicographical order is :
\[
(\ell_1, \ell_2, r) \preceq (\ell_1', \ell_2', r') 
\]
if and only if:
1. \(\ell_1 < \ell_1'\), or
2. \(\ell_1 = \ell_1'\) and \(\ell_2 < \ell_2'\), or
3. \(\ell_1 = \ell_1'\), \(\ell_2 = \ell_2'\), and \(r \leq r'\).

1. Reflexivity:\\ 
For any \((\ell_1, \ell_2, r) \in S\), we have:
\[
(\ell_1, \ell_2, r) \preceq (\ell_1, \ell_2, r)
\]
since \(r \leq r\).

2. Antisymmetry:\\
If 
\[
(\ell_1, \ell_2, r) \preceq (\ell_1', \ell_2', r') \text{ and } (\ell_1', \ell_2', r') \preceq (\ell_1, \ell_2, r),
\]
then by the lexicographical order, all components must be equal, implying:
\[
(\ell_1, \ell_2, r) = (\ell_1', \ell_2', r').
\]

3. Transitivity:\\
If 
\[
(\ell_1, \ell_2, r) \preceq (\ell_1', \ell_2', r') \text{ and } (\ell_1', \ell_2', r') \preceq (\ell_1'', \ell_2'', r''),
\]
then through the lexicographical rules, we can conclude:
\[
(\ell_1, \ell_2, r) \preceq (\ell_1'', \ell_2'', r'').
\]

4. Totality:\\
For any two elements \((\ell_1, \ell_2, r)\) and \((\ell_1', \ell_2', r')\) in \(S\), either:
\[
(\ell_1, \ell_2, r) \preceq (\ell_1', \ell_2', r') \text{ or } (\ell_1', \ell_2', r') \preceq (\ell_1, \ell_2, r).
\]
This is guaranteed by the lexicographical order, which always allows for comparison.

Therefore, the set \(S\) with the defined lexicographical order \(\preceq\) forms a total ordering set, satisfying all necessary properties for a total order.
$\hfill\Box$

Based on Lemma~\ref{lma: tuple is total order}, we have the following Proposition~\ref{prop: GDOI div DOI} to quantify the deviation between the Generalized Degree of Orthogonal Invariance (GDOI) and the conventional Degree of Orthogonal Invariance (DOI) due to the non-Hermitian nature of matrices $\bm{X}_1$ and $\bm{X}_2$

\begin{proposition}\label{prop: GDOI div DOI}
Let $\mu(\bm{X}_1, \bm{X}_2, f)$ be a matrix obtained by Eq.~\eqref{eq:extra nilpotent term in pertf}. The matrix $\mu(\bm{X}_1, \bm{X}_2, f)$ can be categorized as either non-nilpotent or nilpotent. Define $\ell_1$ as the number of non-zero eigenvalues of $\mu(\bm{X}_1, \bm{X}_2, f)$. The following conditions hold:

1. If $\ell_1 > 0$, the matrix $\mu(\bm{X}_1, \bm{X}_2, f)$ is non-nilpotent.\\

2. If $\ell_1 = 0$, the matrix $\mu(\bm{X}_1, \bm{X}_2, f)$ is nilpotent. \\

For the nilpotent case, let $\ell_2$ represent the nilpotent degree of $\mu(\bm{X}_1, \bm{X}_2, f)$, meaning $\ell_2$ is the smallest integer $k$ such that
\[
\mu^{k}(\bm{X}_1, \bm{X}_2, f) = \bm{0}.
\]

For two matrices $\mu(\bm{X}_1, \bm{X}_2, f)$ and $\mu(\bm{X}'_1, \bm{X}'_2, f')$, if $\ell_1 = \ell'_1$ and $\ell_2 = \ell'_2$, the Frobenius norm of the matrices $\mu(\bm{X}_1, \bm{X}_2, f)$ and $\mu(\bm{X}'_1, \bm{X}'_2, f')$ can be used to distinguish between them by the total-ordering structure via triple $(\ell_1, \ell_2, r)$. Moreover, if $(\ell_1, \ell_2, r) = (0,0,0)$, the GDOI reduces to the DOI, indicating that both matrices $\bm{X}_1$ and $\bm{X}_2$ are Hermitian.
\end{proposition}

We post two design questions here regarding conditions to make the matrix $\mu(\bm{X}_1, \bm{X}_2, f)$ become nilpotent. They are:
\begin{eqnarray}\label{eq: design Q1}
\mbox{(1) Given the matrices $\bm{X}_1$ and $\bm{X}_2$, what the analytic function $f$ to make the matrix $\mu(\bm{X}_1, \bm{X}_2, f)$}\nonumber \\
\mbox{become nilpotent,}~~~~~~~~~~~~~~~~~~~~~~~~~~~~~~~~~~~~~~~~~~~~~~~~~~~~~~~~~~~~~~~~~~~~~~~~~~~~~~~~~~~~~~~~~~~~~~~~~~~~~~~~~~~~~~~~~~~~~~~~
\end{eqnarray}
and
\begin{eqnarray}\label{eq: design Q2}
\mbox{(2) Given  the analytic function $f$, what matrices $\bm{X}_1$ and $\bm{X}_2$ to make the matrix $\mu(\bm{X}_1, \bm{X}_2, f)$}\nonumber \\
\mbox{become nilpotent.}~~~~~~~~~~~~~~~~~~~~~~~~~~~~~~~~~~~~~~~~~~~~~~~~~~~~~~~~~~~~~~~~~~~~~~~~~~~~~~~~~~~~~~~~~~~~~~~~~~~~~~~~~~~~~~~~~~~~~~~~
\end{eqnarray}

For the question posted by Eq.~\eqref{eq: design Q1}, if we have $\mathrm{det}(\mu(\bm{X}_1, \bm{X}_2, f) - \lambda \bm{I}) = \lambda^n$, where $n \times n$ is the dimension of the  matrices $\bm{X}_1$ and $\bm{X}_2$, the matrix $\mu(\bm{X}_1, \bm{X}_2, f)$ will be a nilpotent matrix.

\begin{proposition}\label{prop: strict triangular form make mu nilpotent}
As the setting in Theorem~\ref{thm:pert formula}, if both the matrix $\bm{N}_{k_1,i_1}$ and the matrix $\bm{N}_{k_2,i_2}$ are strict lower triangular form or strict upper triangular form, the matrix $\mu(\bm{X}_1, \bm{X}_2, f)$ is nilpotent with degree $\max\limits_{k_1,i_1,k_2,i_2} [m_{k_1,i_1}, m_{k_2,i_2}]$.
\end{proposition}
\textbf{Proof:}
The nilpotency of the matrix $\bm{N}_{k_2,i_2}$ is clear since the linear combination of matrices with strict lower triangular form or strict upper triangular form will still be a matrix with strict lower triangular form or strict upper triangular form, respectively.

The nilpotent degree of the matrix \( \mu(\bm{X}_1, \bm{X}_2, f) \) can be determined from its strict lower triangular or strict upper triangular structure, as the nilpotent degree of a strict triangular matrix corresponds to the position of the first nonzero off-diagonal array.
$\hfill\Box$

The next Proposition~\ref{prop: commutative make mu nilpotent} shows that the matrix $\mu(\bm{X}_1, \bm{X}_2, f)$ is a nilpotent matrix if nilpotent components in $\bm{X}_1$ and $\bm{X}_2$ commute each other. 

\begin{proposition}\label{prop: commutative make mu nilpotent}
As the setting in Theorem~\ref{thm:pert formula}, if $\bm{N}_{k_1,i_1}$ and $\bm{N}_{k_2,i_2}$ commute each other, the matrix $\mu(\bm{X}_1, \bm{X}_2, f)$ is nilpotent. The nilpotent index of the matrix $\mu(\bm{X}_1, \bm{X}_2, f)$ is $\min [n, \max\limits_{k_1,i_1}[m_{k_1,i_1}]+\max\limits_{k_2,i_2}[m_{k_2,i_2}]]$.
\end{proposition}
\textbf{Proof:}
We first prove a fact that $\bm{A}_1$ and $\bm{A}_2$ are commutative nilpotent matrices, then, for any scalars \(a, b \in \mathbb{C}\), the linear combination \(\bm{A} = a\bm{A}_1 + b\bm{A}_2\) is also nilpotent.

Suppose \(\bm{A}_1^{k_1} = 0\) and \(\bm{A}_2^{k_2} = 0\) where \(k_1\) and \(k_2\) are the nilpotent indices of \(\bm{A}_1\) and \(\bm{A}_2\), respectively. Consider the linear combination \(\bm{A} = a\bm{A}_1 + b\bm{A}_2\). We analyze \(\bm{A}^{k_1 + k_2}\) using the binomial theorem:

\[
\bm{A}^{k_1 + k_2} = (a\bm{A}_1 + b\bm{A}_2)^{k_1 + k_2} = \sum_{j=0}^{k_1 + k_2} \binom{k_1 + k_2}{j} a^j b^{k_1 + k_2 - j} \bm{A}_1^j \bm{A}_2^{k_1 + k_2 - j}
\]

Note that if \(j \geq k_1\), then \(\bm{A}_1^j = \bm{0}\). Also, if \(k_1 + k_2 - j \geq k_2\), then \(\bm{A}_2^{k_1 + k_2 - j} = 0\). Since at least one of these conditions is always true for every term in the sum, all terms in the expansion are zero. Therefore:

\[
\bm{A}^{k_1 + k_2} = \bm{0}
\]

Since there exists a finite power of \(\bm{A}\) that is zero, \(\bm{A}\) is nilpotent. The nilpotent index of \(\bm{A}\) is at most \(k_1 + k_2\).

Since we have
\begin{eqnarray}
\mu(\bm{X}_1, \bm{X}_2, f)&=&\sum\limits_{k_1=1}^{K_1}\sum\limits_{i_1=1}^{\alpha_{k_1}^{(\mathrm{G})}}\sum\limits_{q_1=1}^{m_{k_1,i_1}-1}\frac{f^{(q_1)}(\lambda_{k_1})}{q_1!}\bm{N}_{k_1,i_1}^{q_1}-\sum\limits_{k_2=1}^{K_2}\sum\limits_{i_2=1}^{\alpha_{k_2}^{(\mathrm{G})}}\sum\limits_{q_2=1}^{m_{k_2,i_2}-1}\frac{f^{(q_2)}(\lambda_{k_2})}{q_2!}\bm{N}_{k_2,i_2}^{q_2},
\end{eqnarray}
then, the nilpotent index for the first term in the above equation is $\max\limits_{k_1,i_1}[m_{k_1,i_1}]$ and the nilpotent index for the second term in the above equation is $\max\limits_{k_2,i_2}[m_{k_2,i_2}]$. From the above fact just proved, we have the nilpotent index of the matrix $\mu(\bm{X}_1, \bm{X}_2, f)$ is $\min \left[n, \max\limits_{k_1,i_1}[m_{k_1,i_1}]+\max\limits_{k_2,i_2}[m_{k_2,i_2}]\right]$.
$\hfill\Box$

\subsection{Lipschitz Estimations}\label{sec: Lipschitz Estimations}

In this section, Theorem~\ref{thm:Lipschitz Estimations} is given to provide the lower and the upper bounds for Lipschitz Estimations.
\begin{theorem}\label{thm:Lipschitz Estimations}
Given an analytic function $f(z)$ within the domain for $|z| < R$, the first matrix $\bm{X}_1$ with the dimension $m$ and $K_1$ distinct eigenvalues $\lambda_{k_1}$ for $k_1=1,2,\ldots,K_1$ such that
\begin{eqnarray}\label{eq1-1:thm:Lipschitz Estimations}
\bm{X}_1&=&\sum\limits_{k_1=1}^{K_1}\sum\limits_{i_1=1}^{\alpha_{k_1}^{\mathrm{G}}} \lambda_{k_1} \bm{P}_{k_1,i_1}+
\sum\limits_{k_1=1}^{K_1}\sum\limits_{i_1=1}^{\alpha_{k_1}^{\mathrm{G}}} \bm{N}_{k_1,i_1},
\end{eqnarray}
where $\left\vert\lambda_{k_1}\right\vert<R$, and second matrix $\bm{X}_2$ with the dimension $m$ and $K_2$ distinct eigenvalues $\lambda_{k_2}$ for $k_2=1,2,\ldots,K_2$ such that
\begin{eqnarray}\label{eq1-2:thm:Lipschitz Estimations}
\bm{X}_2&=&\sum\limits_{k_2=1}^{K_2}\sum\limits_{i_2=1}^{\alpha_{k_2}^{\mathrm{G}}} \lambda_{k_2} \bm{P}_{k_2,i_2}+
\sum\limits_{k_2=1}^{K_2}\sum\limits_{i_2=1}^{\alpha_{k_2}^{\mathrm{G}}} \bm{N}_{k_2,i_2},
\end{eqnarray}
where $\left\vert\lambda_{k_2}\right\vert<R$. We also assume that $\lambda_{k_1} \neq \lambda_{k_2}$ for any $k_1$ and $k_2$. We define $f^{[1]}(\lambda_1, \lambda_2) \define \frac{f(\lambda_1) - f(\lambda_2)}{\lambda_1 - \lambda_2}$. 

Then, we have the following upper bound for Lipschitz estimation:
\begin{eqnarray}\label{eq2:thm:Lipschitz Estimations}
\lefteqn{\left\Vert f(\bm{X}_1) - f(\bm{X}_2) \right\Vert \leq}\nonumber \\
&& \left[\max\limits_{\lambda_1 \in \Lambda_{\bm{X}_1},\lambda_2 \in \Lambda_{\bm{X}_2}} \left\vert f^{[1]}(\lambda_1,\lambda_2)\right\vert\right]\left\Vert\bm{X}_1-\bm{X}_2\right\Vert \nonumber \\
&&+\sum\limits_{k_2=1}^{K_2}\sum\limits_{i_2=1}^{\alpha_{k_2}^{(\mathrm{G})}}\sum_{q_2=1}^{m_{k_2,i_2}-1}\left[\max\limits_{\lambda_1 \in \Lambda_{\bm{X}_1},\lambda_2 \in \Lambda_{\bm{X}_2}}\left\vert\frac{(f^{[1]})^{(-,q_2)}(\lambda_{1},\lambda_{2})}{q_2!}\right\vert\right]\left\Vert\bm{N}_{k_2,i_2}^{q_2}\right\Vert\left\Vert\bm{X}_1-\bm{X}_2\right\Vert\nonumber \\
&&+\sum\limits_{k_1=1}^{K_1}\sum\limits_{i_1=1}^{\alpha_{k_1}^{(\mathrm{G})}}\sum_{q_1=1}^{m_{k_1,i_1}-1}\left[\max\limits_{\lambda_1 \in \Lambda_{\bm{X}_1},\lambda_2 \in \Lambda_{\bm{X}_2}}\left\vert\frac{(f^{[1]})^{(q_1,-)}(\lambda_{1},\lambda_{2})}{q_1!}\right\vert\right]\left\Vert\bm{N}_{k_1,i_1}^{q_1}\right\Vert\left\Vert\bm{X}_1-\bm{X}_2\right\Vert  \nonumber \\
&&+ \sum\limits_{k_1=1}^{K_1}\sum\limits_{k_2=1}^{K_2}\sum\limits_{i_1=1}^{\alpha_{k_1}^{(\mathrm{G})}}\sum\limits_{i_2=1}^{\alpha_{k_2}^{(\mathrm{G})}}\sum_{q_1=1}^{m_{k_1,i_1}-1}\sum_{q_2=1}^{m_{k_2,i_2}-1}\left[\max\limits_{\lambda_1 \in \Lambda_{\bm{X}_1},\lambda_2 \in \Lambda_{\bm{X}_2}}\left\vert\frac{(f^{[1]})^{(q_1,q_2)}(\lambda_{1},\lambda_{2})}{q_1!q_2!}\right\vert\right]\nonumber\\
&&~~~\times\left\Vert\bm{N}_{k_1,i_1}^{q_1}\right\Vert\left\Vert\bm{N}_{k_2,i_2}^{q_2}\right\Vert\left\Vert\bm{X}_1-\bm{X}_2\right\Vert.
\end{eqnarray}

On the other hand, let us define the following matrices
\begin{eqnarray}\label{eq3:thm:Lipschitz Estimations}
\bm{A}_1&\define& \sum\limits_{k_1=1}^{K_1}\sum\limits_{k_2=1}^{K_2}\sum\limits_{i_1=1}^{\alpha_{k_1}^{(\mathrm{G})}}\sum\limits_{i_2=1}^{\alpha_{k_2}^{(\mathrm{G})}}f^{[1]}(\lambda_{k_1}, \lambda_{k_2})\bm{P}_{k_1,i_1}(\bm{X}_1-\bm{X}_2)\bm{P}_{k_2,i_2}\nonumber \\
\bm{A}_2&\define& \sum\limits_{k_1=1}^{K_1}\sum\limits_{k_2=1}^{K_2}\sum\limits_{i_1=1}^{\alpha_{k_1}^{(\mathrm{G})}}\sum\limits_{i_2=1}^{\alpha_{k_2}^{(\mathrm{G})}}\sum_{q_2=1}^{m_{k_2,i_2}-1}\frac{(f^{[1]})^{(-,q_2)}(\lambda_{k_1},\lambda_{k_2})}{q_2!}\bm{P}_{k_1,i_1}(\bm{X}_1-\bm{X}_2)\bm{N}_{k_2,i_2}^{q_2}\nonumber \\
\bm{A}_3&\define&\sum\limits_{k_1=1}^{K_1}\sum\limits_{k_2=1}^{K_2}\sum\limits_{i_1=1}^{\alpha_{k_1}^{(\mathrm{G})}}\sum\limits_{i_2=1}^{\alpha_{k_2}^{(\mathrm{G})}}\sum_{q_1=1}^{m_{k_1,i_1}-1}\frac{(f^{[1]})^{(q_1,-)}(\lambda_{k_1},\lambda_{k_2})}{q_1!}\bm{N}_{k_1,i_1}^{q_1}(\bm{X}_1-\bm{X}_2)\bm{P}_{k_2,i_2}\nonumber \\
\bm{A}_4&\define&\sum\limits_{k_1=1}^{K_1}\sum\limits_{k_2=1}^{K_2}\sum\limits_{i_1=1}^{\alpha_{k_1}^{(\mathrm{G})}}\sum\limits_{i_2=1}^{\alpha_{k_2}^{(\mathrm{G})}}\sum_{q_1=1}^{m_{k_1,i_1}-1}\sum_{q_2=1}^{m_{k_2,i_2}-1}\frac{(f^{[1]})^{(q_1,q_2)}(\lambda_{k_1},\lambda_{k_2})}{q_1!q_2!}\bm{N}_{k_1,i_1}^{q_1}(\bm{X}_1-\bm{X}_2)\bm{N}_{k_2,i_2}^{q_2}.
\end{eqnarray}
then, we have the lower bound for the Frobenius norm of $f(\bm{X}_1) - f(\bm{X}_2)$, which is given by
\begin{eqnarray}\label{eq4:thm:Lipschitz Estimations}
\left\Vert f(\bm{X}_1) - f(\bm{X}_2) \right\Vert\geq \max\left[0, \left\Vert\bm{A}_{\sigma(1)}\right\Vert - \sum\limits_{i=2}^4 \left\Vert\bm{A}_{\sigma(i)}\right\Vert\right],
\end{eqnarray}
where $\sigma$ is the permutation of matrices $\bm{A}_i$ for $i=1,2,3,4$ such that $\left\Vert\bm{A}_{\sigma(1)}\right\Vert \geq \left\Vert\bm{A}_{\sigma(2)}\right\Vert \geq\left\Vert\bm{A}_{\sigma(3)}\right\Vert \geq\left\Vert\bm{A}_{\sigma(4)}\right\Vert \geq$. 

Further, if we have $\left[\min\limits_{\lambda_1 \in \Lambda_{\bm{X}_1},\lambda_2 \in \Lambda_{\bm{X}_2}} f^{[1]}(\lambda_1,\lambda_2)\right]\left\Vert\bm{X}_1 - \bm{X}_2\right\Vert\geq \left\Vert\bm{A}_{2}\right\Vert+\left\Vert\bm{A}_{3}\right\Vert + \left\Vert\bm{A}_{4}\right\Vert$,  the lower bound for the Frobenius norm of $f(\bm{X}_1) - f(\bm{X}_2)$ can be expressed by
\begin{eqnarray}\label{eq5:thm:Lipschitz Estimations}
\left\Vert f(\bm{X}_1) - f(\bm{X}_2)\right\Vert\geq \left[\min\limits_{\lambda_1 \in \Lambda_{\bm{X}_1},\lambda_2 \in \Lambda_{\bm{X}_2}} f^{[1]}(\lambda_1,\lambda_2)\right]\left\Vert\bm{X}_1-\bm{X}_2\right\Vert - \sum\limits_{i=2}^4 \left\Vert\bm{A}_{i}\right\Vert.
\end{eqnarray}
s%
\end{theorem}
\textbf{Proof:}
From Corollary~\ref{cor:pert formula I},  we have
\begin{eqnarray}
f(\bm{X}_1)-f(\bm{X}_2)&=&T_{\frac{f(x_1) - f(x_2)}{x_1 -x_2}}^{\bm{X}_{1},\bm{X}_{2}}(\bm{X}_{1}-\bm{X}_{2}). 
\end{eqnarray}
This theorem is obtained by applying in $\beta = f^{[1]}$ and $\bm{Y}=\bm{X}_1 - \bm{X}_2$ in Theorem~\ref{thm: GDOI norm}.
$\hfill\Box$

\section{Continuity}\label{sec: Continuity}

In this section, we will prove the continuity property of the GDOI $T_{\beta}^{\bm{X}_1,\bm{X}_2}(\bm{Y})$. Follow the basic spirit of GDOI, we will define Generalized Triple Operator Integral (GTOI) in Eq.~\eqref{eq1:  GTOI def}. We have
\begin{eqnarray}\label{eq1:  GTOI def}
\lefteqn{T_{\beta}^{\bm{X}_1,\bm{X}_2,\bm{X}_3}(\bm{Y}_1, \bm{Y}_2)\define}\nonumber\\
&&\sum\limits_{k_1=1}^{K_1}\sum\limits_{k_2=1}^{K_2}\sum\limits_{k_3=1}^{K_3}\sum\limits_{i_1=1}^{\alpha_{k_1}^{(\mathrm{G})}}\sum\limits_{i_2=1}^{\alpha_{k_2}^{(\mathrm{G})}}\sum\limits_{i_3=1}^{\alpha_{k_3}^{(\mathrm{G})}}\beta(\lambda_{k_1}, \lambda_{k_2}, \lambda_{k_3})\bm{P}_{k_1,i_1}\bm{Y}_1\bm{P}_{k_2,i_2}\bm{Y}_2\bm{P}_{k_3,i_3} \nonumber \\
&&+\sum\limits_{k_1=1}^{K_1}\sum\limits_{k_2=1}^{K_2}\sum\limits_{k_3=1}^{K_3}\sum\limits_{i_1=1}^{\alpha_{k_1}^{(\mathrm{G})}}\sum\limits_{i_2=1}^{\alpha_{k_2}^{(\mathrm{G})}}\sum\limits_{i_3=1}^{\alpha_{k_3}^{(\mathrm{G})}}\sum_{q_3=1}^{m_{k_3,i_3}-1}\frac{\beta^{(-,-,q_3)}(\lambda_{k_1},\lambda_{k_2},\lambda_{k_3})}{q_3!}\bm{P}_{k_1,i_1}\bm{Y}_1\bm{P}_{k_2,i_2}\bm{Y}_2\bm{N}_{k_3,i_3}^{q_3}\nonumber \\
&&+\sum\limits_{k_1=1}^{K_1}\sum\limits_{k_2=1}^{K_2}\sum\limits_{k_3=1}^{K_3}\sum\limits_{i_1=1}^{\alpha_{k_1}^{(\mathrm{G})}}\sum\limits_{i_2=1}^{\alpha_{k_2}^{(\mathrm{G})}}\sum\limits_{i_3=1}^{\alpha_{k_3}^{(\mathrm{G})}}\sum_{q_2=1}^{m_{k_2,i_2}-1}\frac{\beta^{(-,q_2,-)}(\lambda_{k_1},\lambda_{k_2},\lambda_{k_3})}{q_2!}\bm{P}_{k_1,i_1}\bm{Y}_1\bm{N}_{k_2,i_2}^{q_2}\bm{Y}_2\bm{P}_{k_3,i_3} \nonumber \\
&&+\sum\limits_{k_1=1}^{K_1}\sum\limits_{k_2=1}^{K_2}\sum\limits_{k_3=1}^{K_3}\sum\limits_{i_1=1}^{\alpha_{k_1}^{(\mathrm{G})}}\sum\limits_{i_2=1}^{\alpha_{k_2}^{(\mathrm{G})}}\sum\limits_{i_3=1}^{\alpha_{k_3}^{(\mathrm{G})}}\sum_{q_1=1}^{m_{k_1,i_1}-1}\frac{\beta^{(q_1,-,-)}(\lambda_{k_1},\lambda_{k_2},\lambda_{k_3})}{q_1!}\bm{N}_{k_1,i_1}^{q_1}\bm{Y}_1\bm{P}_{k_2,i_2}\bm{Y}_2\bm{P}_{k_3,i_3} \nonumber \\
&&+\sum\limits_{k_1=1}^{K_1}\sum\limits_{k_2=1}^{K_2}\sum\limits_{k_3=1}^{K_3}\sum\limits_{i_1=1}^{\alpha_{k_1}^{(\mathrm{G})}}\sum\limits_{i_2=1}^{\alpha_{k_2}^{(\mathrm{G})}}\sum\limits_{i_3=1}^{\alpha_{k_3}^{(\mathrm{G})}}\sum_{q_2=1}^{m_{k_2,i_2}-1}\sum_{q_3=1}^{m_{k_3,i_3}-1}\frac{\beta^{(-,q_2,q_3)}(\lambda_{k_1},\lambda_{k_2},\lambda_{k_3})}{q_2! q_3!}\bm{P}_{k_1,i_1}\bm{Y}_1\bm{N}_{k_2,i_2}^{q_2}\bm{Y}_2\bm{N}_{k_3,i_3}^{q_3}\nonumber \\
&&+\sum\limits_{k_1=1}^{K_1}\sum\limits_{k_2=1}^{K_2}\sum\limits_{k_3=1}^{K_3}\sum\limits_{i_1=1}^{\alpha_{k_1}^{(\mathrm{G})}}\sum\limits_{i_2=1}^{\alpha_{k_2}^{(\mathrm{G})}}\sum\limits_{i_3=1}^{\alpha_{k_3}^{(\mathrm{G})}}\sum_{q_1=1}^{m_{k_1,i_1}-1}\sum_{q_3=1}^{m_{k_3,i_3}-1}\frac{\beta^{(q_1,-,q_3)}(\lambda_{k_1},\lambda_{k_2},\lambda_{k_3})}{q_1!q_3!}\bm{N}_{k_1,i_1}^{q_1}\bm{Y}_1\bm{P}_{k_2,i_2}\bm{Y}_2\bm{N}_{k_3,i_3}^{q_3} \nonumber \\
&&+\sum\limits_{k_1=1}^{K_1}\sum\limits_{k_2=1}^{K_2}\sum\limits_{k_3=1}^{K_3}\sum\limits_{i_1=1}^{\alpha_{k_1}^{(\mathrm{G})}}\sum\limits_{i_2=1}^{\alpha_{k_2}^{(\mathrm{G})}}\sum\limits_{i_3=1}^{\alpha_{k_3}^{(\mathrm{G})}}\sum_{q_1=1}^{m_{k_1,i_1}-1}\sum_{q_2=1}^{m_{k_2,i_2}-1}\frac{\beta^{(q_1,q_2,-)}(\lambda_{k_1},\lambda_{k_2},\lambda_{k_3})}{q_1!q_2!}\bm{N}_{k_1,i_1}^{q_1}\bm{Y}_1\bm{N}_{k_2,i_2}^{q_2}\bm{Y}_2\bm{P}_{k_3,i_3} \nonumber \\
&&+\sum\limits_{k_1=1}^{K_1}\sum\limits_{k_2=1}^{K_2}\sum\limits_{k_3=1}^{K_3}\sum\limits_{i_1=1}^{\alpha_{k_1}^{(\mathrm{G})}}\sum\limits_{i_2=1}^{\alpha_{k_2}^{(\mathrm{G})}}\sum\limits_{i_3=1}^{\alpha_{k_3}^{(\mathrm{G})}}\sum_{q_1=1}^{m_{k_1,i_1}-1}\sum_{q_2=1}^{m_{k_2,i_2}-1}\sum_{q_3=1}^{m_{k_3,i_3}-1}\frac{\beta^{(q_1,q_2,q_3)}(\lambda_{k_1},\lambda_{k_2},\lambda_{k_3})}{q_1!q_2!q_3!}\nonumber \\
&&~~\times\bm{N}_{k_1,i_1}^{q_1}\bm{Y}_1\bm{N}_{k_2,i_2}^{q_2}\bm{Y}_2\bm{N}_{k_3,i_3}^{q_3}.
\end{eqnarray}

From the definition of GTOI given by Eq.~\eqref{eq1:  GTOI def}, we have the following Theorem~\ref{thm: GTOI norm} about the upper and the lower bounds for GTOI.
\begin{theorem}\label{thm: GTOI norm}
We have the upper bound for the Frobenius norm of $T_{\beta}^{\bm{X}_1,\bm{X}_2,\bm{X}_3}(\bm{Y}_1,\bm{Y}_2)$, which is given by
\begin{eqnarray}\label{eq1: thm: GTOI norm}
\lefteqn{\left\Vert T_{\beta}^{\bm{X}_1,\bm{X}_2,\bm{X}_3}(\bm{Y}_1,\bm{Y}_2)\right\Vert\leq}\nonumber \\
&&\left[\max\limits_{\lambda_1 \in \Lambda_{\bm{X}_1},\lambda_2 \in \Lambda_{\bm{X}_2},\lambda_3 \in \Lambda_{\bm{X}_3}} \left\vert\beta(\lambda_1,\lambda_2,\lambda_3)\right\vert\right]\left\Vert\bm{Y}_1\right\Vert\left\Vert\bm{Y}_2\right\Vert\nonumber \\
&&+\sum\limits_{k_3=1}^{K_3}\sum\limits_{i_3=1}^{\alpha_{k_3}^{(\mathrm{G})}}\sum_{q_3=1}^{m_{k_3,i_3}-1}\left[\max\limits_{\lambda_1 \in \Lambda_{\bm{X}_1},\lambda_2 \in \Lambda_{\bm{X}_2}\lambda_3 \in \Lambda_{\bm{X}_3}}\left\vert\frac{\beta^{(-,-,q_3)}(\lambda_{1},\lambda_{2},\lambda_3)}{q_3!}\right\vert\right]\left\Vert\bm{Y}_1\right\Vert\left\Vert\bm{Y}_2\right\Vert\left\Vert\bm{N}_{k_3,i_3}^{q_3}\right\Vert\nonumber \\
&&+\sum\limits_{k_2=1}^{K_2}\sum\limits_{i_2=1}^{\alpha_{k_2}^{(\mathrm{G})}}\sum_{q_2=1}^{m_{k_2,i_2}-1}\left[\max\limits_{\lambda_1 \in \Lambda_{\bm{X}_1},\lambda_2 \in \Lambda_{\bm{X}_2}\lambda_3 \in \Lambda_{\bm{X}_3}}\left\vert\frac{\beta^{(-,q_2,-)}(\lambda_{1},\lambda_{2},\lambda_3)}{q_2!}\right\vert\right]\left\Vert\bm{Y}_1\right\Vert\left\Vert\bm{N}_{k_2,i_2}^{q_2}\right\Vert\left\Vert\bm{Y}_2\right\Vert\nonumber \\
&&+\sum\limits_{k_1=1}^{K_1}\sum\limits_{i_1=1}^{\alpha_{k_1}^{(\mathrm{G})}}\sum_{q_1=1}^{m_{k_1,i_1}-1}\left[\max\limits_{\lambda_1 \in \Lambda_{\bm{X}_1},\lambda_2 \in \Lambda_{\bm{X}_2}\lambda_3 \in \Lambda_{\bm{X}_3}}\left\vert\frac{\beta^{(q_1,-,-)}(\lambda_{1},\lambda_{2},\lambda_3)}{q_1!}\right\vert\right]\left\Vert\bm{N}_{k_1,i_1}^{q_1}\right\Vert\left\Vert\bm{Y}_1\right\Vert\left\Vert\bm{Y}_2\right\Vert \nonumber \\
&&+\sum\limits_{k_2=1}^{K_2}\sum\limits_{k_3=1}^{K_3}\sum\limits_{i_2=1}^{\alpha_{k_2}^{(\mathrm{G})}}\sum\limits_{i_3=1}^{\alpha_{k_3}^{(\mathrm{G})}}\sum_{q_2=1}^{m_{k_2,i_2}-1}\sum_{q_3=1}^{m_{k_3,i_3}-1}\left[\max\limits_{\lambda_1 \in \Lambda_{\bm{X}_1},\lambda_2 \in \Lambda_{\bm{X}_2},\lambda_3 \in \Lambda_{\bm{X}_3}}\left\vert\frac{\beta^{(-,q_2,q_3)}(\lambda_{1},\lambda_{2},\lambda_3)}{q_2!q_3!}\right\vert\right]\nonumber\\
&&~~~\times\left\Vert\bm{Y}_1\right\Vert\left\Vert\bm{N}_{k_2,i_2}^{q_2}\right\Vert\left\Vert\bm{Y}_2\right\Vert\left\Vert\bm{N}_{k_3,i_3}^{q_3}\right\Vert\nonumber \\
&&+ \sum\limits_{k_1=1}^{K_1}\sum\limits_{k_3=1}^{K_3}\sum\limits_{i_1=1}^{\alpha_{k_1}^{(\mathrm{G})}}\sum\limits_{i_3=1}^{\alpha_{k_3}^{(\mathrm{G})}}\sum_{q_1=1}^{m_{k_1,i_1}-1}\sum_{q_3=1}^{m_{k_3,i_3}-1}\left[\max\limits_{\lambda_1 \in \Lambda_{\bm{X}_1},\lambda_2 \in \Lambda_{\bm{X}_2},\lambda_3 \in \Lambda_{\bm{X}_3}}\left\vert\frac{\beta^{(q_1,-,q_3)}(\lambda_{1},\lambda_{2},\lambda_3)}{q_1!q_3!}\right\vert\right]\nonumber\\
&&~~~\times\left\Vert\bm{N}_{k_1,i_1}^{q_1}\right\Vert\left\Vert\bm{Y}_1\right\Vert\left\Vert\bm{Y}_2\right\Vert\left\Vert\bm{N}_{k_3,i_3}^{q_3}\right\Vert \nonumber \\
&&+\sum\limits_{k_1=1}^{K_1}\sum\limits_{k_2=1}^{K_2}\sum\limits_{i_1=1}^{\alpha_{k_1}^{(\mathrm{G})}}\sum\limits_{i_2=1}^{\alpha_{k_2}^{(\mathrm{G})}}\sum_{q_1=1}^{m_{k_1,i_1}-1}\sum_{q_2=1}^{m_{k_2,i_2}-1}\left[\max\limits_{\lambda_1 \in \Lambda_{\bm{X}_1},\lambda_2 \in \Lambda_{\bm{X}_2},\lambda_3 \in \Lambda_{\bm{X}_3}}\left\vert\frac{\beta^{(q_1,q_2,-)}(\lambda_{1},\lambda_{2},\lambda_3)}{q_1!q_2!}\right\vert\right]\nonumber\\
&&~~~\times\left\Vert\bm{N}_{k_1,i_1}^{q_1}\right\Vert\left\Vert\bm{Y}_1\right\Vert\left\Vert\bm{N}_{k_2,i_2}^{q_2}\right\Vert\left\Vert\bm{Y}_2\right\Vert \nonumber \\
&&+\sum\limits_{k_1=1}^{K_1}\sum\limits_{k_2=1}^{K_2}\sum\limits_{k_3=1}^{K_3}\sum\limits_{i_1=1}^{\alpha_{k_1}^{(\mathrm{G})}}\sum\limits_{i_2=1}^{\alpha_{k_2}^{(\mathrm{G})}}\sum\limits_{i_3=1}^{\alpha_{k_3}^{(\mathrm{G})}}\sum_{q_1=1}^{m_{k_1,i_1}-1}\sum_{q_2=1}^{m_{k_2,i_2}-1}\sum_{q_3=1}^{m_{k_3,i_3}-1}\nonumber \\
&&\left[\max\limits_{\lambda_1 \in \Lambda_{\bm{X}_1},\lambda_2 \in \Lambda_{\bm{X}_2},\lambda_3 \in \Lambda_{\bm{X}_3}}\left\vert\frac{\beta^{(q_1,q_2,q_3)}(\lambda_{1},\lambda_{2},\lambda_3)}{q_1!q_2!q_3!}\right\vert\right]\left\Vert\bm{N}_{k_1,i_1}^{q_1}\right\Vert\left\Vert\bm{Y}_1\right\Vert\left\Vert\bm{N}_{k_2,i_2}^{q_2}\right\Vert\left\Vert\bm{Y}_2\right\Vert\left\Vert\bm{N}_{k_3,i_3}^{q_3}\right\Vert
\end{eqnarray}
where $\Lambda_{\bm{X}_1}, \Lambda_{\bm{X}_2}$ and $\Lambda_{\bm{X}_3}$ are spectrums of the matrix $\bm{X}_1, \bm{X}_2$ and the matrix $\bm{X}_3$, respectively.

On the other hand, let us define the following matrices
\begin{eqnarray}\label{eq1-1: thm: GTOI norm}
\bm{A}_1&\define&\sum\limits_{k_1=1}^{K_1}\sum\limits_{k_2=1}^{K_2}\sum\limits_{k_3=1}^{K_3}\sum\limits_{i_1=1}^{\alpha_{k_1}^{(\mathrm{G})}}\sum\limits_{i_2=1}^{\alpha_{k_2}^{(\mathrm{G})}}\sum\limits_{i_3=1}^{\alpha_{k_3}^{(\mathrm{G})}}\beta(\lambda_{k_1}, \lambda_{k_2}, \lambda_{k_3})\bm{P}_{k_1,i_1}\bm{Y}_1\bm{P}_{k_2,i_2}\bm{Y}_2\bm{P}_{k_3,i_3}\nonumber \\
\bm{A}_2&\define&\sum\limits_{k_1=1}^{K_1}\sum\limits_{k_2=1}^{K_2}\sum\limits_{k_3=1}^{K_3}\sum\limits_{i_1=1}^{\alpha_{k_1}^{(\mathrm{G})}}\sum\limits_{i_2=1}^{\alpha_{k_2}^{(\mathrm{G})}}\sum\limits_{i_3=1}^{\alpha_{k_3}^{(\mathrm{G})}}\sum_{q_3=1}^{m_{k_3,i_3}-1}\frac{\beta^{(-,-,q_3)}(\lambda_{k_1},\lambda_{k_2},\lambda_{k_3})}{q_3!}\bm{P}_{k_1,i_1}\bm{Y}_1\bm{P}_{k_2,i_2}\bm{Y}_2\bm{N}_{k_3,i_3}^{q_3}\nonumber \\
\bm{A}_3&\define&\sum\limits_{k_1=1}^{K_1}\sum\limits_{k_2=1}^{K_2}\sum\limits_{k_3=1}^{K_3}\sum\limits_{i_1=1}^{\alpha_{k_1}^{(\mathrm{G})}}\sum\limits_{i_2=1}^{\alpha_{k_2}^{(\mathrm{G})}}\sum\limits_{i_3=1}^{\alpha_{k_3}^{(\mathrm{G})}}\sum_{q_2=1}^{m_{k_2,i_2}-1}\frac{\beta^{(-,q_2,-)}(\lambda_{k_1},\lambda_{k_2},\lambda_{k_3})}{q_2!}\bm{P}_{k_1,i_1}\bm{Y}_1\bm{N}_{k_2,i_2}^{q_2}\bm{Y}_2\bm{P}_{k_3,i_3}\nonumber \\
\bm{A}_4&\define&\sum\limits_{k_1=1}^{K_1}\sum\limits_{k_2=1}^{K_2}\sum\limits_{k_3=1}^{K_3}\sum\limits_{i_1=1}^{\alpha_{k_1}^{(\mathrm{G})}}\sum\limits_{i_2=1}^{\alpha_{k_2}^{(\mathrm{G})}}\sum\limits_{i_3=1}^{\alpha_{k_3}^{(\mathrm{G})}}\sum_{q_1=1}^{m_{k_1,i_1}-1}\frac{\beta^{(q_1,-,-)}(\lambda_{k_1},\lambda_{k_2},\lambda_{k_3})}{q_1!}\bm{N}_{k_1,i_1}^{q_1}\bm{Y}_1\bm{P}_{k_2,i_2}\bm{Y}_2\bm{P}_{k_3,i_3}\nonumber \\
\bm{A}_5&\define&\sum\limits_{k_1=1}^{K_1}\sum\limits_{k_2=1}^{K_2}\sum\limits_{k_3=1}^{K_3}\sum\limits_{i_1=1}^{\alpha_{k_1}^{(\mathrm{G})}}\sum\limits_{i_2=1}^{\alpha_{k_2}^{(\mathrm{G})}}\sum\limits_{i_3=1}^{\alpha_{k_3}^{(\mathrm{G})}}\sum_{q_2=1}^{m_{k_2,i_2}-1}\sum_{q_3=1}^{m_{k_3,i_3}-1}\frac{\beta^{(-,q_2,q_3)}(\lambda_{k_1},\lambda_{k_2},\lambda_{k_3})}{q_2! q_3!}\bm{P}_{k_1,i_1}\bm{Y}_1\bm{N}_{k_2,i_2}^{q_2}\bm{Y}_2\bm{N}_{k_3,i_3}^{q_3}\nonumber \\
\bm{A}_6&\define&\sum\limits_{k_1=1}^{K_1}\sum\limits_{k_2=1}^{K_2}\sum\limits_{k_3=1}^{K_3}\sum\limits_{i_1=1}^{\alpha_{k_1}^{(\mathrm{G})}}\sum\limits_{i_2=1}^{\alpha_{k_2}^{(\mathrm{G})}}\sum\limits_{i_3=1}^{\alpha_{k_3}^{(\mathrm{G})}}\sum_{q_1=1}^{m_{k_1,i_1}-1}\sum_{q_3=1}^{m_{k_3,i_3}-1}\frac{\beta^{(q_1,-,q_3)}(\lambda_{k_1},\lambda_{k_2},\lambda_{k_3})}{q_1!q_3!}\bm{N}_{k_1,i_1}^{q_1}\bm{Y}_1\bm{P}_{k_2,i_2}\bm{Y}_2\bm{N}_{k_3,i_3}^{q_3} \nonumber \\
\bm{A}_7&\define&\sum\limits_{k_1=1}^{K_1}\sum\limits_{k_2=1}^{K_2}\sum\limits_{k_3=1}^{K_3}\sum\limits_{i_1=1}^{\alpha_{k_1}^{(\mathrm{G})}}\sum\limits_{i_2=1}^{\alpha_{k_2}^{(\mathrm{G})}}\sum\limits_{i_3=1}^{\alpha_{k_3}^{(\mathrm{G})}}\sum_{q_1=1}^{m_{k_1,i_1}-1}\sum_{q_2=1}^{m_{k_2,i_2}-1}\frac{\beta^{(q_1,q_2,-)}(\lambda_{k_1},\lambda_{k_2},\lambda_{k_3})}{q_1!q_2!}\bm{N}_{k_1,i_1}^{q_1}\bm{Y}_1\bm{N}_{k_2,i_2}^{q_2}\bm{Y}_2\bm{P}_{k_3,i_3} \nonumber \\
\bm{A}_8&\define&\sum\limits_{k_1=1}^{K_1}\sum\limits_{k_2=1}^{K_2}\sum\limits_{k_3=1}^{K_3}\sum\limits_{i_1=1}^{\alpha_{k_1}^{(\mathrm{G})}}\sum\limits_{i_2=1}^{\alpha_{k_2}^{(\mathrm{G})}}\sum\limits_{i_3=1}^{\alpha_{k_3}^{(\mathrm{G})}}\sum_{q_1=1}^{m_{k_1,i_1}-1}\sum_{q_2=1}^{m_{k_2,i_2}-1}\sum_{q_3=1}^{m_{k_3,i3}-1}\frac{\beta^{(q_1,q_2,q_3)}(\lambda_{k_1},\lambda_{k_2},\lambda_{k_3})}{q_1!q_2!q_3!}\nonumber \\
&&~~\times\bm{N}_{k_1,i_1}^{q_1}\bm{Y}_1\bm{N}_{k_2,i_2}^{q_2}\bm{Y}_2\bm{N}_{k_3,i_3}^{q_3},
\end{eqnarray}
then, we have the lower bound for the Frobenius norm of $T_{\beta}^{\bm{X}_1,\bm{X}_2,\bm{X}_3}(\bm{Y}_1,\bm{Y}_2)$, which is given by
\begin{eqnarray}\label{eq2: thm: GTOI norm}
\left\Vert T_{\beta}^{\bm{X}_1,\bm{X}_2,\bm{X}_3}(\bm{Y}_1,\bm{Y}_2) \right\Vert\geq \max\left[0, \left\Vert\bm{A}_{\sigma(1)}\right\Vert - \sum\limits_{i=2}^8 \left\Vert\bm{A}_{\sigma(i)}\right\Vert\right],
\end{eqnarray}
where $\sigma$ is the permutation of matrices $\bm{A}_i$ for $i=1,2,3,4,5,6,7,8$ such that $\left\Vert\bm{A}_{\sigma(1)}\right\Vert \geq \ldots \geq \left\Vert\bm{A}_{\sigma(8)}\right\Vert \geq$. 

Further, if we have $\left[\min\limits_{\lambda_1 \in \Lambda_{\bm{X}_1},\lambda_2 \in \Lambda_{\bm{X}_2},\lambda_3 \in \Lambda_{\bm{X}_3}} \beta(\lambda_1,\lambda_2,\lambda_3)\right]\left\Vert\bm{Y}_1\right\Vert \left\Vert\bm{Y}_2\right\Vert \geq \sum\limits_{i=2}^8\left\Vert\bm{A}_{i}\right\Vert$,  the lower bound for the Frobenius norm of $T_{\beta}^{\bm{X}_1,\bm{X}_2,\bm{X}_3}(\bm{Y}_1,\bm{Y}_2)$ can be expressed by
\begin{eqnarray}\label{eq3: thm: GTOI norm}
\left\Vert T_{\beta}^{\bm{X}_1,\bm{X}_2,\bm{X}_3}(\bm{Y}_1,\bm{Y}_2) \right\Vert\geq \left[\min\limits_{\lambda_1 \in \Lambda_{\bm{X}_1},\lambda_2 \in \Lambda_{\bm{X}_2},\lambda_3 \in \Lambda_{\bm{X}_3}} \beta(\lambda_1,\lambda_2,\lambda_3)\right]\left\Vert\bm{Y}_1\right\Vert\left\Vert\bm{Y}_2\right\Vert - \sum\limits_{i=2}^8 \left\Vert\bm{A}_{i}\right\Vert.
\end{eqnarray}
\end{theorem}
\textbf{Proof:}
The proof is similar to the proof in Theorem~\ref{thm: GDOI norm}.
$\hfill\Box$

Below, we will show the telescope property for GTOI similar to conventional MOI that only consider Hermitian or self-adjoint parameter matrices~\cite{skripka2019multilinear}.
\begin{lemma}\label{lma: telescope propt}
We have
\begin{eqnarray}\label{eq1:lma: telescope propt}
T_{f^{[1]}}^{\bm{A},\bm{X}}(\bm{Y}) -  T_{f^{[1]}}^{\bm{B},\bm{X}}(\bm{Y})&=&
T_{f^{[2]}}^{\bm{A},\bm{B},\bm{X}}(\bm{A}-\bm{B},\bm{Y}),
\end{eqnarray}
where $f^{[1]}(x_0,x_1)$ and $^{[2]}(x_0,x_1,x_2)$ are first and second divide differences. 
\end{lemma}
\textbf{Proof:}
We have 
\begin{eqnarray}\label{eq2:lma: telescope propt}
T_{f^{[2]}}^{\bm{A},\bm{B},\bm{X}}(\bm{A}-\bm{B},\bm{Y})&=&T_{f^{[2]}}^{\bm{A},\bm{B},\bm{X}}(\bm{A},\bm{Y}) - T_{f^{[2]}}^{\bm{A},\bm{B},\bm{X}}(\bm{B},\bm{Y})\nonumber \\
&=&T_{x_0 f^{[2]}}^{\bm{A},\bm{B},\bm{X}}(\bm{I},\bm{Y}) - T_{x_1 f^{[2]}}^{\bm{A},\bm{B},\bm{X}}(\bm{I},\bm{Y})\nonumber \\
&=&T_{x_0 f^{[2]}-x_1 f^{[2]}}^{\bm{A},\bm{B},\bm{X}}(\bm{I},\bm{Y})\nonumber \\
&=_1&T_{f^{[1]}(x_0,x_2)}^{\bm{A},\bm{X}}(\bm{Y}) - T_{f^{[1]}(x_1,x_2)}^{\bm{B},\bm{X}}(\bm{Y})\nonumber \\
&=&T_{f^{[1]}}^{\bm{A},\bm{X}}(\bm{Y}) -  T_{f^{[1]}}^{\bm{B},\bm{X}}(\bm{Y}),
\end{eqnarray}
where we have $x_0 f^{[2]}-x_1 f^{[2]}=f^{[1]}(x_0,x_2)-f^{[1]}(x_1,x_2)$ in $=_1$.
$\hfill\Box$

The following Theorem~\ref{thm:GDOI continuity dd} will show the continuity property of the GDOI $T_{f^{[1]}}^{\bm{X}_1,\bm{X}_2}(\bm{Y})$. 
\begin{theorem}\label{thm:GDOI continuity dd}
Given two sequence of matrices $\bm{X}_{1,\ell_1}$ and $\bm{X}_{2,\ell_2}$ satisfying $\bm{X}_{1,\ell_1} \rightarrow \bm{X}_{1}$ and $\bm{X}_{2,\ell_2} \rightarrow \bm{X}_{2}$, respectively, then, we have 
\begin{eqnarray}\label{eq1:thm:GDOI continuity dd}
T_{f^{[1]}}^{\bm{X}_{1,\ell_1},\bm{X}_{2,\ell_2}}(\bm{Y}) \rightarrow T_{f^{[1]}}^{\bm{X}_{1},\bm{X}_{2}}(\bm{Y}),
\end{eqnarray}
where $\rightarrow$ is in the sense of Frobenius norm, i.e., 
\begin{eqnarray}\label{eq2:thm:GDOI continuity dd}
\lim\limits_{\ell_1, \ell_2 \rightarrow \infty}\left\Vert T_{f^{[1]}}^{\bm{X}_{1,\ell_1},\bm{X}_{2,\ell_2}}(\bm{Y}) - T_{f^{[1]}}^{\bm{X}_{1},\bm{X}_{2}}(\bm{Y})\right\Vert=0
\end{eqnarray}
\end{theorem}
\textbf{Proof:}
Because we have
\begin{eqnarray}\label{eq3:thm:GDOI continuity dd}
\lefteqn{\left\Vert T_{f^{[1]}}^{\bm{X}_{1,\ell_1},\bm{X}_{2,\ell_2}}(\bm{Y}) - T_{f^{[1]}}^{\bm{X}_{1},\bm{X}_{2}}(\bm{Y})\right\Vert}\nonumber \\
&\leq&\left\Vert T_{f^{[1]}}^{\bm{X}_{1,\ell_1},\bm{X}_{2,\ell_2}}(\bm{Y}) - T_{f^{[1]}}^{\bm{X}_{1},\bm{X}_{2,\ell_2}}(\bm{Y})\right\Vert + \left\Vert T_{f^{[1]}}^{\bm{X}_{1},\bm{X}_{2,\ell_2}}(\bm{Y}) - T_{f^{[1]}}^{\bm{X}_{1},\bm{X}_{2}}(\bm{Y})\right\Vert\nonumber \\
&=_1&\left\Vert T_{f^{[2]}}^{\bm{X}_{1,\ell_1},\bm{X}_{1},\bm{X}_{2,\ell_2}}(\bm{X}_{1,\ell_1}-\bm{X}_1,\bm{Y})\right\Vert + \left\Vert T_{f^{[2]}}^{\bm{X}_{1},\bm{X}_{2,\ell_2},\bm{X}_2}(\bm{Y}, \bm{X}_{2,\ell_2}-\bm{X}_2)\right\Vert\nonumber \\
&\leq_2& \epsilon/2 +\epsilon/2 = \epsilon
\end{eqnarray}
where we apply Lemma~\ref{lma: telescope propt} in $=_1$, and both $\bm{X}_{1,\ell_1}, \bm{X}_{2,\ell_2}$ converge to $\bm{X}_1, \bm{X}_2$ with Theorem~\ref{thm: GTOI norm} used to establish $\leq_2$.
$\hfill\Box$

Theorem~\ref{thm:GDOI continuity dd} is srestrictive to the underlying function as divide difference only. We wish to have continuity property for more genereal function instead just $f^{[1]}$. We have to prepare   the following Lemma~\ref{lma:continuity ii and iii} about the continuity of the variable matrix $\bm{Y}$ and the underlying function $\beta$ in GDOI $T_{\beta}^{\bm{X}_1, \bm{X}_2}(\bm{Y})$.
\begin{lemma}\label{lma:continuity ii and iii}
(i) Given a sequence of matrix $\bm{Y}_\ell$ such that $\bm{Y}_\ell \rightarrow \bm{Y}$, we assume that 
those terms involving $\max\limits_{\lambda_1 \in \Lambda_{\bm{X}_1},\lambda_2 \in \Lambda_{\bm{X}_2}}$ in Eq.~\eqref{eq1: thm: GDOI norm} are finite. Then, we have
\begin{eqnarray}\label{eq1:lma:continuity ii and iii}
T_{\beta}^{\bm{X}_1, \bm{X}_2}(\bm{Y}_\ell)&\rightarrow&T_{\beta}^{\bm{X}_1, \bm{X}_2}(\bm{Y}).
\end{eqnarray}

(ii) Given a sequence of function $\beta_\ell$ such that supremum norm $\left\Vert \beta_\ell - \beta \right\Vert \rightarrow 0$, we have
\begin{eqnarray}\label{eq2:lma:continuity ii and iii}
T_{\beta_\ell}^{\bm{X}_1, \bm{X}_2}(\bm{Y})&\rightarrow&T_{\beta}^{\bm{X}_1, \bm{X}_2}(\bm{Y}).
\end{eqnarray}
\end{lemma}
\textbf{Proof:}
Since we have
\begin{eqnarray}\label{eq3:lma:continuity ii and iii}
\left\Vert T_{\beta}^{\bm{X}_1, \bm{X}_2}(\bm{Y}_\ell) - T_{\beta}^{\bm{X}_1, \bm{X}_2}(\bm{Y}) \right\Vert
&=&\left\Vert T_{\beta}^{\bm{X}_1, \bm{X}_2}(\bm{Y}_\ell - \bm{Y}) \right\Vert \nonumber \\
&\leq& \epsilon,
\end{eqnarray}
where $\epsilon$ is any positive number and this inequality comes from Theorem~\ref{thm: GDOI norm}and $\bm{Y}_\ell \rightarrow \bm{Y}$. This proves Part (i).

For Part (ii), we have
\begin{eqnarray}\label{eq4:lma:continuity ii and iii}
\left\Vert T_{\beta_\ell}^{\bm{X}_1, \bm{X}_2}(\bm{Y}) - T_{\beta}^{\bm{X}_1, \bm{X}_2}(\bm{Y}) \right\Vert
&=&\left\Vert T_{\beta_\ell-\beta}^{\bm{X}_1, \bm{X}_2}( \bm{Y}) \right\Vert \nonumber \\
&\leq& \epsilon,
\end{eqnarray}
where $\epsilon$ is any positive number and this inequality comes from Theorem~\ref{thm: GDOI norm}and $\left\Vert \beta_\ell - \beta \right\Vert \rightarrow 0$. 
$\hfill\Box$

\begin{theorem}\label{thm:GDOI continuity}
Given two sequence of matrices $\bm{X}_{1,\ell_1}$ and $\bm{X}_{2,\ell_2}$ satisfying $\bm{X}_{1,\ell_1} \rightarrow \bm{X}_{1}$ and $\bm{X}_{2,\ell_2} \rightarrow \bm{X}_{2}$, respectively, then, we have 
\begin{eqnarray}\label{eq1:thm:GDOI continuity}
T_{\beta}^{\bm{X}_{1,\ell_1},\bm{X}_{2,\ell_2}}(\bm{Y}) \rightarrow T_{\beta}^{\bm{X}_{1},\bm{X}_{2}}(\bm{Y}),
\end{eqnarray}
where $\rightarrow$ is in the sense of Frobenius norm, i.e., 
\begin{eqnarray}\label{eq2:thm:GDOI continuity}
\lim\limits_{\ell_1, \ell_2 \rightarrow \infty}\left\Vert T_{\beta}^{\bm{X}_{1,\ell_1},\bm{X}_{2,\ell_2}}(\bm{Y}) - T_{\beta}^{\bm{X}_{1},\bm{X}_{2}}(\bm{Y})\right\Vert=0.
\end{eqnarray}
\end{theorem}
\textbf{Proof:}
For any bivariate polynomial function of degree $\ell$, we can exactly represent it using a linear combination of the divided differences of monomials and the second variable. We have
\begin{eqnarray}\label{eq2-1:thm:GDOI continuity}
x_1^{k_1}x_2^{k_2}&=&\left[f^{[1]}_{k_1 + 1}(x_1,x_2) - x_2 f^{[1]}_{k_1}(x_1,x_2)\right]x_2^{k_2}\nonumber \\
&=&f^{[1]}_{k_1 + 1}(x_1,x_2)x_2^{k_2}-  f^{[1]}_{k_1}(x_1,x_2)x_2^{k_2 + 1},
\end{eqnarray}
where $f^{[1]}_{k_1}(x_1,x_2) \define \frac{x_1^{k_1} - x_2^{k_2}}{x_1 - x_2}$. Then, we have
\begin{eqnarray}\label{eq2-2:thm:GDOI continuity}
\sum\limits_{k_1=0,k_2=0}^{k_1+k_2=\ell} c_{k_1,k_2}x_1^{k_1}x_2^{k_2}&=&\sum\limits_{k_1=0,k_2=0}^{k_1+k_2=\ell}c_{k_1,k_2}\left[f^{[1]}_{k_1 + 1}(x_1,x_2)x_2^{k_2}-  f^{[1]}_{k_1}(x_1,x_2)x_2^{k_2 + 1}\right]\nonumber \\
&=&\sum\limits_{k'_1=1,k'_2=0}^{k'_1+k'_2=\ell+1} d_{k'_1,k'_2}f^{[1]}_{k'_1 }(x_1,x_2)x_2^{k'_2},
\end{eqnarray}
which shows that any bivariate polynomial can be expressed as linear combination of divided difference and its product with the variable $x_2$.

By Weierstrass approximation theorem, let us use $\tilde{\beta}\define\sum\limits_{k_1,k_2} c_{k_1,k_2}x_1^{k_1}x_2^{k_2}$ as a bivariate polynomial function to approximate $\beta$, i.e., $\tilde{\beta}\rightarrow\beta$. Therefore, we have
\begin{eqnarray}\label{eq3:thm:GDOI continuity}
\tilde{\beta}&=&\sum\limits_{k'_1=1,k'_2=0}^{k'_1+k'_2=\ell+1} d_{k'_1,k'_2}f^{[1]}_{k'_1 }(x_1,x_2)x_2^{k'_2},
\end{eqnarray}
where $d_{k'_1,k'_2}$ are complex scalars.

Then, we have
\begin{eqnarray}\label{eq4:thm:GDOI continuity}
\left\Vert T_{\beta}^{\bm{X}_{1,\ell_1},\bm{X}_{2,\ell_2}}(\bm{Y}) - T_{\beta}^{\bm{X}_{1},\bm{X}_{2}}(\bm{Y})\right\Vert&\leq&\underbrace{\left\Vert T_{\beta}^{\bm{X}_{1,\ell_1},\bm{X}_{2,\ell_2}}(\bm{Y}) - T_{\tilde{\beta}}^{\bm{X}_{1,\ell_1},\bm{X}_{2,\ell_2}}(\bm{Y})\right\Vert}_{\mbox{Part I}}\nonumber \\
&&+\underbrace{\left\Vert T_{\tilde{\beta}}^{\bm{X}_{1,\ell_1},\bm{X}_{2,\ell_2}}(\bm{Y}) - T_{\tilde{\beta}}^{\bm{X}_{1},\bm{X}_{2}}(\bm{Y})\right\Vert}_{\mbox{Part II}}\nonumber \\
&&+\underbrace{\left\Vert T_{\tilde{\beta}}^{\bm{X}_{1},\bm{X}_{2}}(\bm{Y}) - T_{\beta}^{\bm{X}_{1},\bm{X}_{2}}(\bm{Y})\right\Vert}_{\mbox{Part III}}.
\end{eqnarray}

Let us analyze each term in the R.H.S. of Eq.~\eqref{eq4:thm:GDOI continuity}. For Part I, we have
\begin{eqnarray}\label{eq5:thm:GDOI continuity}
\left\Vert T_{\beta}^{\bm{X}_{1,\ell_1},\bm{X}_{2,\ell_2}}(\bm{Y}) - T_{\tilde{\beta}}^{\bm{X}_{1,\ell_1},\bm{X}_{2,\ell_2}}(\bm{Y})\right\Vert \leq \epsilon/3,
\end{eqnarray}
by part (ii) in Lemma~\ref{lma:continuity ii and iii}.  For Part III, we also have
\begin{eqnarray}\label{eq6:thm:GDOI continuity}
\left\Vert T_{\tilde{\beta}}^{\bm{X}_{1},\bm{X}_{2}}(\bm{Y}) - T_{\beta}^{\bm{X}_{1},\bm{X}_{2}}(\bm{Y})\right\Vert\leq \epsilon/3,
\end{eqnarray}
by part (ii) in Lemma~\ref{lma:continuity ii and iii} again. For Part II, we have
\begin{eqnarray}\label{eq7:thm:GDOI continuity}
\lefteqn{\left\Vert T_{\tilde{\beta}}^{\bm{X}_{1,\ell_1},\bm{X}_{2,\ell_2}}(\bm{Y}) - T_{\tilde{\beta}}^{\bm{X}_{1},\bm{X}_{2}}(\bm{Y})\right\Vert}\nonumber \\
&=&\left\Vert T_{\sum\limits_{k'_1=1,k'_2=0}^{k'_1+k'_2=\ell+1} d_{k'_1,k'_2}f^{[1]}_{k'_1 }(x_1,x_2)x_2^{k'_2}}^{\bm{X}_{1,\ell_1},\bm{X}_{2,\ell_2}}(\bm{Y}) - T_{\sum\limits_{k'_1=1,k'_2=0}^{k'_1+k'_2=\ell+1} d_{k'_1,k'_2}f^{[1]}_{k'_1 }(x_1,x_2)x_2^{k'_2}}^{\bm{X}_{1},\bm{X}_{2}}(\bm{Y})\right\Vert\nonumber \\
&=_1&\left\Vert \sum\limits_{k'_1=1,k'_2=0}^{k'_1+k'_2=\ell+1} d_{k'_1,k'_2} T_{ f^{[1]}_{k'_1}(x_1,x_2)}^{\bm{X}_{1,\ell_1},\bm{X}_{2,\ell_2}}(T_{x_2^{k'_2}}^{\bm{X}_{1,\ell_1},\bm{X}_{2,\ell_2}}(\bm{Y})) -\sum\limits_{k'_1=1,k'_2=0}^{k'_1+k'_2=\ell+1} d_{k'_1,k'_2} T_{f^{[1]}_{k'_1}(x_1,x_2)}^{\bm{X}_{1},\bm{X}_{2}}(T_{x_2^{k'_2}}^{\bm{X}_{1},\bm{X}_{2}}(\bm{Y}))\right\Vert\nonumber \\
&\leq&\left\Vert \sum\limits_{k'_1=1,k'_2=0}^{k'_1+k'_2=\ell+1} d_{k'_1,k'_2} T_{ f^{[1]}_{k'_1}(x_1,x_2)}^{\bm{X}_{1,\ell_1},\bm{X}_{2,\ell_2}}(T_{x_2^{k'_2}}^{\bm{X}_{1,\ell_1},\bm{X}_{2,\ell_2}}(\bm{Y})) -\sum\limits_{k'_1=1,k'_2=0}^{k'_1+k'_2=\ell+1} d_{k'_1,k'_2} T_{f^{[1]}_{k'_1}(x_1,x_2)}^{\bm{X}_{1,\ell_1},\bm{X}_{2,\ell_2}}(T_{x_2^{k'_2}}^{\bm{X}_{1},\bm{X}_{2}}(\bm{Y}))\right\Vert \nonumber \\
&&+\left\Vert \sum\limits_{k'_1=1,k'_2=0}^{k'_1+k'_2=\ell+1} d_{k'_1,k'_2} T_{ f^{[1]}_{k'_1}(x_1,x_2)}^{\bm{X}_{1,\ell_1},\bm{X}_{2,\ell_2}}(T_{x_2^{k'_2}}^{\bm{X}_{1},\bm{X}_{2}}(\bm{Y})) -\sum\limits_{k'_1=1,k'_2=0}^{k'_1+k'_2=\ell+1} d_{k'_1,k'_2} T_{f^{[1]}_{k'_1}(x_1,x_2)}^{\bm{X}_{1},\bm{X}_{2}}(T_{x_2^{k'_2}}^{\bm{X}_{1},\bm{X}_{2}}(\bm{Y}))\right\Vert\nonumber \\
&\leq&\sum\limits_{k'_1=1,k'_2=0}^{k'_1+k'_2=\ell+1} |d_{k'_1,k'_2}|\left\Vert T_{ f^{[1]}_{k'_1}(x_1,x_2)}^{\bm{X}_{1,\ell_1},\bm{X}_{2,\ell_2}}(T_{x_2^{k'_2}}^{\bm{X}_{1,\ell_1},\bm{X}_{2,\ell_2}}(\bm{Y})) -T_{f^{[1]}_{k'_1}(x_1,x_2)}^{\bm{X}_{1,\ell_1},\bm{X}_{2,\ell_2}}(T_{x_2^{k'_2}}^{\bm{X}_{1},\bm{X}_{2}}(\bm{Y}))\right\Vert \nonumber \\
&&+\sum\limits_{k'_1=1,k'_2=0}^{k'_1+k'_2=\ell+1}|d_{k'_1,k'_2}|\left\Vert T_{ f^{[1]}_{k'_1}(x_1,x_2)}^{\bm{X}_{1,\ell_1},\bm{X}_{2,\ell_2}}(T_{x_2^{k'_2}}^{\bm{X}_{1},\bm{X}_{2}}(\bm{Y})) -T_{f^{[1]}_{k'_1}(x_1,x_2)}^{\bm{X}_{1},\bm{X}_{2}}(T_{x_2^{k'_2}}^{\bm{X}_{1},\bm{X}_{2}}(\bm{Y}))\right\Vert \nonumber \\
&\leq_2&\epsilon/6 + \epsilon/6 =\epsilon/3,
\end{eqnarray}
where we apply Lemma~\ref{lma: GDOI linear homomorphism} in $=_1$, apply (i) from Lemma~\ref{lma:continuity ii and iii} to obtain the first $\epsilon/6$ in $\leq_2$, and apply Theorem~\ref{thm:GDOI continuity dd} to obtain the second $\epsilon/6$ in $\leq_2$. 

Finally, this theorem is proved by combining Eq.~\eqref{eq5:thm:GDOI continuity}, Eq.~\eqref{eq6:thm:GDOI continuity}, and Eq.~\eqref{eq7:thm:GDOI continuity}. 
$\hfill\Box$

\section{Applications}\label{sec: Applications}

In this section, we discuss two applications of the proposed GDOI. The first application concerns the tail behavior of random matrices, as presented in Section~\ref{sec:Tail Bounds for Lipschitz Estimation for Random Matrices}. The second application extends the upper bound in Theorem~\ref{thm:Lipschitz Estimations} from \emph{Lipschitz estimations} to \emph{H\"older estimations}, as detailed in Section~\ref{sec:Matrix Holder Estimation}.

\subsection{Tail Bounds for Lipschitz Estimation for Random Matrices}\label{sec:Tail Bounds for Lipschitz Estimation for Random Matrices}

In the study of random matrix theory, particularly for Gaussian ensemble random matrices—including the Gaussian Orthogonal Ensemble (GOE), Gaussian Unitary Ensemble (GUE), and Gaussian Symplectic Ensemble (GSE)—the possibility of repeated (or degenerate) eigenvalues is a topic of significant theoretical interest. These ensembles are characterized by their eigenvalue distributions, which exhibit strong level repulsion due to the underlying statistical properties of these matrices.  

For Gaussian random matrices, the joint probability distribution of eigenvalues reveals a  repulsion phenomenon, meaning that eigenvalues tend to avoid clustering. Mathematically, in the limit of large matrices, the probability density function of eigenvalues contains a Vandermonde determinant squared (or raised to some power depending on the ensemble), which introduces a repelling force that suppresses the occurrence of degenerate eigenvalues. Specifically, in the GUE case, the probability of any two eigenvalues being exactly equal is zero, as the eigenvalues are governed by a repelling potential analogous to a  Coulomb gas model.  

A more formal argument arises from  perturbation theory: for an $n \times n $  matrix drawn from a Gaussian ensemble, if two or more eigenvalues were exactly identical, this would require the determinant of a highly structured polynomial system to vanish in a continuous probability space, which occurs with probability  zero. Thus, for finite-sized Gaussian random matrices, the probability of exact duplicate eigenvalues is zero in an idealized mathematical sense. However, in numerical computations, due to finite precision, near-duplicate eigenvalues may appear, though they do not constitute true degeneracies in the theoretical sense~\cite{edelman2005random}.  

In contrast, certain structured random matrix models—such as  Wishart matrices  or  non-Gaussian ensembles —may allow for degeneracies with nonzero probability. However, for GOE, GUE, and GSE, the probability of repeated eigenvalues remains strictly  zero  due to level repulsion, reinforcing the  universality  of eigenvalue distributions in these ensembles.  Given this, we consider alternative random specifications where the random matrix \( \bm{X} \) possesses  duplicate eigenvalues with nilpotent parts. For random matrices without a nilpotent component, such as  Hermitian matrices  (which can be viewed as a special case of tensors), we refer readers to our previous works~\cite{chang2022randomMOI,chang2022randomDTI,chang2023tailTRP,chang2023tailMulRanTenMeans,chang2023algebraicConn,chang2023BiRanTenPartI,chang2022generalMaj,chang2022randomPDT}.

The Jordan decomposition theorem states that any square matrix \(\bm{X} \in \mathbb{C}^{m \times m}\) can be decomposed as follows~\cite{gohberg1996simple}:  
\begin{equation}\label{eq: Jordan decomposition}
\bm{X} = \bm{U} \left( \bigoplus\limits_{k=1}^{K} \bigoplus\limits_{i=1}^{\alpha_k^{(\mathrm{G})}} \bm{J}_{m_{k,i}}(\lambda_k) \right) \bm{U}^{-1},
\end{equation}
where \(\bm{U} \in \mathbb{C}^{n \times n}\) is an invertible matrix, and \(\alpha_k^{(\mathrm{G})}\) represents the geometric multiplicity corresponding to the \(k\)-th eigenvalue \(\lambda_k\). Given $K,\alpha_k^{(\mathrm{G})}, m_{k,i}$, the randomness of the random matrix $\bm{X}$ comes from $K$ random eigenvalues $\lambda_k$ and random invertible matrix $\bm{U}$. For those random matrices with eigenvalues $\lambda_k \neq 0$, they are \emph{analogous} each other if they share same $K,\alpha_k^{(\mathrm{G})}$ and $m_{k,i}$~\cite{chang2024operatorChar}. In this paper, all random matrices are assumed to follow such randomness specifications.

\begin{theorem}\label{thm:Lipschitz Estimations RM}
Given an analytic function $f(z)$ within the domain for $|z| < R$, the first rndom matrix $\bm{X}_1$ with the dimension $m$ and $K_1$ distinct eigenvalues $\lambda_{k_1}$ for $k_1=1,2,\ldots,K_1$ such that
\begin{eqnarray}\label{eq1-1:thm:Lipschitz Estimations RM}
\bm{X}_1&=&\sum\limits_{k_1=1}^{K_1}\sum\limits_{i_1=1}^{\alpha_{k_1}^{\mathrm{G}}} \lambda_{k_1} \bm{P}_{k_1,i_1}+
\sum\limits_{k_1=1}^{K_1}\sum\limits_{i_1=1}^{\alpha_{k_1}^{\mathrm{G}}} \bm{N}_{k_1,i_1},
\end{eqnarray}
where $\left\vert\lambda_{k_1}\right\vert<R$, and second random matrix $\bm{X}_2$ with the dimension $m$ and $K_2$ distinct eigenvalues $\lambda_{k_2}$ for $k_2=1,2,\ldots,K_2$ such that
\begin{eqnarray}\label{eq1-2:thm:Lipschitz Estimations RM}
\bm{X}_2&=&\sum\limits_{k_2=1}^{K_2}\sum\limits_{i_2=1}^{\alpha_{k_2}^{\mathrm{G}}} \lambda_{k_2} \bm{P}_{k_2,i_2}+
\sum\limits_{k_2=1}^{K_2}\sum\limits_{i_2=1}^{\alpha_{k_2}^{\mathrm{G}}} \bm{N}_{k_2,i_2},
\end{eqnarray}
where $\left\vert\lambda_{k_2}\right\vert<R$. We also assume that $\lambda_{k_1} \neq \lambda_{k_2}$ for any $k_1$ and $k_2$. The random matrices $\bm{X}_1$ and $\bm{X}_2$ are independent each other. We define $f^{[1]}(\lambda_1, \lambda_2) \define \frac{f(\lambda_1) - f(\lambda_2)}{\lambda_1 - \lambda_2}$. 

We also assume that 
\begin{eqnarray}\label{eq1-3:thm:Lipschitz Estimations RM}
\max \left\Vert\bm{N}_{k_1,i_1}^{q_1}\right\Vert&\leq&\Gamma_{k_1,i_1}^{q_1},\nonumber \\
\max \left\Vert\bm{N}_{k_2,i_2}^{q_2}\right\Vert&\leq&\Gamma_{k_2,i_2}^{q_2}.
\end{eqnarray}

Then, we have the following tail bound for the random variable $\left\Vert f(\bm{X}_1) - f(\bm{X}_2) \right\Vert$, which is
\begin{eqnarray}\label{eq2:thm:Lipschitz Estimations RM}
\mathrm{Pr}\left(\left\Vert f(\bm{X}_1) - f(\bm{X}_2) \right\Vert\geq\delta\right)\leq(B_1 + B_2 + B_3 + B_4)\frac{\mathbb{E}\left[\left\Vert\bm{X}_1-\bm{X}_2\right\Vert\right]}{\delta},
\end{eqnarray}
where $\delta$ is any positive number, $\mathbb{E}$ is the expectation operation, and terms $B_1, B_2, B_3, B_4$ are expressed by
\begin{eqnarray}\label{eq2-1:thm:Lipschitz Estimations RM}
B_1&=& \left[\max\limits_{\lambda_1 \in \Lambda_{\bm{X}_1},\lambda_2 \in \Lambda_{\bm{X}_2}} \left\vert f^{[1]}(\lambda_1,\lambda_2)\right\vert\right], \nonumber \\
B_2&=& \sum\limits_{k_2=1}^{K_2}\sum\limits_{i_2=1}^{\alpha_{k_2}^{(\mathrm{G})}}\sum_{q_2=1}^{m_{k_2,i_2}-1}\left[\max\limits_{\lambda_1 \in \Lambda_{\bm{X}_1},\lambda_2 \in \Lambda_{\bm{X}_2}}\left\vert\frac{(f^{[1]})^{(-,q_2)}(\lambda_{1},\lambda_{2})}{q_2!}\right\vert\right]\Gamma_{k_2,i_2}^{q_2}, \nonumber \\
B_3&=&\sum\limits_{k_1=1}^{K_1}\sum\limits_{i_1=1}^{\alpha_{k_1}^{(\mathrm{G})}}\sum_{q_1=1}^{m_{k_1,i_1}-1}\left[\max\limits_{\lambda_1 \in \Lambda_{\bm{X}_1},\lambda_2 \in \Lambda_{\bm{X}_2}}\left\vert\frac{(f^{[1]})^{(q_1,-)}(\lambda_{1},\lambda_{2})}{q_1!}\right\vert\right]\Gamma_{k_1,i_1}^{q_1}, \nonumber \\
B_4&=& \sum\limits_{k_1=1}^{K_1}\sum\limits_{k_2=1}^{K_2}\sum\limits_{i_1=1}^{\alpha_{k_1}^{(\mathrm{G})}}\sum\limits_{i_2=1}^{\alpha_{k_2}^{(\mathrm{G})}}\sum_{q_1=1}^{m_{k_1,i_1}-1}\sum_{q_2=1}^{m_{k_2,i_2}-1}\left[\max\limits_{\lambda_1 \in \Lambda_{\bm{X}_1},\lambda_2 \in \Lambda_{\bm{X}_2}}\left\vert\frac{(f^{[1]})^{(q_1,q_2)}(\lambda_{1},\lambda_{2})}{q_1!q_2!}\right\vert\right]\Gamma_{k_1,i_1}^{q_1}\Gamma_{k_2,i_2}^{q_2}.
\end{eqnarray}
\end{theorem}
\textbf{Proof:}
From Markov inequality, we have
\begin{eqnarray}\label{eq4:thm:Lipschitz Estimations RM}
\mathrm{Pr}\left(\left\Vert f(\bm{X}_1) - f(\bm{X}_2) \right\Vert\geq\delta\right)\leq\frac{\mathbb{E}\left[\left\Vert f(\bm{X}_1) - f(\bm{X}_2) \right\Vert\right]}{\delta}.
\end{eqnarray}

From Theorem~\ref{thm:Lipschitz Estimations}, we have the following bound for the term $\mathbb{E}\left[\left\Vert f(\bm{X}_1) - f(\bm{X}_2) \right\Vert\right]$:
\begin{eqnarray}\label{eq5:thm:Lipschitz Estimations RM}
\lefteqn{\mathbb{E}\left[\left\Vert f(\bm{X}_1) - f(\bm{X}_2) \right\Vert\right]}\nonumber \\
&\leq& \left[\max\limits_{\lambda_1 \in \Lambda_{\bm{X}_1},\lambda_2 \in \Lambda_{\bm{X}_2}} \left\vert f^{[1]}(\lambda_1,\lambda_2)\right\vert\right]\mathbb{E}\left[\left\Vert\bm{X}_1-\bm{X}_2\right\Vert\right] \nonumber \\
&&+\sum\limits_{k_2=1}^{K_2}\sum\limits_{i_2=1}^{\alpha_{k_2}^{(\mathrm{G})}}\sum_{q_2=1}^{m_{k_2,i_2}-1}\left[\max\limits_{\lambda_1 \in \Lambda_{\bm{X}_1},\lambda_2 \in \Lambda_{\bm{X}_2}}\left\vert\frac{(f^{[1]})^{(-,q_2)}(\lambda_{1},\lambda_{2})}{q_2!}\right\vert\right]\mathbb{E}\left[\left\Vert\bm{N}_{k_2,i_2}^{q_2}\right\Vert\left\Vert\bm{X}_1-\bm{X}_2\right\Vert\right]\nonumber \\
&&+\sum\limits_{k_1=1}^{K_1}\sum\limits_{i_1=1}^{\alpha_{k_1}^{(\mathrm{G})}}\sum_{q_1=1}^{m_{k_1,i_1}-1}\left[\max\limits_{\lambda_1 \in \Lambda_{\bm{X}_1},\lambda_2 \in \Lambda_{\bm{X}_2}}\left\vert\frac{(f^{[1]})^{(q_1,-)}(\lambda_{1},\lambda_{2})}{q_1!}\right\vert\right]\mathbb{E}\left[\left\Vert\bm{N}_{k_1,i_1}^{q_1}\right\Vert\left\Vert\bm{X}_1-\bm{X}_2\right\Vert\right]\nonumber \\
&&+ \sum\limits_{k_1=1}^{K_1}\sum\limits_{k_2=1}^{K_2}\sum\limits_{i_1=1}^{\alpha_{k_1}^{(\mathrm{G})}}\sum\limits_{i_2=1}^{\alpha_{k_2}^{(\mathrm{G})}}\sum_{q_1=1}^{m_{k_1,i_1}-1}\sum_{q_2=1}^{m_{k_2,i_2}-1}\left[\max\limits_{\lambda_1 \in \Lambda_{\bm{X}_1},\lambda_2 \in \Lambda_{\bm{X}_2}}\left\vert\frac{(f^{[1]})^{(q_1,q_2)}(\lambda_{1},\lambda_{2})}{q_1!q_2!}\right\vert\right]\nonumber\\
&&~~~\times\mathbb{E}\left[\left\Vert\bm{N}_{k_1,i_1}^{q_1}\right\Vert\left\Vert\bm{N}_{k_2,i_2}^{q_2}\right\Vert\left\Vert\bm{X}_1-\bm{X}_2\right\Vert\right]\nonumber \\
&\leq_1& \left[\max\limits_{\lambda_1 \in \Lambda_{\bm{X}_1},\lambda_2 \in \Lambda_{\bm{X}_2}} \left\vert f^{[1]}(\lambda_1,\lambda_2)\right\vert\right]\mathbb{E}\left[\left\Vert\bm{X}_1-\bm{X}_2\right\Vert\right] \nonumber \\
&&+\sum\limits_{k_2=1}^{K_2}\sum\limits_{i_2=1}^{\alpha_{k_2}^{(\mathrm{G})}}\sum_{q_2=1}^{m_{k_2,i_2}-1}\left[\max\limits_{\lambda_1 \in \Lambda_{\bm{X}_1},\lambda_2 \in \Lambda_{\bm{X}_2}}\left\vert\frac{(f^{[1]})^{(-,q_2)}(\lambda_{1},\lambda_{2})}{q_2!}\right\vert\right]\Gamma_{k_2,i_2}^{q_2}\mathbb{E}\left[\left\Vert\bm{X}_1-\bm{X}_2\right\Vert\right]\nonumber \\
&&+\sum\limits_{k_1=1}^{K_1}\sum\limits_{i_1=1}^{\alpha_{k_1}^{(\mathrm{G})}}\sum_{q_1=1}^{m_{k_1,i_1}-1}\left[\max\limits_{\lambda_1 \in \Lambda_{\bm{X}_1},\lambda_2 \in \Lambda_{\bm{X}_2}}\left\vert\frac{(f^{[1]})^{(q_1,-)}(\lambda_{1},\lambda_{2})}{q_1!}\right\vert\right]\Gamma_{k_1,i_1}^{q_1}\mathbb{E}\left[\left\Vert\bm{X}_1-\bm{X}_2\right\Vert\right]\nonumber \\
&&+ \sum\limits_{k_1=1}^{K_1}\sum\limits_{k_2=1}^{K_2}\sum\limits_{i_1=1}^{\alpha_{k_1}^{(\mathrm{G})}}\sum\limits_{i_2=1}^{\alpha_{k_2}^{(\mathrm{G})}}\sum_{q_1=1}^{m_{k_1,i_1}-1}\sum_{q_2=1}^{m_{k_2,i_2}-1}\left[\max\limits_{\lambda_1 \in \Lambda_{\bm{X}_1},\lambda_2 \in \Lambda_{\bm{X}_2}}\left\vert\frac{(f^{[1]})^{(q_1,q_2)}(\lambda_{1},\lambda_{2})}{q_1!q_2!}\right\vert\right]\nonumber\\
&&~~~\times\Gamma_{k_1,i_1}^{q_1}\Gamma_{k_2,i_2}^{q_2}\mathbb{E}\left[\left\Vert\bm{X}_1-\bm{X}_2\right\Vert\right],
\end{eqnarray}
where we apply the assumption provided by Eq.~\eqref{eq1-3:thm:Lipschitz Estimations RM} to $\leq_1$

Finally, this theorem is proved by applying Eq.~\eqref{eq5:thm:Lipschitz Estimations RM} to Eq.~\eqref{eq4:thm:Lipschitz Estimations RM}.
$\hfill\Box$

\subsection{Matrix H\"older Estimation}\label{sec:Matrix Holder Estimation}

The purpose of this section is to extend Theorem~\ref{thm:Lipschitz Estimations} upper bound from Lipschitz estimations to H\"older Estimation. Recall the set of H\"older function is defined as
\begin{eqnarray}\label{eq:Holder func def}
S_{\mbox{H\"older}}&\define&\{f: \mathbb{C} \rightarrow \mathbb{C} \mbox{~~such that}\sup\limits_{x,y \in \mathbb{R}, x \neq y}\frac{\left\vert f(x)-f(y) \right\vert}{\left\vert x-y \right\vert^{\omega}} < \infty\}
\end{eqnarray}
where $\omega$ is a positive number. 

Given two matrices $\bm{X}_1$ and $\bm{X}_2$ with spectrum in a bounded domain $\Lambda$, we have the following Lemma~\ref{lma:matrix diff bound by omega power} and Lemma~\ref{lma:Lip est bound by Holder est}. These two lemmas are used to quantify the effects of $\omega$ with respect to norms. 
\begin{lemma}\label{lma:matrix diff bound by omega power}
Given two matrices $\bm{X}_1$ and $\bm{X}_2$ with spectrum in a bounded domain $\Lambda$,  we have
\begin{eqnarray}\label{eq1:lma:matrix diff bound by omega power}
\left\Vert \bm{X}_1 - \bm{X}_2 \right\Vert \leq C_{\Lambda,\omega}\left\Vert \bm{X}_1 - \bm{X}_2 \right\Vert^{\omega}
\end{eqnarray}
where $C_{\Lambda,\omega}$ is a constant depending only on the bounded domain $\Lambda$ and a positive real number $\omega$. We assume that $\left\Vert \bm{X}_1 - \bm{X}_2 \right\Vert \geq \nu$ for some positive number $\nu$.
\end{lemma}
\textbf{Proof:}
Let \( t = \|\bm{X}_1 - \bm{X}_2\| \), and we aim to prove the inequality  

   \[
   t \leq C t^\omega
   \]

   for some constant \( C \) depending only on the bounded domain \( \Lambda \) and any positive real number \( \omega \).  This inequality can be rewritten as     

   \[
   t^{1-\omega} \leq C.
   \]

Since we assume that \( t \) has a lower bound \( \nu > 0 \), we have  $t \geq \nu$.  Consider the function \( f(t) = t^{1-\omega} \). The behavior of this function depends on \( \omega \):  
\begin{enumerate}
\item If \( \omega = 1 \), then \( f(t) = 1 \), and we can choose \( C = 1 \), making the inequality trivially hold.  
\item If \( \omega > 1 \), then \( f(t) = t^{1-\omega} \) is decreasing in \( t \). Since \( t \geq \nu \), the maximum value of \( f(t) \) in this range is attained at \( t = \nu \), giving  

     \[
     t^{1-\omega} \leq \nu^{1-\omega}.
     \]

     Thus, we can take \( C = \nu^{1-\omega} \), ensuring the inequality holds.  
\item If \( 0 < \omega < 1 \), then \( f(t) = t^{1-\omega} \) is increasing in \( t \). Since \( t \geq \nu \), the maximum occurs at the largest possible \( t \), which is bounded due to the spectral constraints of \( \bm{X}_1 \) and \( \bm{X}_2 \). Suppose there exists an upper bound \( M \), then  

     \[
     t^{1-\omega} \leq M^{1-\omega}.
     \]

     In this case, we can choose \( C = M^{1-\omega} \).  
\end{enumerate}

In all cases, there exists a constant \( C \) that depends only on \( \Lambda \), ensuring that the inequality  

   \[
   t \leq C t^\omega
   \]

   holds for all \( t \geq \nu \).  By setting \( C_{\Lambda,\omega}=C\) derived above from each case, we can choose $C_{\Lambda,\omega}$ to satisfy the given inequality. Thus, the proof is complete. 
$\hfill\Box$

\begin{lemma}\label{lma:Lip est bound by Holder est}
We define Lipschitz seminorm with respect to a function $f$ as
\begin{eqnarray}\label{eq1:lma:Lip est bound by Holder est}
D_1(f)\define \sup\limits_{x,y \in \Lambda, x \neq y}\frac{\left\vert f(x)-f(y) \right\vert}{\left\vert x-y \right\vert},
\end{eqnarray}
and H\"older seminorm as
\begin{eqnarray}\label{eq2:lma:Lip est bound by Holder est}
D_\omega(f)\define \sup\limits_{x,y \in \Lambda, x \neq y}\frac{\left\vert f(x)-f(y) \right\vert}{\left\vert x-y \right\vert^{\omega}}
\end{eqnarray}
where $\omega$ is a positive number. 

Then, we have
\begin{eqnarray}\label{eq3:lma:Lip est bound by Holder est}
D_1(f) \leq C'_{\Lambda,\omega} D_\omega(f),
\end{eqnarray}
where $C'_{\Lambda,\omega}$ is a constant depending only on the bounded domain $\Lambda$. We assume that $\left\vert x-y \right\vert \geq \nu'$ for some positive number $\nu'$.
\end{lemma}
\textbf{Proof:}
By the definition of the H\"older seminorm \( D_\omega(f) \), we have  
\begin{eqnarray}
\left\vert f(x)-f(y) \right\vert \leq D_\omega(f) \left\vert x-y \right\vert^\omega, \quad \forall x, y \in \Lambda, x \neq y.
\end{eqnarray}
Dividing both sides by \( \left\vert x - y \right\vert \) (which is nonzero by assumption), we obtain  
\begin{eqnarray}
\frac{\left\vert f(x)-f(y) \right\vert}{\left\vert x-y \right\vert} \leq D_\omega(f) \left\vert x-y \right\vert^{\omega - 1}.
\end{eqnarray}
Taking the supremum over all \( x, y \in \Lambda \), we get  
\begin{eqnarray}
D_1(f) = \sup\limits_{x,y \in \Lambda, x \neq y} \frac{\left\vert f(x)-f(y) \right\vert}{\left\vert x-y \right\vert}
\leq D_\omega(f) \sup\limits_{x,y \in \Lambda, x \neq y} \left\vert x - y \right\vert^{\omega - 1}.
\end{eqnarray}
Since we assume \( \left\vert x - y \right\vert \geq \nu' \) for some positive number \( \nu' \), it follows that  
\begin{eqnarray}
\sup\limits_{x,y \in \Lambda, x \neq y} \left\vert x - y \right\vert^{\omega - 1} \leq (\nu')^{\omega - 1},
\end{eqnarray}
if $\omega < 1$. On the other hand, if $\omega \geq 1$, we have
\begin{eqnarray}
\sup\limits_{x,y \in \Lambda, x \neq y} \left\vert x - y \right\vert^{\omega - 1} \leq |\Lambda|^{\omega - 1},
\end{eqnarray}
where $|\Lambda|$ is the range norm of the domain $\Lambda$.

Thus, we can set \( C'_{\Lambda, \omega} = (\nu')^{\omega - 1} \) if $\omega < 1$ or \( C'_{\Lambda, \omega} = |\Lambda|^{\omega - 1} \) if $\omega \geq 1$, leading to  
\begin{eqnarray}
D_1(f) \leq C'_{\Lambda,\omega} D_\omega(f),
\end{eqnarray}
which proves the desired inequality.
$\hfill\Box$

We are ready to present the H\"older estimations for the proposed GDOI.

\begin{theorem}\label{thm:Lipschitz Estimations HO}
Given an analytic function $f(z)$ within the domain for $|z| < R$, the first matrix $\bm{X}_1$ with the dimension $m$ and $K_1$ distinct eigenvalues $\lambda_{k_1}$ for $k_1=1,2,\ldots,K_1$ such that
\begin{eqnarray}\label{eq1-1:thm:Lipschitz Estimations HO}
\bm{X}_1&=&\sum\limits_{k_1=1}^{K_1}\sum\limits_{i_1=1}^{\alpha_{k_1}^{\mathrm{G}}} \lambda_{k_1} \bm{P}_{k_1,i_1}+
\sum\limits_{k_1=1}^{K_1}\sum\limits_{i_1=1}^{\alpha_{k_1}^{\mathrm{G}}} \bm{N}_{k_1,i_1},
\end{eqnarray}
where $\left\vert\lambda_{k_1}\right\vert<R$, and second matrix $\bm{X}_2$ with the dimension $m$ and $K_2$ distinct eigenvalues $\lambda_{k_2}$ for $k_2=1,2,\ldots,K_2$ such that
\begin{eqnarray}\label{eq1-2:thm:Lipschitz Estimations HO}
\bm{X}_2&=&\sum\limits_{k_2=1}^{K_2}\sum\limits_{i_2=1}^{\alpha_{k_2}^{\mathrm{G}}} \lambda_{k_2} \bm{P}_{k_2,i_2}+
\sum\limits_{k_2=1}^{K_2}\sum\limits_{i_2=1}^{\alpha_{k_2}^{\mathrm{G}}} \bm{N}_{k_2,i_2},
\end{eqnarray}
where $\left\vert\lambda_{k_2}\right\vert<R$. We also assume that $\lambda_{k_1} \neq \lambda_{k_2}$ for any $k_1$ and $k_2$. We define $f^{[1]}(\lambda_1, \lambda_2) \define \frac{f(\lambda_1) - f(\lambda_2)}{\lambda_1 - \lambda_2}$. 

We also assume that 
\begin{eqnarray}\label{eq1-3:thm:Lipschitz Estimations HO}
\max \left\Vert\bm{N}_{k_1,i_1}^{q_1}\right\Vert&\leq&\Gamma_{k_1,i_1}^{q_1},\nonumber \\
\max \left\Vert\bm{N}_{k_2,i_2}^{q_2}\right\Vert&\leq&\Gamma_{k_2,i_2}^{q_2}.
\end{eqnarray}

Then, we have the following upper bound for H\"older estimation:
\begin{eqnarray}\label{eq2:thm:Lipschitz Estimations HO}
\lefteqn{\left\Vert f(\bm{X}_1) - f(\bm{X}_2) \right\Vert \leq}\nonumber \\
&& C_{\Lambda, \omega}C'_{\Lambda, \omega}D_\omega(f)\left\Vert\bm{X}_1-\bm{X}_2\right\Vert ^\omega\nonumber \\
&&+ C_{\Lambda, \omega}\sum\limits_{k_2=1}^{K_2}\sum\limits_{i_2=1}^{\alpha_{k_2}^{(\mathrm{G})}}\sum_{q_2=1}^{m_{k_2,i_2}-1}\left[\max\limits_{\lambda_1 \in \Lambda_{\bm{X}_1},\lambda_2 \in \Lambda_{\bm{X}_2}}\left\vert\frac{(f^{[1]})^{(-,q_2)}(\lambda_{1},\lambda_{2})}{q_2!}\right\vert\right]\Gamma_{k_2,i_2}^{q_2}\left\Vert\bm{X}_1-\bm{X}_2\right\Vert^\omega\nonumber \\
&&+ C_{\Lambda, \omega}\sum\limits_{k_1=1}^{K_1}\sum\limits_{i_1=1}^{\alpha_{k_1}^{(\mathrm{G})}}\sum_{q_1=1}^{m_{k_1,i_1}-1}\left[\max\limits_{\lambda_1 \in \Lambda_{\bm{X}_1},\lambda_2 \in \Lambda_{\bm{X}_2}}\left\vert\frac{(f^{[1]})^{(q_1,-)}(\lambda_{1},\lambda_{2})}{q_1!}\right\vert\right]\Gamma_{k_1,i_1}^{q_1}\left\Vert\bm{X}_1-\bm{X}_2\right\Vert^\omega \nonumber \\
&&+  C_{\Lambda, \omega}\sum\limits_{k_1=1}^{K_1}\sum\limits_{k_2=1}^{K_2}\sum\limits_{i_1=1}^{\alpha_{k_1}^{(\mathrm{G})}}\sum\limits_{i_2=1}^{\alpha_{k_2}^{(\mathrm{G})}}\sum_{q_1=1}^{m_{k_1,i_1}-1}\sum_{q_2=1}^{m_{k_2,i_2}-1}\left[\max\limits_{\lambda_1 \in \Lambda_{\bm{X}_1},\lambda_2 \in \Lambda_{\bm{X}_2}}\left\vert\frac{(f^{[1]})^{(q_1,q_2)}(\lambda_{1},\lambda_{2})}{q_1!q_2!}\right\vert\right]\nonumber\\
&&~~~\times\Gamma_{k_1,i_1}^{q_1}\Gamma_{k_2,i_2}^{q_2}\left\Vert\bm{X}_1-\bm{X}_2\right\Vert^\omega.
\end{eqnarray}
\end{theorem}
\textbf{Proof:}
From Theorem~\ref{thm:Lipschitz Estimations}, we have the following upper bound for Lipschitz estimation:
\begin{eqnarray}\label{eq3:thm:Lipschitz Estimations HO}
\lefteqn{\left\Vert f(\bm{X}_1) - f(\bm{X}_2) \right\Vert \leq}\nonumber \\
&& \underbrace{\left[\max\limits_{\lambda_1 \in \Lambda_{\bm{X}_1},\lambda_2 \in \Lambda_{\bm{X}_2}} \left\vert f^{[1]}(\lambda_1,\lambda_2)\right\vert\right]\left\Vert\bm{X}_1-\bm{X}_2\right\Vert}_{\mbox{Part I}}\nonumber \\
&&+\underbrace{\sum\limits_{k_2=1}^{K_2}\sum\limits_{i_2=1}^{\alpha_{k_2}^{(\mathrm{G})}}\sum_{q_2=1}^{m_{k_2,i_2}-1}\left[\max\limits_{\lambda_1 \in \Lambda_{\bm{X}_1},\lambda_2 \in \Lambda_{\bm{X}_2}}\left\vert\frac{(f^{[1]})^{(-,q_2)}(\lambda_{1},\lambda_{2})}{q_2!}\right\vert\right]\left\Vert\bm{N}_{k_2,i_2}^{q_2}\right\Vert\left\Vert\bm{X}_1-\bm{X}_2\right\Vert}_{\mbox{Part II}}\nonumber \\
&&+\underbrace{\sum\limits_{k_1=1}^{K_1}\sum\limits_{i_1=1}^{\alpha_{k_1}^{(\mathrm{G})}}\sum_{q_1=1}^{m_{k_1,i_1}-1}\left[\max\limits_{\lambda_1 \in \Lambda_{\bm{X}_1},\lambda_2 \in \Lambda_{\bm{X}_2}}\left\vert\frac{(f^{[1]})^{(q_1,-)}(\lambda_{1},\lambda_{2})}{q_1!}\right\vert\right]\left\Vert\bm{N}_{k_1,i_1}^{q_1}\right\Vert\left\Vert\bm{X}_1-\bm{X}_2\right\Vert}_{\mbox{Part III}}  \nonumber \\
&&+\sum\limits_{k_1=1}^{K_1}\sum\limits_{k_2=1}^{K_2}\sum\limits_{i_1=1}^{\alpha_{k_1}^{(\mathrm{G})}}\sum\limits_{i_2=1}^{\alpha_{k_2}^{(\mathrm{G})}}\sum_{q_1=1}^{m_{k_1,i_1}-1}\sum_{q_2=1}^{m_{k_2,i_2}-1}\left[\max\limits_{\lambda_1 \in \Lambda_{\bm{X}_1},\lambda_2 \in \Lambda_{\bm{X}_2}}\left\vert\frac{(f^{[1]})^{(q_1,q_2)}(\lambda_{1},\lambda_{2})}{q_1!q_2!}\right\vert\right]\nonumber\\
&&~~~ \underbrace{\times\left\Vert\bm{N}_{k_1,i_1}^{q_1}\right\Vert\left\Vert\bm{N}_{k_2,i_2}^{q_2}\right\Vert\left\Vert\bm{X}_1-\bm{X}_2\right\Vert}_{\mbox{Part IV}}.
\end{eqnarray} 

This theorem is proved by applying the following 
\begin{itemize}
\item For Part I, we use Lemma~\ref{lma:matrix diff bound by omega power} and Lemma~\ref{lma:Lip est bound by Holder est}.
\item For Part II, Part III and Part IV, we use the assumption given by Eq.~\eqref{eq1-3:thm:Lipschitz Estimations HO} and Lemma~\ref{lma:matrix diff bound by omega power}.
\end{itemize}
$\hfill\Box$

\bibliographystyle{IEEETran}
\bibliography{SpecialCase_and_DOI_Bib}

\end{document}